\documentclass[smallextended,nospthms,
natbib]{svjour3}

\usepackage{graphicx}
\usepackage{amsmath}
\usepackage{amsfonts}
\usepackage{amssymb}
\usepackage{amsthm}
\usepackage{mathtools}
\usepackage{pdfsync}

\usepackage{times}
\usepackage{mathptmx}
\usepackage{epic,graphpap,graphics,float,tabularx,comment}

\usepackage{bm}
\usepackage[colorlinks,citecolor=blue]{hyperref}

\numberwithin{equation}{section}

\newcommand{\bX}{\mathbf{X}}
\newcommand{\bU}{\mathbf{U}}
\newcommand{\be}{\mathbf{e}}
\newcommand{\q}[1]{\langle{#1}\rangle} %
\renewcommand{\d}{\mathrm{d}}  
\newcommand{\D}[1][x_1]{\frac{\partial}{\partial #1}}
\newcommand{\DD}[1][x_1]{\frac{\partial^2}{\partial #1^2}}
\newcommand{\cL}{\mathcal{L}}
\newcommand\indep{\protect\mathpalette{\protect\independenT}{\perp}}
\def\independenT#1#2{\mathrel{\rlap{$#1#2$}\mkern2mu{#1#2}}}

\newcommand{\p}{\mathbb{P}}            %
\newcommand{\R}{\mathbb{R}}            %
\newcommand{\N}{\mathbb{N}}            %
\newcommand{\cA}{{\mathcal A}}

\newcommand{\F}{{\cal F}}         %

\newcommand{\ra}{\rightarrow}           %

\newcommand{\Om}{\Omega}

\DeclareMathOperator{\sign}{sgn}

\bmdefine\mub{\mu}
\bmdefine\etab{\eta}
\bmdefine\thetab{\vartheta}
\bmdefine\betab{\beta}
\bmdefine\sigmab{\sigma}
\bmdefine\varsigmab{\varsigma}
\bmdefine\taub{\tau}
\bmdefine\nub{\nu}
\bmdefine\rhob{\rho}
\bmdefine\varrhob{\varrho}
\bmdefine\gammab{\gamma}
\bmdefine\Gammab{\Gamma}

\makeatletter
\newcommand{\zjel}[1]{\left({#1}\right)}
\newcommand{\br}[1]{\left[{#1}\right]}
\newcommand{\cbr}[1]{\left\{{#1}\right\}}
\newcommand{\abs}[1]{\left|{#1}\right|}
\newcommand{\norm}[1]{\left\|{#1}\right\|}
\newcommand{\set}[2]{\cbr{#1\,:\,#2}}

\def\isparamm{\@ifnextchar{\bgroup}}
\def\subscr#1{_{\{#1\}}}
\def\subnozjel#1{_{#1}}
\newcommand{\I}[1][]{\mathbf{1}%
  \@ifempty{#1}{\isparamm\subscr\subnozjel}{\subnozjel{#1}}}
\newcommand{\E}[1][]{\mathbb{E}\@ifempty{#1}{}{_{#1}}\isparamm\zjel\relax}
\renewcommand{\P}[1][]{\mathbb{P}\@ifempty{#1}{}{_{#1}}\isparamm\zjel\relax}
\newcommand{\PQ}[1][]{\mathbb{Q}\@ifempty{#1}{}{_{#1}}\isparamm\zjel\relax}
\makeatother

\newcounter{enumr}
\newenvironment{rlist}[1][]{%
  \begin{list}{%
      \hbox to 0pt{%
        \hss(\theenumr)%
      }\unskip\ignorespaces
    }{%
      \usecounter{enumr}%
    }%
    \def\theenumr{\roman{enumr}#1}%
  }{%
  \end{list}%
}

\theoremstyle{plain}
\newtheorem{theorem}{Theorem}[section]
\newtheorem{proposition}{Proposition}[section]

\theoremstyle{definition}
\newtheorem{example}{Example}[section]

\theoremstyle{remark}
\newtheorem{remark}{Remark}[section]

\title{Planar Diffusions with Rank-Based Characteristics: Transition Probabilities, Time Reversal, Maximality and Perturbed Tanaka equations} 

\titlerunning{Planar Diffusions with Rank-Based Characteristics}

\author{E. Robert Fernholz \and Tomoyuki Ichiba \and  Ioannis Karatzas \and
  Vilmos Prokaj}
\institute{E. Robert Fernholz \at \textsc{Intech} Investment Management LLC, One Palmer Square, Suite
  441,  Princeton, NJ 08542 (\email{bob@enhanced.com}) \and
  Tomoyuki Ichiba \at
  Department of Statistics and Applied Probability, South Hall,
  University of California, Santa Barbara, CA 93106 (\email{ichiba@pstat.ucsb.edu}) 
  \and
  Ioannis Karatzas \at
  \textsc{Intech} Investment Management LLC, One Palmer Square, Suite
  441,  Princeton, NJ 08542 (\email{ik@enhanced.com}), and
  Department of Mathematics, Columbia University, New York, NY 10027 (\email{ik@math.columbia.edu}). 
  \and
  Vilmos Prokaj \at
  Department of Probability Theory and Statistics, E\"otv\"os
  Lor\'and University, 1117 Budapest, P\'azm\'any P\'eter s\'et\'any 1/C,
  Hungary, and Computer and Automation Institute of the Hungarian Academy of
  Sciences, 1111 Budapest, Kende utca 13-17,  Hungary (\email{prokaj@cs.elte.hu}). 
}  
\date{}

\begin{document}

\maketitle

\begin{abstract}
  For   given nonnegative constants $\, g  \,$, $\,
  h  \,$, $\, \rho  \,$, $\, \sigma \,$ with $\, \rho^2 + \sigma^2 =1\,$ and
  $\, g+h>0\,$, we construct a  diffusion process $\, (X_1 (\cdot), X_2
  (\cdot))\,$ with values in the   plane and  infinitesimal generator 
  \begin{multline}
    \label{0.1}
    \cL=\I{x_1 > x_2}\zjel{\frac{\rho^2}2 \DD+ \frac{\sigma^2}2
      \DD[x_2]-h\D+g\D[x_2]}+\\
    \I{x_1 \le x_2}\zjel{\frac{\sigma^2}2 \DD+ \frac{\rho^2}2
      \DD[x_2]+g\D-h\D[x_2]}.
  \end{multline}
  We  compute the  transition probabilities of this process, discuss its
  realization in terms of appropriate systems of stochastic differential
  equations,   study its  dynamics   under a  time reversal, and note that
  these     involve singularly continuous components governed by local
  time. Crucial in our analysis are  properties of Brownian and
  semimartingale local time; properties of the   {\it generalized perturbed Tanaka     equation} 
  $$ 
  \d  Z (t)\,=\,
  f \big(Z (t)\big) \, \d M( t) + \d  N(t) \,, ~ ~~~~Z(0)=\xi \,
  $$
  driven by suitable continuous, orthogonal semimartingales $M(\cdot)\,$ and
  $\,  N(\cdot)\,$ and with $\,f(\cdot)\,$ of bounded variation, which we
  study here in detail;   and  those of a one-di\-men\-sional diffusion   $\, 
  Y(\cdot)\,$  with   {\it bang-bang} drift     $    \,\d  Y (t)=- \lambda  \,
  \sign \big(  Y (t) \big) \, \d  t $ $+ \, \d  W  (t)$,  $Y(0)=y\,\,$
  driven by a standard Brownian motion $\,  W(\cdot)\,$.   
  
 We also show that the planar diffusion $\, (X_1 (\cdot), X_2 (\cdot))\,$  
  can be represented  in terms of this process $\, Y(\cdot)\,$, its local time $\,
  L^Y (\cdot)\,$ at the origin, and an independent standard Brownian motion
  $\, Q(\cdot)\,$,  in a form  which can be construed as a two-dimensional
  analogue of the stochastic equation satisfied by the so-called {\it skew
    Brownian motion.}      
 
\end{abstract}
\section{ \textsc{Introduction} }

For   given nonnegative constants $\, g  \,$, $\, h  \,$, $\, \rho  \,$, $\,
\sigma \,$ with $\, \rho^2 + \sigma^2 >0\,$ and $\, g+h>0\,$,  and for a given
vector $\, (x_1, x_2 ) \in \R^2\,$, we shall consider the question of
constructing a two-dimensional diffusion process $\, (X_1 (\cdot), X_2
(\cdot))\,$ with dynamics 
\begin{align}
  \label{1.1}
  \d  X_1 (t) &=   
  \left(g\I{X_1 (t) \le  X_2 (t)} 
    - h\I{ X_1 (t) >  X_2 (t)}\right)  \d  t + 
  \left(\rho\I{X_1(t)>X_2 (t)}+\sigma\I{X_1(t)\le X_2(t)}\right)\d  B_1 (t)\,,   
 \\
 \label{1.2}  
 \d  X_2 (t) &= 
 \left(g\I{X_1 (t) >  X_2 (t)}- h\I{X_1 (t)\le X_2 (t)}\right)\d t + 
  \left(\rho\I{X_1(t)\le X_2(t)} +\sigma\I{X_1(t)>X_2(t)}\right)\d B_2(t)\,,
\end{align}
initial condition $(X_1 (0), X_2 (0) )= (x_1, x_2 )$, and $B_1(\cdot)$, $B_2
(\cdot)$  two independent, standard Brownian motions. For simplicity, we shall
use throughout the normalization  
\begin{equation}
  \label{2.0}
  \rho^2 + \sigma^2 \,=\,1 
\end{equation}
and refer to the case $\rho \sigma =0$ as   ``degenerate".  

Speaking informally and a bit imprecisely for the moment about the system of
(\ref{1.1})-(\ref{1.2}), imagine you run two Brownian-like particles on the
real line. At any given time, you assign positive  drift $\, g \,$ and
diffusion $\, \sigma \,$ to the laggard; and you assign negative drift $\,
-h\,$ and diffusion $\, \rho \,$ to the leader.   What is the probabilistic
structure of the resulting two-dimensional diffusion process? Can it be
realized as the solution of a system of stochastic differential equations
other than (\ref{1.1}), (\ref{1.2})? What are its transition probabilities?
How does it look like, when time is reversed?

It has been known for some time now, at least for the non-degenerate case (cf.
\cite{MR532498}, pages 193-194; 
\cite{MR917679}; or \cite{MR2078536}, page 45)  that a unique probability
measure 
$\mub$ can be constructed on the canonical filtered measurable space
$(\mathfrak{W}, {\cal G}), \, {\bf G}=  \{ {\cal G} (t) \}_{0 \le t <
  \infty}\,$ of continuous functions $\, \mathfrak{w}: [0,\infty)\to\R^2\,$
endowed with the topology of uniform convergence on compact 
intervals, such that the process $\, f ( \mathfrak{w} (t)) - \int_0^t {\cal L}
f (\mathfrak{w} (s))\, \d  s\,$, $\,0 \le t <\infty\,$ is a   local
martingale under $\, \mub\,$  for every function $\, f : \R^2\to\R\,$ of class
$\, {\cal C}^2  \,$ (here $\, {\cal L}\,$ is the second-order 
partial differential operator in (\ref{0.1}), the infinitesimal generator of
the resulting diffusion). Our goal in this paper is  to describe this
probability measure $\, \mub\,$ as explicitly as possible; to study its
behavior under time-reversal; and to understand the solvability of   systems
of stochastic differential equations, such as (\ref{1.1}), (\ref{1.2}) above,
which correspond to this martingale problem and help realize its  solution.

\smallskip
We shall show in sections \ref{sec2}, \ref{sec3} and \ref{sec4}  that the
system of stochastic differential equations (\ref{1.1}), (\ref{1.2}) has a
solution which is unique in the sense of the probability distribution (thus
the above martingale problem is indeed well-posed). This solution is shown  to
be strong in section \ref{sec4.3}; it is  characterized   in terms of a
one-dimensional diffusion process $\, Y(\cdot)\,$, which has ``bang-bang"
drift with intensity $\, \lambda = g + h >0\,$ and  is  driven by yet another
standard Brownian motion process $\, W(\cdot)\,$, namely  
\begin{equation}
  \label{1.3}
  \d  Y (t) = - \lambda   \sign \big( Y(t) \big)  \d  t +   \d  W  (t) .
\end{equation}
Here and in what   follows  we shall use the convention for the signum
function 
\begin{equation}
  \label{1.3.a}
  \sign (y) \coloneqq%
  \I[(0, \infty)](y) -  \I[{(-\infty, 0]}] (y), \qquad y \in \R.
\end{equation}
The one-dimensional diffusion of (\ref{1.3})  was studied in some detail by 
\cite{MR744236}%
, who found its transition probabilities as in equations
(\ref{3.2}), (\ref{3.3}) below. We shall use this analysis to compute, in
section \ref{sec5},  the transition probabilities of the two-dimensional
process $\, (X_1 (\cdot), X_2 (\cdot))\,$. As in that earlier paper, a crucial
r\^ole will be played here again  by the local time  
\begin{equation}
  \label{1.4}
  L^Y (t) \coloneqq%
  \lim_{\varepsilon \downarrow 0}\, { 1 \over \, 4\, \varepsilon\,} \int_0^t \mathbf{ 1}_{ \{- \varepsilon \, < \,Y(s) \, < \, \varepsilon \,\} }   \, \d  s 
\end{equation}
accumulated at the origin during the interval $\, [0, t]\,$ by the diffusion
$\, Y(\cdot)\,$ of (\ref{1.3}). In terms of this diffusion, its local
time (\ref{1.4}), and   an independent standard Brownian motion $\,
Q(\cdot)\,$, the unique-in-distribution weak solution of the system in
(\ref{1.1}), (\ref{1.2}) will be shown to admit  the {\it skew representation} of (\ref{2.23}), (\ref{2.24}) below. 
  
\smallskip  
The diffusion process of (\ref{1.3}) is also instrumental in a  time-reversal
analysis we carry out in section \ref{sec6}, where  the dynamics of
time-reversed versions of the processes $\, X_1 (\cdot)\,$, $\, X_2 (\cdot)\,$
are derived in the spirit of %
\cite{MR866342} and
\cite{MR1329103}. We were   quite surprised, at first, that these
reverse-time dynamics should feature  terms involving singularly continuous
components such as local times, as indeed they do.  To the best of our
knowledge, this is the first instance where such a structure is observed in
the context of a ``purely forward"  system of  stochastic differential
equations such as (\ref{1.1})-(\ref{1.2}) that does not involve reflection;
see Remark \ref{Werner} in this regard. We also study the forward and backward
dynamics of the ranks $\, R_1 (\cdot) = \max (X_1 (\cdot), X_2 (\cdot))\,$,
$\, R_2 (\cdot) = \min (X_1 (\cdot), X_2 (\cdot))\,$ in this two-dimensional
diffusion, in sections \ref{sec4} and \ref{sec6}, respectively.   

\smallskip
The planar diffusion process with infinitesimal generator (\ref{0.1}) has
local covariance matrix 
\begin{equation}
  \label{Alpha}
  \mathcal{A}(x_1, x_2) \,=\,
  \begin{pmatrix}
    ~ \rho^2 \,  \mathbf{ 1}_{  \{ x_1 > x_2   \} } +  \sigma^2 \,  \mathbf{ 1}_{  \{ x_1 \le x_2   \} } ~   &    ~0~ \\
    ~0  ~    & ~ \rho^2 \,  \mathbf{ 1}_{ \{ x_1 \le x_2  \}  } +  \sigma^2 \,  \mathbf{ 1}_{  \{ x_1 > x_2   \} }~
  \end{pmatrix}\,.
\end{equation} 
There is a continuum of real square roots for this matrix, of the form  
$\,\Sigma(x_{1},x_{2}) = \Sigma_+\,{\bf 1}_{\{ x_{1} > x_{2}\}} + 
\Sigma_-\,{\bf 1}_{\{ x_{1} \le  x_{2}\}}\,$  with  
\begin{equation} 
  \label{Sigma1234} 
  \Sigma_+ \, := \,
  \left( \begin{array}{cc} 
      \rho \cos \varphi & - \rho \sin \varphi \\
      \varepsilon \sigma \sin \varphi & \varepsilon 
      \sigma \cos \varphi \\ 
    \end{array} \right)\,,\qquad  
  \Sigma_- \, := \,
  \left( \begin{array}{cc} 
      \sigma \cos \vartheta & - \sigma\sin \vartheta \\
      \rho \delta \sin \vartheta & \rho\delta \cos \vartheta \\ 
    \end{array} \right)  \, , 
\end{equation}
parametrized by $\, \varepsilon = \pm 1\,$, $\, \delta = \pm 1\,$, $\, 0 \le
\varphi, \vartheta \le 2 \pi\,$. All such configurations lead to systems of
stochastic differential equations that admit a (unique in distribution) weak
solution. We show in subsection \ref{sec9} that {\it those solutions that
  correspond to configurations  with  }
\begin{equation} 
  \label{WEAK} 
  (\sigma^{2} \varepsilon - \rho^{2} \delta) \sin ( \vartheta - \varphi) +  
  \rho \sigma ( 1+ \varepsilon \delta) \cos (\vartheta - \varphi) = -1
\end{equation}
{\it are not strong;} see the system (\ref{2.14}), (\ref{2.15}) for an example. Whereas
{\it all other configurations lead to strongly solvable systems;}   one  such
system  appears in (\ref{1.1}), (\ref{1.2}), and another one in (\ref{2.12}),
(\ref{2.13}).  

Crucial in this analysis of strength and weakness is the
following recent result  by %
\cite{tanaka2009} %
on the pathwise uniqueness of the   ``perturbed Tanaka equation"
\begin{equation}
  \label{Prokaj}
  \Upsilon (t) \,=\, y \,+\,\int_0^t \sign\big( \Upsilon (s) \big) \, \d  M (s)\,+
  \, N (t)\,, \qquad 0 \le t < \infty \,.
\end{equation} 

\begin{theorem}[%
  \cite{tanaka2009}%
  ]
  \label{Theorem 0}  
  Suppose that $\, M (\cdot)\, , \, N(\cdot)\,$ are  continuous local
  martingales  with $\, M(0)=N(0)=0\,$ and quadratic and cross-variations that satisfy the  conditions of  orthogonality $\,  $   and   domination 
  \begin{equation}
    \label{dom}
    \langle M, N \rangle ( t)\, =\,0\,, \qquad    \langle M \rangle ( t) \,= \,\int_0^t q (s) \, \d  \langle N \rangle    (s)\,; 
    \qquad 0 \le t < \infty\,,
  \end{equation}   
  respectively,  for some   progressively measurable process $\, q(\cdot) $ with values in a  compact interval $\, [0, c]\,$.    Under these assumptions, pathwise
  uniqueness holds for the perturbed Tanaka equation (\ref{Prokaj}). 
\end{theorem}

In  section \ref{App}  we shall use the local time techniques introduced by
\cite{MR658680} and further developed by 
\cite{0527.60062}, to provide a   simple 
proof of a considerably  more general result of this type, Theorem
\ref{Theorem 4} (see also Proposition \ref{Prop}), in which the signum
function is replaced in (\ref{Prokaj}) by an arbitrary   function of finite
variation,  and $\,M(\cdot)\,$, $ N(\cdot)\,$ by  continuous semimartingales
such that (\ref{dom}) is satisfied.        

Multidimensional processes of the type (\ref{1.1}), (\ref{1.2}) were
introduced by 
\cite{MR1894767}, and their ergodic behavior was studied
by 
\cite{MR2187296}, \cite{MR2473654}, \cite{2009arXiv0909.0065I}, among others. Here we focus on the
two-dimensional case, and concentrate on the precise probabilistic structure
of the resulting diffusions governed by stochastic differential equations such
as (\ref{1.1}), (\ref{1.2}).  

\section{ \textsc{Analysis} }
\label{sec2}

Let us assume that the system of stochastic differential equations
(\ref{1.1}), (\ref{1.2}) has a weak solution: To wit, that there exists a
filtered probability space $\, ( \Om,   \mathfrak{F}, \p)$,
  $\mathbf{F} = \{
\mathfrak{F} (t) \}_{0 \le t < \infty }\,$ and on it two pairs $\, (B_1
(\cdot), B_2 (\cdot))\,$ and $\, (X_1 (\cdot), X_2 (\cdot))\,$ of continuous,
$\, \mathbf{F}-$adapted processes, such that $\, B_1 (\cdot)\,$ and $\, B_2
(\cdot)\,$ are   independent standard Brownian motions and (\ref{1.1}),
(\ref{1.2}), $\, X_1 (0) =x_1\,$ and $\, X_2 (0) =x_2\,$ hold. We shall 
fix the   nonnegative constants $\, g  \,$, $\, h  \,$, $\, \rho  \,$, $\,
\sigma \,$ with   $\, g+h>0\,$, and   impose the
normalization (\ref{2.0}). 

We shall assume throughout the paper, and without
further mention, that the filtrations we are dealing with are in their
right-continuous versions   and have been augmented by sets of
$\,\mathbb{P}-$measure zero. We shall also use the convention 
$$
\mathbf{F}^{\,\Xi} = \{ \mathfrak{F}^{\,\Xi } (t) \}_{0 \le t < \infty }\,\,,
\qquad \mathfrak{F}^{\,\Xi } (t) \,:=\, \sigmab \big( \Xi  (s)\,, ~ 0 \le s \le t\big)
$$
for the $\,\mathbb{P}-$augmentation of the filtration generated by a given
process $\,\Xi  : [0, \infty) \times \Omega \rightarrow \R^d\,$ with values in
some Euclidean space and RCLL paths.

With such a setup, and with the notation  
\begin{equation}
  \label{2.1}
  \lambda \,=\, g+h\,, \qquad \nu \,=\,   g-h\,, 
  \qquad y\,=\, x_1 - x_2 \,, \qquad z\,=\, x_1 + x_2\,,
\end{equation}
we note that the difference
\begin{equation}
  \label{2.2}
  Y(t)\,:=\, X_1 (t) - X_2 (t)\,, \qquad 0 \le t < \infty
\end{equation}
satisfies the integral version 
\begin{equation}
  \label{3.1}
  Y(t) \,=\, y - \lambda \int_0^t \sign \big( Y(s)\big)\, \d  s + W(t)\,, 
  \qquad 0 \le t < \infty
\end{equation}
of the equation   (\ref{1.3}). The equation (\ref{3.1}) is  driven by the process
\begin{equation}
  \label{2.3}
  W (t) \,:=\, \rho\, W_1 (t) \,+\, \sigma \, W_2 (t)\,, \qquad 0 \le t < \infty \,,
\end{equation} 
where we have set 
\begin{align}
  \label{2.4}
  W_1 (t) \,&:=\, \int_0^t \mathbf{ 1}_{ \{ Y (s) > 0 \} } \, \d  B_1
  (s) - \int_0^t \mathbf{ 1}_{ \{ Y (s) \le 0 \} } \, \d  B_2 (s)  \,,
  \\
  \label{2.5}
  W_2 (t) \,&:=\, \int_0^t \mathbf{ 1}_{ \{ Y (s) \le 0 \} } \, \d  B_1 (s) - \int_0^t \mathbf{ 1}_{ \{ Y (s) > 0 \} } \, \d  B_2 (s)  \,.
\end{align} 
It is also seen from (\ref{1.1}) and  (\ref{1.2}) that the sum of the two
processes $\, X_1 (\cdot)\,$, $\, X_2 (\cdot)\,$ is of the form  
\begin{equation}
  \label{2.6}
  X_1 ( t)+ X_2 ( t)\,=\, z + \nu \, t + V(t)\,,
  \qquad V(  t) \,:=\, \rho\, V_1 ( t) + \sigma \, V_2 ( t)\,,
\end{equation}
where, by analogy with (\ref{2.4}), (\ref{2.5}) we have set 
\begin{align}
  \label{2.7}
  V_1(t)\,  &:=\,  \int_0^t    \mathbf{ 1}_{ \{ Y (s) > 0 \} } \, \d 
  B_1 (s) +   \int_0^t \mathbf{ 1}_{ \{ Y (s) \le 0 \} }  \, \d  B_2
  (s)   \,, 
  \\
  \label{2.8}
  V_2(t)  \,&:= \, \int_0^t    \mathbf{ 1}_{ \{ Y (s) \le 0 \} } \, \d  B_1 (s) + \int_0^t  \mathbf{ 1}_{ \{ Y (s) > 0 \} }  \, \d  B_2 (s)   \,.
\end{align}

\noindent
$\bullet~$ 
The processes $\, W_1 (\cdot)\,$,  $\, W_2 (\cdot)\,$,  $\, V_1 (\cdot)\,$,
$\, V_2 (\cdot)\,$ are continuous $\mathbf{F}-$martingales   with quadratic
variations  $  \langle W_1 \rangle (t) =  \langle W_2 \rangle (t) =  \langle
V_1 \rangle (t) =  \langle V_2 \rangle (t) =  t$, i.e., Brownian motions by
the \textsc{P. L\'evy} theorem (e.g., 
\cite{MR1121940}%
, p.$\,$157);   also $  \langle V_1, V_2 \rangle
(\cdot)= \langle W_1, W_2 \rangle (\cdot)=\langle V_1, W_2 \rangle (\cdot)=
\langle W_1, V_2 \rangle (\cdot)=0 $.

The pairs $\big( W_1 (\cdot),W_2(\cdot) \big)$ and  $\big(V_1(\cdot),V_2(\cdot) \big)$, as well as  $\big( V_1 (\cdot), W_2 (\cdot)
\big) $ and $\big( W_1 (\cdot), V_2 (\cdot) \big) $,  are thus
two-dimensional Brownian motions,    by the \textsc{P. L\'evy} theorem. In
conjunction with (\ref{2.0}) this implies, in particular, that     $\, W
(\cdot)\,$ and $\, V(\cdot)\,$ in (\ref{2.3}), (\ref{2.6}) are   standard,
one-dimensional Brownian motions.  

\smallskip
\noindent
$\bullet~$ On the other hand, we see from  (\ref{2.4}), (\ref{2.7})   and with
the notation of (\ref{1.3.a}) the intertwinements  
\begin{align}
  \label{2.9}
  V_1(t) \,&=\, \int_0^t \sign \big( Y (s) \big) \, \d  W_1 (s)\,,  
  &%
  W_1(t) \,&=\, \int_0^t \sign \big( Y (s) \big) \, \d  V_1 (s)  \,;
  \shortintertext{and from (\ref{2.5}), (\ref{2.8}) we obtain  the intertwinements} 
  \label{2.10}
  V_2(t) \,&=\, -\int_0^t  \sign \big( Y (s) \big) \, \d  W_2 (s)
  \,,  &%
  W_2(t) \,&=\, -\int_0^t  \sign \big( Y (s) \big) \, \d  V_2 (s) \,.
\end{align}
$\bullet~$ 
The  equation (\ref{3.1}) admits a pathwise unique, strong solution; in
particular, uniqueness in the sense of the probability distribution holds as
well (cf. Proposition 5.5.17, page 341 and the \textsc{Yamada-Watanabe}
Proposition 5.3.20, page 309 in 
\cite{MR1121940}), and we
have the identity $\, \mathbf{ F}^{\,Y} \equiv \mathbf{ F}^{\,W}\,$.

\subsection{ \textsc{Two Auxiliary Systems} }
\label{sec2.1}

The relations 
$\, 
X_1 ( t)+ X_2 ( t)\,=\, x_1 + x_2 +\big( g - h \big)  \, t +  
\rho\, V_1 ( t) + \sigma \, V_2 ( t)
\,$ 
from (\ref{2.6}), and 
\begin{multline*}
X_1 ( t)- X_2 ( t)=\\ 
 x_1 - x_2 +\big( g + h \big) \int_0^t  \left( 
  \I{ X_1 (s) \le  X_2 (s)} - \I{X_1 (s) >  X_2 (s)} \right)\d s + \rho W_1 ( t) + \sigma  W_2 ( t)   
\end{multline*}
from   (\ref{1.3}), (\ref{2.2}), (\ref{2.3}) lead, in conjunction with the
intertwinements  of (\ref{2.9}) and (\ref{2.10}), to the system of  stochastic
integral equations 
\begin{align}
  \nonumber
  X_1 (t) &=x_1 + \int_0^t   \left(g\,\I{X_1 (s) \le  X_2(s)} - h\,\I{X_1 (s)> X_2 (s)} \right) \d s\\
  &\hphantom{=\,x_1+}
  + \rho  \int_0^t \I{X_1 (s) >  X_2 (s)}\d W_1(s)+ \sigma  \int_0^t \I{X_1
    (s) \le   X_2 (s)} \d W_2(s),   \label{2.12} \\
  \nonumber
  X_2 (t) \,&=\,x_2 + \int_0^t   \left(\, g\,\I{X_1 (s) >  X_2(s)} - h\,\I{X_1 (s) \le  X_2 (s)} \right) \d  s 
  \\
  &\hphantom{=\,x_1+}
  - \rho  \int_0^t \I{X_1 (s) \le  X_2 (s)} \d W_1 (s)
  - \sigma \int_0^t \I{X_1 (s) >   X_2 (s)} \d W_2 (s) 
 \label{2.13}  
\end{align}
\noindent
driven by the planar Brownian motion  $\, \big( W_1 (\cdot), W_2 (\cdot) \big)
$ of (\ref{2.4}), (\ref{2.5}).   

\smallskip
Furthermore, we observe that the intertwinements of (\ref{2.9}), (\ref{2.10})
allow us to recast this system as driven by the planar Brownian motion  $\,
\big( V_1 (\cdot), V_2 (\cdot) \big) $ of (\ref{2.7}), (\ref{2.8}), namely  
\begin{align}
  \nonumber 
  X_1 (t) \,&=\,x_1 + \int_0^t   \left(\, g\,\mathbf{ 1}_{ \{ X_1 (s) \le  X_2
      (s) \} } - h\,\mathbf{ 1}_{ \{ X_1 (s) >  X_2 (s) \} } \right) \,
  \d  s \\  \label{2.14}
  &\hphantom{= x_1 +} 
  +\, \rho  \int_0^t \mathbf{ 1}_{ \{ X_1 (s) >  X_2 (s) \} }\, \d  V_1
  (s) \,+\,    \sigma  \int_0^t \mathbf{ 1}_{ \{ X_1 (s) \le   X_2 (s) \} }\,
  \d  V_2 (s) \,,\\ 
  \nonumber
  X_2 (t) \,&=\,x_2 + \int_0^t   \left(\, g\,\mathbf{ 1}_{ \{ X_1 (s) >  X_2
      (s) \} } - h\,\mathbf{ 1}_{ \{ X_1 (s) \le  X_2 (s) \} } \right) \,
  \d  s \\ \label{2.15}
  &\hphantom{= x_1 +} 
  + \, \rho  \int_0^t \mathbf{ 1}_{ \{ X_1 (s) \le  X_2 (s) \} }\, \d  V_1 (s) \,+\,    \sigma  \int_0^t \mathbf{ 1}_{ \{ X_1 (s) >   X_2 (s) \} }\, \d  V_2 (s) \,.
\end{align}
\smallskip
It is quite clear, though perhaps worth noting, that  the systems of
stochastic   equations (\ref{1.1})-(\ref{1.2}),   as well as
(\ref{2.12})-(\ref{2.13}) and (\ref{2.14})-(\ref{2.15}), give rise and
correspond to the same martingale problem  -- namely,  the one with
infinitesimal generator (\ref{0.1}).

\subsection{\textsc{Skew and Integral Representations}}
\label{sec2.2}

By analogy with (\ref{2.3}), (\ref{2.6}) let us introduce the standard
Brownian motions  
\begin{equation}
  \label{2.17}
  W^{\,\mathbf{ \flat}}( \cdot) := \rho\, W_1 (\cdot) - \sigma \, W_2
  (\cdot)\, , 
  \qquad   V^{\,\flat}( \cdot) := \rho\, V_1 (\cdot) - \sigma \, V_2 (\cdot)\,   
\end{equation}
and note  the new intertwinements
\begin{equation}
  \label{2.18}
  V (t) \,=\, \int_0^t \,\sign \big( Y(s) \big)\,\d  W^{\,\mathbf{
      \flat}}  (s)\,,  
  \qquad  V^{\,\flat} (t) \,=\, \int_0^t\,\sign \big( Y(s) \big)\,\d  W  (s)\,.   
\end{equation}
It will be convenient to cast the Brownian motion $\, W^{\,\mathbf{ \flat}}(
\cdot)\,$ of (\ref{2.17}) in the decomposition   
\begin{equation}
  \label{2.19}
  W^{\,\flat} (\cdot) = \gamma \, W(\cdot) + \delta \, U^{\,\flat}
  (\cdot)\,, 
  \quad \mathrm{where}\quad 
  \gamma \coloneqq \rho^2 - \sigma^2\,, \quad 
  \delta \coloneqq \sqrt{1 - \gamma^2\,\,} = 2 \rho\sigma\,,
\end{equation}
in terms of  the independent Brownian motions  
\begin{equation}
  \label{2.20}
  U^{\flat} ( \cdot) := \sigma \, W_1 (\cdot) - \rho\, W_2 (\cdot)     
  \qquad \hbox{and} \qquad 
  W( \cdot) = \rho\, W_1 (\cdot) + \sigma \, W_2 (\cdot) 
\end{equation}
as in (\ref{2.3}). With this setup, the  Brownian motion $\, V(\cdot)\,$
defined  in (\ref{2.6}) takes   the form $\, \,V (t) = \int_0^t \,  \sign
\big( Y(s) \big) \,$ $   [ \, \gamma\,  \d  W(s) \,+\, \delta \,
\d  U^{\flat} (s)\, ]\,$, equivalently 
\begin{align}
  \nonumber
  V(t)\,&=\, \gamma \int_0^t \,  \sign \big( Y(s) \big) \, \Big[ \,
  \d  Y  (s)+  \lambda \, \sign \big( Y(s) \big)\,   \d s \,
  \Big]+ \delta\, Q (t )    \\
  \label{2.21}
  \,&=\, \gamma \, \Big( \big| Y (t)   \big| -
  |y| \,+\, \lambda\, t\,-\,  2 \, L^Y (t) \Big) + \delta\,Q (t )\,, 
  \qquad 0 \le t < \infty  
\end{align}
thanks to (\ref{2.18}), (\ref{2.19}), (\ref{1.3}) and the \textsc{Tanaka}
formulas (e.g., 
\cite{MR1121940}%
, page 220). Here  
\begin{equation}
  \label{2.22}
  Q (t) := \int_0^t   \sign \big( Y(s) \big) \d  U^{\flat}  (s) = \sigma V_1 (t) + \rho V_2 (t), \qquad 0 \le t < \infty 
\end{equation}
is yet another Brownian motion with $\q{Q, W} =  \q{U^{\flat}, W} = 0$. Thus,
$\, Q(\cdot)\,$ is   independent of the 
Brownian motion $\,W (\cdot)\,$   and  of the diffusion process $Y(\cdot)$
that $W (\cdot)$ engenders via (\ref{3.1}).    

\smallskip
\noindent
$ \bullet~$ We can combine now the last expression in (\ref{2.21})  with
$X_1 (t) - X_2 (t) = Y (t)$ of (\ref{2.2}) and with $X_1 (t) + X_2 (t) = x_1 +
x_2+ \nu  t +  V(t)$ of (\ref{2.6}),  and  arrive at the {\it skew
  representations}      
\begin{align}
  \label{2.23}
  X_1 (t)&=x_1 + \mu t + \rho^2 \big(Y^{+}(t) - y^{+}\big) - 
  \sigma^2 \big(Y^{-}(t) - y^{-}\big)- \gamma\, L^Y(t) + \rho\sigma Q (t),\\
  \label{2.24}
  X_2 (t)&=x_2 + \mu t - \sigma^2 \big(Y^{+} (t) - y^{+} \big) + \rho^2 \big(
  Y^{-} (t) - y^{-} \big) - \gamma\, L^Y (t)+   \rho\sigma Q (t) 
\end{align}
 for the components of the two-dimensional diffusion $\, ( X_1 (\cdot), X_2 (\cdot))\,$   of the system (\ref{1.1}), (\ref{1.2});   we have set  
\begin{equation}
  \label{2.25}
  \mu := \frac12 \big( \nu + \lambda\gamma \big)= g\rho^2 - h \sigma^2 \,.
\end{equation} 

These formulas involve the positive and negative parts of the current value of
the one-dimensional diffusion process $Y(\cdot)$ in (\ref{1.3}), the
current value of its local time $L^Y (\cdot)$ at the origin, and the
current value of the {\it independent} Brownian motion $Q(\cdot)$ of
(\ref{2.22}). Two cases stand out.  

\medskip
\noindent
{\bf (A)}~ In the {\it Equal Variance}  (Isotropic) case $\, \rho = \sigma = 1
/ \sqrt{2\,}\,$, the local times disappear from the expressions (\ref{2.23}),
(\ref{2.24}), which then take the very simple form 
\begin{equation}
  \label{2.26}
  X_1 (t) = x_1 + \frac12\big( \nu t + Y(t) - y + Q (t) \big)\,, \quad 
  X_2 (t) = x_2 + \frac12\big( \nu t - Y(t) + y + Q (t) \big)\,.
\end{equation}

\smallskip
\noindent
{\bf (B)}~  In the {\it Degenerate}    case $\rho\sigma = 0$ the
independent Brownian motion $\,Q(\cdot)\,$ disappears from these expressions;
for instance, with $\, \sigma =0\,$ and $\, \rho =1\,$, they become  
\begin{equation}
  \label{2.27}
  X_1 (t)\,=\, x_1 - y^+ + g \,t + Y^+ (t ) - L^Y (t)\,, \qquad 
  X_2 (t)\,=\, x_2 - y^- + g \,t + Y^- (t ) - L^Y (t)\,.
\end{equation}

\subsection{ \textsc{Uniqueness in Distribution} }
\label{sec2.3}

The analysis of this section shows that, given {\it any} weak solution of the
system of stochastic equations (\ref{1.1}), (\ref{1.2}), its  vector process
$\, (X_1 (\cdot), X_2 (\cdot))\,$ can be cast in the form
(\ref{2.23})--(\ref{2.24}). Here the diffusion process $\, Y(\cdot)\,$ is the
pathwise unique, strong solution of the stochastic integral  equation   
$$\, 
Y(t) = y - \lambda \int_0^t \sign \big( Y(s)\big)\, \d  s + W(t)\, , \qquad \,  0 \le t < \infty\,$$
as in (\ref{3.1}) and with the notation of (\ref{1.3.a}), driven by the
Brownian motion $\, W(\cdot)\,$ of (\ref{2.3}); whereas the Brownian motion
$\, Q(\cdot)\,$ is independent of $\, W(\cdot)\,$, thus also of $\,
Y(\cdot)$.  In other words, the joint distribution of the pair $\, (
Y(\cdot), Q(\cdot)) \,$ is determined uniquely -- and thus, from
(\ref{2.23})--(\ref{2.24}), so is the joint distribution of the vector process
$\, ( X_1 (\cdot), X_2 (\cdot)) \,$.  

To put it a bit more succinctly: {\it uniqueness in distribution holds for the
  system of equations (\ref{1.1}), (\ref{1.2}), as well as for the systems of
  equations (\ref{2.12}), (\ref{2.13}) and (\ref{2.14}), (\ref{2.15}). } 

\section{ \textsc{Synthesis} }
\label{sec3}

Let us begin now to reverse  the steps of the preceding analysis. We start
with a filtered probability space $\, ( \Om,   \mathfrak{F},
\p)\,$, $
\mathbf{F} = \{ \mathfrak{F} (t) \}_{0 \le t < \infty }\,$ rich enough to
support two independent, standard Brownian motion $\,  W_1  (\cdot) \,$ and
$\,  W_2  (\cdot) \,$; and without sacrificing generality, we shall assume
$\, \mathbf{F} \equiv \mathbf{F}^{\, (W_1, W_2)}\,$, i.e., that the filtration
is generated by this planar Brownian motion.  

With     given nonnegative constants $\, g  \,$, $\, h  \,$, $\, \rho  \,$,
$\, \sigma \,$ that satisfy (\ref{2.0}) and   $\, g+h>0\,$, with  a given
vector $\, (x_1, x_2 ) \in \R^2\,$, and with the notation of (\ref{2.1}), we
construct the pairs of independent Brownian motions  
\begin{align}
  \label{3.5.a}
  W(\cdot) \,:&=\, \rho\, W_1 (\cdot) \,+\, \sigma \, W_2 (\cdot)\,,&
  U^{\, \flat}(\cdot) \,:&=\, \sigma \, W_1 (\cdot) \,-\, \rho  \, W_2 (\cdot) \\
  \shortintertext{and} 
  \label{3.5.ab}
  U(\cdot) \,:&=\, \sigma \, W_1 (\cdot) \,+\, \rho \, W_2 (\cdot)\,,&
  W^{\, \flat}(\cdot) \,:&=\, \rho\, W_1 (\cdot) \,-\, \sigma \, W_2 (\cdot)  
\end{align}
in accordance with (\ref{2.20}), (\ref{2.17}). Clearly, $\,   \mathbf{F}^{\, (W_1, W_2)} \equiv \mathbf{F}^{\, (W , U^{\flat})} \equiv \mathbf{F}^{\, (U, W^{  \flat})}\,$.  

We construct also the pathwise unique, strong solution $\, Y(\cdot)\,$ of the
stochastic equation (\ref{3.1}) driven by the Brownian motion $\, W(\cdot)\,$
in (\ref{3.5.a}).   This    is a strong \textsc{Markov} and \textsc{Feller}
process, whose transition probabilities can be computed explicitly; see (\ref{3.2})-(\ref{3.4}) below. 

\smallskip
With the process $\, Y(\cdot)\,$ thus in place, we introduce the continuous,
$\mathbf{F}-$adapted processes 
\begin{equation}
  \label{3.5}
  V_1(t) \,=\, \int_0^t \sign \big( Y (s) \big) \, \d  W_1 (s)\,,  \qquad  
  V_2(t) \,=\, -\int_0^t  \sign \big( Y (s) \big) \, \d  W_2 (s)
\end{equation}
in accordance with (\ref{2.9}) and (\ref{2.10}). These are martingales   with
$\, \langle V_1 \rangle (t) =  \langle V_2 \rangle (t) = t\,$ and $\, \langle
V_1, V_2 \rangle (t) =0\,$ for all $ 0 \le t < \infty\,$, thus independent
Brownian motions in their own right. We construct from them, and by analogy
with (\ref{3.5.a}) and (\ref{3.5.ab}), two additional pairs of independent,
standard Brownian motions, namely  
\begin{align}
  \label{3.5.b}
  V(\cdot) \,:&=\, \rho\, V_1 (\cdot) \,+\, \sigma \, V_2 (\cdot)\,, &
  Q^{\, \flat}(\cdot) \,:&=\, \sigma \, V_1 (\cdot) \,-\, \rho  \, V_2 (\cdot) 
  \shortintertext{and} 
  \label{3.5.bb}
  Q(\cdot) \,:&=\, \sigma \, V_1 (\cdot) \,+\, \rho \, V_2 (\cdot)\,, &
  V^{\, \flat}(\cdot) \,:&=\, \rho\, V_1 (\cdot) \,-\, \sigma \, V_2 (\cdot)\,.
\end{align}
We note the  intertwinements (\ref{2.18}), (\ref{2.22}) and $\,
Q^{\flat} (\cdot) = \int_0^\cdot  \,  \sign \big( Y(t) \big) \, \d  U  (t)    \,$, as well as the filtration identities $\,   \mathbf{F}^{\, (V_1, V_2)} \equiv \mathbf{F}^{\, (V , Q^{\flat})} \equiv \mathbf{F}^{\, (Q, V^{  \flat})}\,$.  
 
\smallskip
Finally, we introduce   the continuous, $\, \mathbf{F}-$adapted processes 
\begin{align}
  \label{3.6}
  X_1(t) \,:&=\,x_1 + 
  \int_0^t  \left( \,g\, \mathbf{ 1}_{ \{ Y (s) \le 0 \} } - 
    h\,\mathbf{ 1}_{ \{ Y (s) > 0 \} } \right) \, \d  s \,+\,M_1 (t)\\
  \label{3.7}
  X_2(t) \,:&=\,x_2 +    \int_0^t  \left( \,g\, \mathbf{ 1}_{ \{ Y (s) > 0 \} } - h\,\mathbf{ 1}_{ \{ Y (s) \le 0 \} } \right) \, \d  s \,+\, M_2 (t)  
\end{align}
for $\,0 \le t < \infty \,$, where we set 
\begin{align}
  \label{3.8}
  M_1 (t) \,:&=\, 
  \int_0^t  \Big( \, \rho \, \mathbf{ 1}_{ \{ Y (s) > 0 \} } \, \d  W_1
  (s)  + \sigma \,  \mathbf{ 1}_{ \{ Y (s) \le 0 \} }  \, \d  W_2 (s)
  \Big)\,,\\
  \label{3.9}
  M_2 (t) \,:&=\, 
  \int_0^t  \Big(  - \rho \, \mathbf{ 1}_{ \{ Y (s) \le  0 \} } \, \d  W_1 (s)  - \sigma \,  \mathbf{ 1}_{ \{ Y (s) > 0 \} }  \, \d  W_2 (s) \Big)\,,  
\end{align}
by analogy with the equations of (\ref{2.12}), (\ref{2.13}). We have for these processes 
\begin{equation}
  \label{3.10}
  X_1 ( t) - X_2 ( t)\, =\, Y ( t)\,,      \qquad 
  X_1 ( t)+ X_2 ( t)\,=\, x_1 + x_2 + \nu  \, t +  
  V  ( t) 
\end{equation}
in accordance with  (\ref{2.2}), (\ref{2.6}), and note that $\, M_1
(\cdot)\,$, $\, M_2 (\cdot)\,$  are continuous $\, \mathbf{F}-$mar\-tin\-gales,
with $ \, \langle M_1, M_2 \rangle$ $ ( \cdot) \equiv   0\,$ and quadratic
variations 
\begin{align*}
\q{M_1} (t) &= \int_0^t \big(\rho^2 \I{Y(s) > 0} + \sigma^2 \I{Y (s) \le 0}
\big) \d  s\,, \\
\q{M_2} (t) &= \int_0^t \big(\rho^2  \I{ Y (s) \le 0} + \sigma^2\I{Y (s) > 0}  \big) \d  s\,. 
\end{align*}
There exist then independent Brownian motions $B_1 (\cdot)$, $B_2(\cdot)$ on
our filtered probability space $(\Om,\mathfrak{F},\P)$, $\mathbf{F} = \{
\mathfrak{F} (t) \}_{0 \le t < \infty }$, so   the continuous martingales of
(\ref{3.8}), (\ref{3.9}) are cast in their \textsc{Doob}  representations  
\begin{equation}
  \label{3.11} 
  \begin{aligned}
  M_1 (t) &= \int_0^t\Big(\rho \I{Y (s) > 0 } + \sigma \I{Y(s) \le 0}
  \Big)\d  B_1 (s)\,, \\ \quad 
  M_2 (t) &= \int_0^t\Big(\rho \I{Y (s)\le 0} + \sigma \I{Y(s)> 0}\Big)\d B_2 (s)  
  \end{aligned}
\end{equation}
for $\, 0 \le t < \infty\,$; for instance, we can take the Brownian motions 
\begin{align}
  \label{3.12} 
  B_1 (t) \, &=  \,\hphantom{-}\int_0^t  \Big(  \, \mathbf{ 1}_{ \{ Y (s) > 0 \} } \,
  \d  W_1 (s)  +  \mathbf{ 1}_{ \{ Y (s) \le 0 \} }  \, \d  W_2
  (s) \Big)\,, \\
  \label{3.13} 
  B_2 (t)  \, & =    \, -\int_0^t  \Big(  \mathbf{ 1}_{ \{ Y (s) \le 0 \} } \,
  \d  W_1 (s)  +  \mathbf{ 1}_{ \{ Y (s) > 0 \} }  \, \d  W_2
  (s) \Big)   
\end{align}
that one gets by disentangling $\, (B_1 (\cdot), B_2 (\cdot))\,$ from $\, (W_1
(\cdot), W_2 (\cdot))\,$ in (\ref{2.4}), (\ref{2.5}).

\subsection{ \textsc{Taking Stock} }
\label{sec3.1}

To recapitulate: {\it we have constructed a weak solution for the system of
  stochastic differential equations (\ref{1.1}), (\ref{1.2})}, as is seen
clearly from (\ref{3.6}), (\ref{3.7}) and (\ref{3.11}); and as we stressed in
subsection \ref{sec2.3}, {\it this solution is unique in distribution.}  We
remarked in the Introduction that this is in accordance with general results
of  %
\cite{MR532498}%
, pages 193-194; %
\cite{MR917679}%
; or %
\cite{MR2078536}%
, page 45.

It develops from (\ref{3.6})--(\ref{3.9})   that {\it we have also constructed a
  weak solution to the system   (\ref{2.12}), (\ref{2.13}), once
  again unique in distribution.} On the other hand,  we can express the
martingales of (\ref{3.8}), (\ref{3.9}) as  
\begin{align*}
  M_1 (t) &= 
  \int_0^t  \Big( \, \rho \, \mathbf{ 1}_{ \{ Y (s) > 0 \} } \, \d  V_1 (s)  + \sigma \,  \mathbf{ 1}_{ \{ Y (s) \le 0 \} }  \, \d  V_2 (s) \Big)\,,  \\
  M_2 (t) &=
  \int_0^t  \Big(    \rho \, \mathbf{ 1}_{ \{ Y (s) \le  0 \} } \, \d  V_1 (s)  + \sigma \,  \mathbf{ 1}_{ \{ Y (s) > 0 \} }  \, \d  V_2 (s) \Big)\,,  
\end{align*}
in terms of the independent, standard Brownian motions of (\ref{3.5}). Back
into (\ref{3.6}) and (\ref{3.7}), these expressions show that {\it we have
  also constructed a weak solution for the system of stochastic equations
  (\ref{2.14}), (\ref{2.15}), once again unique in the sense of the
  probability distribution.  }

In the next two sections we shall discuss in detail properties of
strength/weakness for the  solutions to these systems of equations, namely,
(\ref{2.14})-(\ref{2.15}), (\ref{2.12})-(\ref{2.13}) and
(\ref{1.1})-(\ref{1.2}). For the moment, let us remark that the Brownian
motion $\, V(\cdot) = \rho\, V_1 (\cdot) \,+\, \sigma \, V_2 (\cdot)\,$
determines the {\it sum} of $\, X_1 (\cdot),\, X_2 (\cdot)\,$ as in
(\ref{2.6}) or (\ref{3.10}); and that the Brownian motion $\, W(\cdot) =
\rho\, W_1 (\cdot) \,+\, \sigma \, W_2 (\cdot)\,$ determines the {\it
  difference}   $\, Y(\cdot) = X_1 (\cdot)- X_2 (\cdot)\,$ via the solution of
the stochastic equation (\ref{3.1}).

\subsection{ \textsc{A Two-Dimensional Analogue of the Skew Brownian Motion} }

The equations of (\ref{2.23})-(\ref{2.24}) can be cast in the form 
\begin{align}
  \label{2.23.a}
  X_1 (t) \,&=\, x_1 + \mu\, t + \rho^2 \,\Big( \,\big( X_1 (t) - X_2 (t)
  \big)^+  - \big( x_1 - x_2 \big)^+ \Big) \\ \nonumber
  &\hphantom{=x_1} 
  - \,\sigma^2 \Big( \,\big( X_1 (t) - X_2 (t) \big)^-   - \big( x_1 - x_2
  \big)^- \Big)
  - (\rho^2 - \sigma^2)\, L^{X_1-X_2} (t) \,+  \, \rho\, \sigma\, Q (t)\,,\\
  \label{2.24.a}
  X_2 (t) \,&=\, x_2 + \mu\, t - \sigma^2 \,\Big( \,\big( X_1 (t) - X_2 (t)
  \big)^+   - \big( x_1 - x_2 \big)^+ \Big) \\ \nonumber
  &\hphantom{=x_1}
  +\, \rho^2 \,\Big( \,\big( X_1 (t) - X_2 (t) \big)^-   - \big( x_1 - x_2 \big)^- \Big) - (\rho^2 - \sigma^2)\, L^{X_1-X_2} (t) \,+  \, \rho\, \sigma\, Q (t)\,,
\end{align}
namely, as a system of stochastic differential equations involving the local
time $\, L^{X_1- X_2} (\cdot)\,$ at the origin of the difference $\, X_1
(\cdot)- X_2(\cdot)\,$. {\it We have constructed a   solution to this system,
  subject to the requirement that the difference $\, X_1 (\cdot)-
  X_2(\cdot)\,$ be independent of the driving Brownian motion $\,
  Q(\cdot)\,$.}      

\smallskip
The system   (\ref{2.23.a}), (\ref{2.24.a}) can be construed as a
two-dimensional analogue of the stochastic   equation derived by
\cite{MR606993} %
for the %
\cite{Walsh1978} %
{\it skew Brownian motion;} see the recent survey by 
\cite{MR2280299}%
, and the references cited there,  for the various constructions and
properties of this process.

\section{ \textsc{Ranks and Filtrations} }
\label{sec4}

Let us look now at the ranked versions 
\begin{equation}
  \label{4.1}
  R_1 (\cdot) \,=\, \max ( X_1 (\cdot), X_2 (\cdot) )\, , 
  \qquad \, 
  R_2 (\cdot)\, = \,\min ( X_1 (\cdot), X_2 (\cdot) )
\end{equation}
of the components of the diffusion process constructed in the previous
section. We have $\, R_1 (t) + R_2 (t)=  X_1 (t) + X_2 (t) = x_1 + x_2 + \nu
\,t + V(t)\,$ from  (\ref{3.10}), as well as  
\begin{equation}
  \label{4.2}
  \begin{aligned}[b]
    R_1 (t) - R_2 (t) \,&=\, \big|\, X_1 (t) - X_2 (t)\, \big| \,=\,  
    \big| Y (t)   \big|\,=\, |y| +  \int_0^t \sign \big(  Y (s)  \big) \, \d Y ( s) \,+\, 2 \, L^Y (t)\\
    &=\, |x_1 - x_2| \,-\, \lambda\, t\,+\,V^{\,\flat} (t) \,+\, 2 \, L^Y (t)\,, 
    \qquad 0 \le t < \infty  
  \end{aligned}
\end{equation}
from the \textsc{Tanaka} formulas and (\ref{3.1}), (\ref{2.18}); and from the
theory of the \textsc{Skorokhod} reflection problem (e.g., 
\cite{MR1121940}%
, page 210), we note  
\begin{equation}
  \label{4.5}
  2 \, L^Y (t)\, =\, \max_{0 \le s \le t} \Big( - \big( \,|y| +  V^{\,\flat}
  (s) -  \lambda \, s \,\big) \Big)^+\,, \qquad 0 \le t < \infty\,. 
\end{equation}
We also note from (\ref{4.2}), (\ref{4.5}) and (\ref{1.4}) the filtration
relations 
\begin{equation}
  \label{4.3.a}
\mathfrak{ F}^{\,V^\flat} (t)\,=\,  \mathfrak{ F}^{\,|Y|} (t) \,  \subsetneqq \,\mathfrak{F}^Y(t)  \,, 
  \qquad 0 < t < \infty\,;
\end{equation}
the strict inclusion comes from the fact that $\, Y(\cdot)\,$ has a zero on $\, [0,t]\,$ with positive probability, so that the random variable  $\,\sign(Y(t))$    is not    $\,\mathfrak{F}^{\,\abs{Y}}(t)-$measurable. Alternatively, one can note from (\ref{5.2}) below  that both $\, \P \, \big(Y(t)>0 \,|\, \mathfrak{ F}^{\,|Y|} (t) \big)\,$ and $\, \P \, \big(Y(t)<0 \,|\, \mathfrak{ F}^{\,|Y|} (t) \big)\,$ are non-trivial. 

With   $\, r_1 = \max (x_1, x_2)\,$, $\, r_2 = \min (x_1, x_2)\,$ we deduce
from these equations and (\ref{3.5.b}), (\ref{3.5.bb}) the dynamics for the
ranks  
\begin{align}
  \label{4.3}
  R_1 (t)\,&=\, r_1 - h\, t + \rho \, V_1(t) + L^Y (t)\,,  \qquad 
  0 \le t < \infty\\
  \label{4.4}
  \,
  R_2 (t)\,&=\, r_2 + g\, t + \sigma \, V_2(t)- L^Y (t)\,,  \qquad 
  0 \le t < \infty \,.
\end{align}
Equations (\ref{4.3}), (\ref{4.4})  identify the processes $\, V_1(\cdot)\,$
and $\, V_2(\cdot)\,$ of (\ref{3.5}) as  the independent Brownian motions
associated with   individual {\it ranks;} such an interpretation is also
possible using the equations of (\ref{2.14}), (\ref{2.15}). On the other hand,
the independent, standard Brownian motions $\, B_1 (\cdot)\,$, $\, B_2
(\cdot)\,$ of (\ref{3.12}) and (\ref{3.13}),  are those  associated with the
{\it ``names"} (indices, identities) of the individual particles.

\subsection{ \textsc{The Degenerate Case} }
\label{sec4.1}

We   embark now on a detailed study of the solutions to the systems of
equations (\ref{1.1}), (\ref{1.2}), as well as of those of (\ref{2.12}),
(\ref{2.13}) and (\ref{2.14}), (\ref{2.15}). In this effort it is instructive
to focus first on the degenerate case $\, \rho \, \sigma =0\,$, which has some
interesting features of its own. 

\begin{proposition}
  \label{Proposition 1}
  In the degenerate case with $\,   \sigma =0\,$, thus $\, \rho   =1\,$ in
  light of (\ref{2.0}), we have   the relations 
  \begin{equation}
    \label{4.6}
    \mathfrak{F}^{\,(R_1, R_2)} (t)=\mathfrak{F}^{\,V} (t)  = \mathfrak{F}^{\,X_1 + X_2} (t)  \,   \subsetneqq  \, \mathfrak{F}^{\,(X_1, X_2)} (t) =\mathfrak{F}^{\,Y} (t)= \mathfrak{F}^{\,X_1 - X_2} (t)=\mathfrak{F}^{\,W} (t)  
  \end{equation}
  for every $\, 0 < t < \infty\,$, where the inclusion is strict.
\end{proposition}
\begin{proof}  
  With $\, \rho =1\,$, $\, \sigma =0\,$ we have $\, V_1 (\cdot) =
V(\cdot) = V^{\, \flat} (\cdot)\,$ in (\ref{3.5.b}), (\ref{3.5.bb});   the
first of the claims in (\ref{4.6}) follows from (\ref{4.5}) and
(\ref{4.3})-(\ref{4.4}), which read now 
\begin{equation}
  \label{4.7.a}
  R_1 (t)\,=\, r_1 - h\, t + V(t) + L^Y (t)\,,  \qquad  ~\,
  R_2 (t)\,=\, r_2 + g\, t  - L^Y (t)\,,  \qquad 0 \le t < \infty\,.
\end{equation}
The   second claim is immediate from (\ref{3.10}). 

The third (strict inclusion) is fairly obvious from (\ref{2.27}). Indeed, the  equations of (\ref{2.27}) show 
$$
X_1 (t) + X_2 (t) \,=\, z - |y| + 2\,g \,t + \big| Y(t) \big| - 2 \, L^Y (t)\,,
$$
thus $\, \mathfrak{F}^{\,X_1 + X_2} (t) \subseteq \mathfrak{F}^{\,|Y|} (t)
\subseteq   \mathfrak{F}^{\,Y} (t)\,$ for $\, 0 \le t < \infty\,$; and, as we remarked in (\ref{4.3.a}),  the second inclusion is strict  for $\, 0 < t < \infty\,$.   

The fourth and fifth  relations in (\ref{4.6}) follow  from the equations in
(\ref{2.27}) and $\, X_1 (\cdot) - X_2 (\cdot) = Y(\cdot)\,$; whereas the
sixth relation is a consequence of the strong solvability of the stochastic
equation (\ref{3.1}).
\end{proof}

The strictness  of the inclusion in $\, \mathfrak{F}^{\,(R_1, R_2)} (t)  =
\mathfrak{F}^{\,V} (t)   \subsetneqq  \mathfrak{F}^{\,W} (t)=
\mathfrak{F}^{\,(X_1, X_2)} (t)   \,$  for $\, 0 < t < \infty\,$ reflects the fact that some information
is inevitably lost  when one passes from the ``names" $\, (X_1 (\cdot) , X_2
(\cdot)) \,$ to the ``ranks"   $\, (R_1 (\cdot) , R_2 (\cdot)) \,$.  
In the case $\, \rho =1\,$, $\, \sigma =0\,$ that we are studying here, the
equations (\ref{2.12}), (\ref{2.13}) read 
\begin{align}
  \label{4.7}
  X_1 (t) \,&=\,  x_1 + \int_0^t \left(\, g\,\mathbf{ 1}_{ \{ X_1 (s) \le  X_2 (s) \} } - h\,\mathbf{ 1}_{ \{ X_1 (s) >  X_2 (s) \} } \right)   \d  s \,+\, \int_0^t\mathbf{ 1}_{ \{ X_1 (s) >  X_2 (s) \} }\, \d  W  (s)\,,\\
  \label{4.8} 
  X_2 (t) \,&=\,  x_2+ \int_0^t \left( \,g\,\mathbf{ 1}_{ \{ X_1 (s) >  X_2 (s) \} } - h\,\mathbf{ 1}_{ \{ X_1 (s) \le  X_2 (s) \} } \right)   \d  s \,-\, \int_0^t \mathbf{ 1}_{ \{ X_1 (s) \le   X_2 (s) \} }\, \d  W  (s)
\end{align}
with $\, W(\cdot) \equiv W_1 (\cdot)\,$ in accordance with (\ref{2.3}). The
construction (synthesis) of section \ref{sec3} shows that the solution
constructed there is {\it strong}, that is, adapted to the filtration $
\mathbf{ F}^{\,W}$ generated by the Brownian motion $\, W(\cdot)\,$ that
drives this system. Repeating the analysis of section \ref{sec2}, we see that
{\it any} solution $\, (X_1 (\cdot), X_2 (\cdot))\,$ of the system
(\ref{4.7}), (\ref{4.8}) can be cast in the form (\ref{2.27}), where
$\,Y(\cdot)\,$ is the pathwise unique, strong solution of the diffusion
equation (\ref{3.1}) with $\,    \mathbf{ F}^{\,Y} \equiv  \mathbf{ F}^{\,W}
\,$  and $\, L^Y (\cdot)\,$ is the local time  in (\ref{1.4}).  In particular,
any such solution of the system (\ref{4.7}), (\ref{4.8}) is a non-anticipative
functional  (given as in (\ref{2.27}))  of its driving Brownian motion $\,
W(\cdot)\,$. Thus,  pathwise uniqueness holds for the system of equations  
(\ref{4.7}), (\ref{4.8}).

On the other hand, and always in the case $\, \rho =1\,$, $\, \sigma =0\,$,
we can cast the equations (\ref{2.14}), (\ref{2.15}) as 
\begin{align}
  \label{4.9}
  X_1 (t) &=  x_1 + \int_0^t \left(\, g\,\mathbf{ 1}_{ \{ X_1 (s) \le  X_2 (s) \} } - h\,\mathbf{ 1}_{ \{ X_1 (s) >  X_2 (s) \} } \right)   \d  s \,+\, \int_0^t\mathbf{ 1}_{ \{ X_1 (s) >  X_2 (s) \} }\, \d  V  (s)\,,\\
  \label{4.10} 
  X_2 (t) &=  x_2+ \int_0^t \left( \,g\,\mathbf{ 1}_{ \{ X_1 (s) >  X_2 (s) \} } - h\,\mathbf{ 1}_{ \{ X_1 (s) \le  X_2 (s) \} } \right)   \d  s \,+\, \int_0^t \mathbf{ 1}_{ \{ X_1 (s) \le   X_2 (s) \} }\, \d  V  (s)
\end{align}
now driven by the Brownian motion $\, V(\cdot) \equiv V_1 (\cdot)\,$ in
accordance with (\ref{3.5.b}). We know from  Proposition \ref{Proposition 1}
that this system does {\it not} admit a strong solution. 

\smallskip
Collecting the results of this section and of the preceding one, we see that
we have established the following  theorem; this can be seen as a
two-dimensional analogue of %
\cite{MR1008231}.  

\begin{theorem}
  \label{Theorem 1} 
  The system of stochastic differential equations (\ref{1.1}), (\ref{1.2})
  with $\, \rho \, \sigma =0\,$ has a weak solution which is unique in the
  sense of the probability distribution. The same is true for each of the
  systems of equations (\ref{4.7}), (\ref{4.8}) and (\ref{4.9}),
  (\ref{4.10}).  
  
  On the other hand, the system of stochastic differential equations
  (\ref{4.7}), (\ref{4.8})   admits a     strong solution,  which is therefore pathwise unique;
  whereas the system  (\ref{4.9}), (\ref{4.10}) does not admit a strong
  solution.  
\end{theorem}

We have deduced here pathwise uniqueness from the   ``obverse  \textsc{Yamada
  \& Watanabe}" results of %
\cite{MR1128494} and \cite{MR1978664}.

\smallskip 
\begin{figure}[H]
  \centerline{ $$\boxed{\hbox{~~The following Figure is a Simulation of the
        Diffusion with $\, \rho =1\,$, $\, \sigma =0\,$ and $\, g = h
        =1\,$.~\,}   } $$ }
  \vspace{4ex}
  \begin{center}
    \includegraphics[width=84mm,keepaspectratio=true]{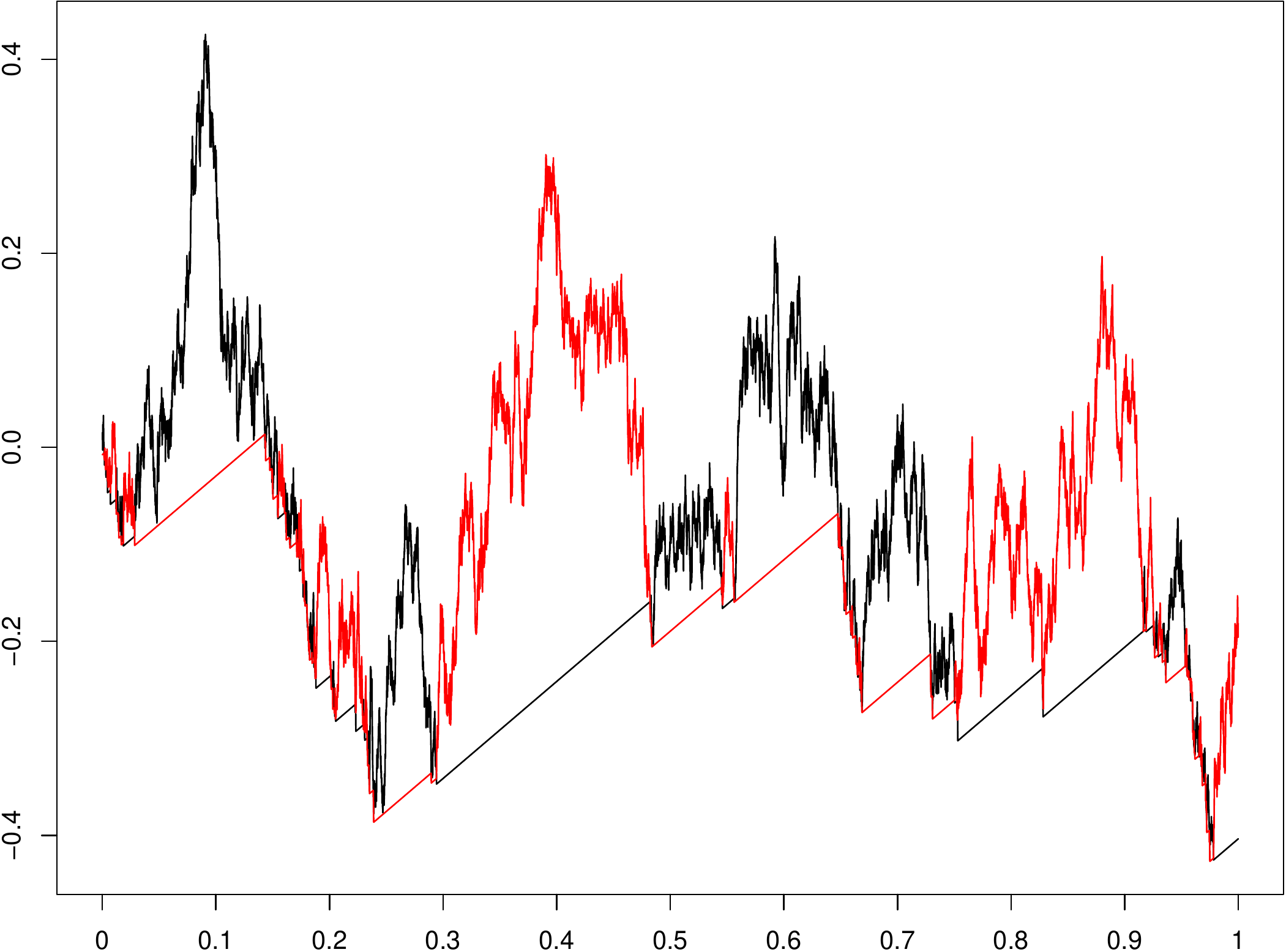}
    \caption{Simulated processes; Black $\,=X_1(\cdot)\,$, Red $\,=X_2 (\cdot)\,$. }
  \end{center}
\end{figure}

\subsection{ \textsc{The Non-Degenerate Case} }
\label{sec4.2}

Let us move now on to the non-degenerate case $\, \rho \, \sigma >0\,$. 
In  contrast to the  situation we encountered in Proposition
\ref{Proposition 1}, now both pairs $\, \big( X_1 (\cdot), X_2 (\cdot) \big)
\,$ (of positions by ``name")  and $\, \big( R_1 (\cdot), R_2 (\cdot) \big)
\,$ (of positions by ``rank") generate two-dimensional Brownian filtrations,
as our next result shows.

\smallskip
\begin{proposition}
  \label{Proposition 2}
  In the non-degenerate case $\, \rho \, \sigma >0\,$, we have for every $\,
  0 < t < \infty\,$   the filtration relations 
  \begin{multline}
    \label{4.11}
    \mathfrak{F}^{\,(V_1, V_2)} (t) \,=\,\mathfrak{F}^{\,(R_1,\, R_2)} (t)  \, =\,\mathfrak{F}^{\,(|Y|, V)} (t)  \, =\,\mathfrak{F}^{\,(|Y|, Q)} (t)
    \,\subsetneqq   \, \mathfrak{F}^{\,(X_1, X_2)}
    (t)\,=\,\mathfrak{F}^{\,(W_1, W_2)} (t)\\ =\,\mathfrak{F}^{\,(Y, V)} (t)  \,
    =\,\mathfrak{F}^{\,(W, U^{\flat})} (t) \,=\,\mathfrak{F}^{\,(Y, U^{\flat})} (t)\,=\,\mathfrak{F}^{\,(W, Q)} (t) \,=\,\mathfrak{F}^{\,(Y, Q)} (t)    
    \, ,
  \end{multline}
  where the inclusion is strict.  
\end{proposition}

\begin{proof}
  We have clearly $\, \mathfrak{F}^{\,L^Y} (t) \subseteq
  \mathfrak{F}^{\,V^{\flat}} (t) \subseteq \mathfrak{F}^{\,(V_1, V_2)} (t) \,$
  by virtue of (\ref{4.5}), (\ref{3.5.b}); back into the equations of
  (\ref{4.3}), (\ref{4.4}), this implies $\, \, \mathfrak{F}^{\,(R_1, \,R_2)}
  (t)\subseteq \mathfrak{F}^{\,(V_1, V_2)} (t) \,$. For the reverse inclusion,
  we note that   $\, \mathfrak{F}^{\,L^Y} (t) \subseteq  \mathfrak{F}^{\,|Y|}
  (t) \subseteq \mathfrak{F}^{\,(R_1, R_2)} (t) \,$ holds, thanks to
  (\ref{1.4})  and $\, \big| Y(\cdot) \big| =   R_1 (\cdot)- R_2 (\cdot)
  \,$. Back into (\ref{4.3}), (\ref{4.4})  this gives $\,
  \mathfrak{F}^{\,(V_1, V_2)} (t) \subseteq \mathfrak{F}^{\,(R_1, \,R_2)} (t)
  \,$, and the first equality in (\ref{4.11}) is proved.
  
  In conjunction with (\ref{2.6}) and (\ref{2.21}), these considerations also give $\, \mathfrak{F}^{\,(R_1, \,R_2)} (t) = \mathfrak{F}^{\,(|Y|, \,X_1+X_2)} (t)= \mathfrak{F}^{\,(|Y|, V)} (t) =\mathfrak{F}^{\,(|Y|, Q)} (t)\,$ for every $\, 0 \le t < \infty\,$, justifying the second and third equalities in (\ref{4.11}).

  For the fourth equality in (\ref{4.11}), we have noted already the inclusion
  $ \, \mathfrak{F}^{\,(X_1, \,X_2)} (t)\subseteq\mathfrak{F}^{\,(W_1, W_2)} (t)
  \,$ from the construction of (\ref{3.6})-(\ref{3.9}). In order to argue the
  reverse inclusion, let us note that (\ref{3.8}), (\ref{3.9}) imply  
  \begin{align*}   
    W_1 (t)\,&=\, \big( 1 / \rho \big)  \int_0^t  \left( \mathbf{ 1}_{ \{ X_1
        (s) >  X_2 (s) \} }\, \d  M_1 (s) - \mathbf{ 1}_{ \{ X_1 (s) \le
        X_2 (s) \} }\, \d  M_2 (s) \right)\,,\\
    W_2 (t)\,&=\, \big( 1 / \sigma \big)  \int_0^t  \left( \mathbf{ 1}_{ \{ X_1 (s) \le  X_2 (s) \} }\, \d  M_1 (s) - \mathbf{ 1}_{ \{ X_1 (s) >   X_2 (s) \} }\, \d  M_2 (s) \right)\,.
  \end{align*}
  On the strength of (\ref{3.6}) and  (\ref{3.7}), the $\, \mathbf{ F}^{(W_1,
    W_2)}-$martingales $\, M_1 (\cdot)\,$ and $\, M_2 (\cdot)\,$ are also
  martingales of  the filtration $\, \mathbf{F}^{\, (X_1, X_2)}$   generated
  by the state process $\,   (X_1  (\cdot) , \,   X_2  (\cdot)  ) \,$, so we
  deduce $\, \mathbf{F}^{\, (W_1, W_2)} \subseteq \mathbf{F}^{\, (X_1, X_2)}\,
  $ and the  fourth equality in (\ref{4.11}) follows.  

  \smallskip
  For the fifth equality we note that $\, \sigmab ( X_1 (t), X_2 (t)) =
  \sigmab ( Y (t), V (t))\,$   holds for every $\, t \ge 0\,$, courtesy of
  (\ref{3.10}), and gives $ \, \mathfrak{F}^{\,(X_1, X_2)} (t) =
  \mathfrak{F}^{\,(Y, V)} (t) \,$.  
  The sixth  equality in (\ref{4.11}) is a consequence  of  (\ref{2.20}),
  which gives actually the stronger statement   
  $\, \sigmab ( W_1 (t), W_2 (t)) = \sigmab ( W (t), U^{\flat} (t))\,$    for every $\, t \ge 0\,$.

  The seventh and eighth equalities in (\ref{4.11})   are  straightforward
  consequences of $\, \mathbf{ F}^W = \mathbf{ F}^Y\,$ and of the equation
  (\ref{2.22}), along with  its ``twin"  $\,U^{\flat}  (\cdot) = \int_0^\cdot
  \sign ( Y (t)  ) \, \d  Q  (t)\,$. The ninth equality is   now   obvious.  

 Finally, the strictness of the inclusion in (\ref{4.11}) follows now from the strictness of the inclusion in
$$
\mathfrak{F}^{\,(R_1,\, R_2)} (t)     \, =\,\mathfrak{F}^{\,(|Y|, Q)} (t)
    \,\subsetneqq   \, \mathfrak{F}^{\,(Y, Q)} (t)\,=\, \mathfrak{F}^{\,(X_1,\, X_2)} (t) \,, \qquad 0 < t < \infty\,,
$$
itself a consequence of the strictness of the inclusion in (\ref{4.3.a}) and of the independence of the processes $\, Y(\cdot)\,$, $\, Q(\cdot)\,$.
\end{proof}

Arguing as in the previous subsection  and using Proposition
\ref{Proposition 2},  we obtain the following   analogue of Theorem
\ref{Theorem 1} for the non-degenerate case.    

\begin{theorem}
  \label{Theorem 2} 
  The system of stochastic differential equations (\ref{1.1}), (\ref{1.2})
  with $\, \rho \, \sigma >0\,$ has a weak solution which is unique in the
  sense of the probability distribution. The same is true for each of the
  systems of equations (\ref{2.12}), (\ref{2.13}) and (\ref{2.14}),
  (\ref{2.15}). 
  
  On the other hand,   the system of stochastic differential equations
  (\ref{2.12}), (\ref{2.13})  admits a     strong solution,  which is therefore pathwise unique;     whereas the system   (\ref{2.14}), (\ref{2.15})  admits no strong    solution. 
\end{theorem}

The existence and uniqueness-in-distribution of a weak solution to
(\ref{1.1}), (\ref{1.2}) in the non-degenerate case  $\, \rho \, \sigma
>0\,$ follow also from the general, multidimensional results of 
\cite{MR917679}%
. The approach here is more direct and concrete,
capitalizing on the special two-dimensional nature of our setting; the
\textsc{Bass-Pardoux} results cannot, however, be applied to the  degenerate
case   $\, \rho \, \sigma =0\,$, so it does not seem possible to obtain even
the first part of Theorem \ref{Theorem 1} from them.  

We shall see in Theorem \ref{Theorem 3} of the next section that the system
of equations  (\ref{1.1}), (\ref{1.2}) admits a strong solution; this too
goes beyond the 
\cite{MR917679} %
 results. As far as we know
it is an open issue, whether such strength might obtain in their  setting as
well (to wit, where $\, \R^d\,$, $\,d \ge 2\,$ is partitioned into disjoint
polyhedral chambers, diffusion characteristics   are   constant in each
chamber, and strong non-degeneracy prevails).  

\subsection{ \textsc{Discussion} }
\label{sec4.4}

Theorem  \ref{Theorem 2}   shows  that the system of    (\ref{2.14}),
(\ref{2.15}) with $\, \rho \, \sigma >0\,$ can be thought of as a {\it
  genuinely  two-dimensional \textsc{Tanaka} example:} a system of
stochastic differential equations that admits a weak solution which is
unique in the sense of the probability distribution, but   no strong
solution. 

Theorem \ref{Theorem 2}    might usefully be compared also  with the
strictly one-dimensional results by 
\cite{MR1008231} and \cite{MR0326840}, as well as with Theorem 2 in
\cite{MR673915}; see also \cite{MR675177}.   
\cite{MR1008231} %
 shows, in particular, that with any given constants
$\, x \in \R\,$, $ \, \rho >0\,$ and $\, \sigma >0\,$, the one-dimensional
stochastic differential equation 
$$
X(t) \,=\, x + \int_0^t \, \taub \big( X(s) \big)\, \d  B(s)\,, \qquad 0 \le t < \infty\,,
$$
which has a pathwise unique, strong solution for  $\, \taub (\cdot) \,$ given as 
$\,   \taub_1   ( y ):=  \rho \,  \mathbf{ 1}_{  (0, \infty)    } (y) +
\sigma \, \mathbf{ 1}_{ (-\infty, 0] } (y)\,$ by virtue of the
\cite{MR0326840} %
results,  admits only a weak solution (and  no strong
solution) for $\, \taub (\cdot) \,$ equal to  
$$\,
\taub_2   ( y )\,=\,  \rho \,  \mathbf{ 1}_{  (0, \infty)    } (y) -  \sigma \, \mathbf{ 1}_{ (-\infty, 0] } (y)\,.
$$ 
In both these cases the variance structure is the same, namely $\, \taub_1^2
( y )= \taub_2^2   ( y )= \rho^2 \,  \mathbf{ 1}_{  (0, \infty)    } (y) +
\sigma^2 \, \mathbf{ 1}_{ (-\infty, 0] } (y)$. 

\smallskip
A rather similar situation obtains in Theorem \ref{Theorem 2} (and with
obvious minor chan\-ges, in Theorem \ref{Theorem 1} as well). To wit, both  systems of
equations (\ref{2.12}), (\ref{2.13}) and (\ref{2.14}), (\ref{2.15}) have the
same covariance structure, namely $\, \mathcal{A}(x_1, x_2) \,$ in
(\ref{Alpha}),   
as the original system (\ref{1.1}), (\ref{1.2}); whereas the  diffusion
matrices for (\ref{2.12}), (\ref{2.13}) and (\ref{2.14}), (\ref{2.15}) are
given  respectively by  
\begin{equation}
  \label{SigmaWV}
  \Sigma_W (x_1, x_2) \,=\,
  \begin{pmatrix}
     \rho  \,  \I{x_1 > x_2}   &   \sigma  \,\I{x_1 \le x_2 } \\
    - \rho  \,  \I{x_1 \le x_2}&  - \sigma \,\I{x_1 > x_2}
  \end{pmatrix}
  ,\quad   
  \Sigma_V (x_1, x_2) \,=\, \begin{pmatrix}
    \rho  \,  \I{x_1 > x_2 }   & \sigma\,\I{x_1 \le x_2}  \\
    \rho  \,  \I{x_1 \le x_2}  & \sigma\,\I{x_1 > x_2} 
  \end{pmatrix}.
\end{equation}
These two matrices  differ only by changes of  signs in the second row; yet
the system (\ref{2.12}), (\ref{2.13}) of stochastic differential equations
induced by the first of them is strongly solvable, whereas the system
(\ref{2.14}), (\ref{2.15}) induced by the second diffusion matrix is
solvable only   weakly. To paraphrase 
\cite{MR673915}%
, page
448, in this latter case the diffusion generated by the system (\ref{2.14}),
(\ref{2.15}) ``cannot leave the diagonal  $\, \{ x_1 = x_2\}\,$ strongly".

We shall see in the next section that the original system (\ref{1.1}),
(\ref{1.2}), with diffusion matrix of the form 
\begin{equation}
  \label{SigmaB}
  \Sigma_B (x_1, x_2) \,=\,
  \begin{pmatrix}
    ~ \rho  \,  \mathbf{ 1}_{  \{ x_1 > x_2   \} } + \sigma  \,  \mathbf{ 1}_{  \{ x_1 \le x_2   \} }  &   ~  0 ~     \\
    ~0~      &    \rho  \,  \mathbf{ 1}_{  \{ x_1 \le x_2   \} }+ \sigma  \,  \mathbf{ 1}_{  \{ x_1 > x_2   \} } ~
  \end{pmatrix}\,,
\end{equation}
is also strongly solvable. A thorough discussion of weak and strong
solutions, for all possible real square roots of the covariance matrix $\,
\mathcal{A}(x_1, x_2) \,$ in (\ref{Alpha}), appears in    subsection
\ref{sec9}.

\section{ \textsc{Strength, Restored} }
\label{sec4.3}

In conjunction with the identity $ \, \mathfrak{F}^{\,(X_1, X_2)} (t) =
\mathfrak{F}^{\,(W_1, W_2)} (t) \,$ from (\ref{4.11}), the expressions in
(\ref{3.12}), (\ref{3.13}) imply  
\begin{equation}
  \label{4.12}
  \, \mathfrak{F}^{\,(B_1, B_2)} (t) \, \subseteq \, \mathfrak{F}^{\,(X_1,
    X_2)} (t) \,,\qquad \forall ~~ ~0 \le t < \infty  
\end{equation}
in the non-degenerate case $\, \rho \, \sigma >0\,$.  It is an interesting
question to settle, whether   the reverse inclusion might also hold in
(\ref{4.12}), thus implying  the strength of the solution to the system
(\ref{1.1}), (\ref{1.2}) constructed in  section \ref{sec3}.

We do have such strength in the isotropic case   $\, \rho^2 = \sigma^2 =
1/2\,$ of equal variances, as has been well known since
\cite{MR532447,MR568986,MR673915}. To obtain this   property from
first principles in our context, let us observe from (\ref{2.3})-(\ref{2.5})
the representation $\, W(t) = (B_1 (t) - B_2 (t))/ \sqrt{2\,}\,$ in the
isotropic case, thus also the inclusion  
\begin{equation}
  \label{4.13}
  \mathfrak{F}^{\,Y} (t)\, =\, \mathfrak{F}^{\, W } (t) \, \subseteq  \,\mathfrak{F}^{\,(B_1, B_2)} (t) 
\end{equation}
for each $\, t \ge 0\,$; back into (\ref{2.4}), (\ref{2.5}) and in
conjunction with (\ref{4.11}), this gives $ \, \mathfrak{F}^{\,(X_1, X_2)}
(t) = \mathfrak{F}^{\,(W_1, W_2)} (t) \, \subseteq \, \mathfrak{F}^{\,(B_1,
  B_2)} (t)  \,$, showing that    (\ref{4.12}) holds with equality  in the
isotropic case.

More generally, the inclusion in (\ref{4.13}) will follow  whenever one is
able  to   show  that the {\it modified \textsc{Tanaka} equation} 
\begin{multline}
  \label{Tanaka}
  Y (t) \,=\, y - \lambda   \int_0^t \sign \big( Y(s) \big)\, \d  s
  \, +  \int_0^t   \, \Big(       \, \rho \, \mathbf{ 1}_{ \{ Y (s) > 0 \} }
  + \sigma \,  \mathbf{ 1}_{ \{ Y (s) \le 0 \} }   \,\Big)\,\d  B_1
  (s)\\
  -\int_0^t  \Big( \,  \rho \, \mathbf{ 1}_{ \{ Y (s) \le  0 \} }   + \sigma \,
  \mathbf{ 1}_{ \{ Y (s) > 0 \} }   \Big)\,   \d  B_2 (s)\,,\qquad 0
  \le t < \infty\,, 
\end{multline}
a consequence of (\ref{3.1})-(\ref{2.5}),  has a strong (that is, $\,
\mathbf{ F}^{\,(B_1, B_2)}-$adapted)  solution; 
then   {\it this strong solvability implies the inclusion  
  \begin{equation}
    \label{4.12.a}
    \, \mathfrak{F}^{\,(X_1, X_2)} (t) \, \subseteq \, \mathfrak{F}^{\,(B_1,
      B_2)} (t) \,,\qquad \forall ~~ ~0 \le t < \infty \,, 
  \end{equation}
  that is, the strong solvability of the system (\ref{1.1}), (\ref{1.2}). } 
Indeed, as 
  Johannes \textsc{Ruf} suggests, the strong solvability of (\ref{Tanaka}) implies the inclusion $\,
  \mathfrak{F}^{\,(Y, B_1, B_2)} (t)    \subseteq   \mathfrak{F}^{\,(B_1, B_2)} (t)\,$, and   (\ref{4.12.a}) follows then directly from the equations (\ref{1.1}), (\ref{1.2}).

\begin{theorem}
  \label{Theorem 3} 
  The system of stochastic differential equations (\ref{1.1}), (\ref{1.2})
  admits a pathwise unique, strong solution; in particular, the inclusion
  (\ref{4.12.a}) holds.    
\end{theorem}

\begin{proof}    
  We have seen already that  the weak solution of (\ref{1.1}),
  (\ref{1.2})     constructed in section \ref{sec3} is actually strong in the
  isotropic case   $\, \rho^2 = \sigma^2 =  1/2\,$ of equal variances; so
  we   focus   on the   case $\, \rho \neq \sigma\,$.  
  
   From the   results of \textsc{Yamada \& Watanabe} (e.g.,
    p. 310 in %
    \cite{MR1121940}%
    ) and their ``obverse"   counterparts due to 
    \cite{MR1128494} %
  and %
  \cite{MR1978664} we know that, in the
    presence of uniqueness in distribution, strong existence   and
    pathwise uniqueness are equivalent; and in this case every solution is
    strong. We have already seen the that the solution of \eqref{1.1},
    \eqref{1.2} is unique in distribution, so   we need only  
     show that there is a strong solution; and by the argument
    before the statement of the Theorem, it is enough to prove that the solution    of \eqref{Tanaka} is pathwise unique. Furthermore, it suffices to argue
    such pathwise uniqueness on an arbitrary but fixed time-horizon $[0, T ]$
    of finite length $T \in (0,\infty)$.

  A further reduction is  that we need to prove  such pathwise uniqueness only
  as regards the  {\it ``driftless" modified \textsc{Tanaka} equation}  
  \begin{multline}
    \label{TanakaThree}
    Y (t) \,=\, y   +  \int_0^t   \, \Big(       \, \rho \, \mathbf{ 1}_{ \{ Y
      (s) > 0 \} }   +    \sigma \,  \mathbf{ 1}_{ \{ Y (s) \le 0 \} }
    \,\Big)\,\d  \betab_1 (s)         \\- 
    \int_0^t  \Big( \,  \rho \, \mathbf{ 1}_{ \{ Y (s) \le  0 \} }   + \sigma
    \,  \mathbf{ 1}_{ \{ Y (s) > 0 \} }   \Big)\,   \d  \betab_2 (s)~~~ 
  \end{multline}
  for $\, 0 \le t \le T \,$, with $\, \betab_1 (\cdot)\,$ and $\, \betab_2
  (\cdot)\,$ independent, standard Brownian motions under a suitable
  equivalent probability measure $\, \mathbb{Q}\,$; because then a
  \textsc{Cameron-Martin-Girsanov}   change of   measure brings us back to
  (\ref{Tanaka}),  with  
  $$
  B_1 (t) \,= \, \betab_1 (t) + { \lambda\, t \over \, \rho - \sigma\,}\,\,,
  \quad ~~ B_2 (t) \, = \, \betab_2 (t) + { \lambda\, t \over \, \rho -
    \sigma\,}\,\,, \qquad 0 \le t \le T
  $$
  independent Brownian motions under the original probability measure $\,
  \mathbb{P}\,$.  Let us observe also, that uniqueness in distribution holds
  for the equation \eqref{TanakaThree}: every solution is a standard Brownian
  motion starting  at $\, Y(0)=y\,$, under the auxiliary, equivalent
  probability measure $\,\mathbb{Q}\,$.

  Introducing the independent, standard $\,\mathbb{Q}-$Brownian motions 
  $$   \betab (\cdot)\,
  := \,  \frac{\, \betab_1 (\cdot) + \betab_2 (\cdot)  \,} { \sqrt{2\,}}\, ,\qquad 
  \thetab (\cdot)\,
  := \,  \frac{\,\betab_1 (\cdot) - \betab_2 (\cdot)  \,} { \sqrt{2\,}}\, ,
  $$
 we can write (\ref{TanakaThree})  in the  equivalent form 
  $$
  \, Y (t) \,=\, y  + \frac{\,\rho    - \sigma\,}{\sqrt{2\,}}  \int_0^t   \mathrm{sgn}\big(  Y (s)  \big)\,  \d   \betab  (s)    
 -   \frac{\,\rho    + \sigma\,}{\sqrt{2\,}} \, \thetab   (t)\,.   
  $$
  Pathwise uniqueness for this equation follows now directly from 
  Theorem \ref{Theorem 0}.  
  \end{proof}

\subsection{ \textsc{The Algebra and Geometry of Strength}  }
\label{sec9}

The matrices of (\ref{SigmaWV}), (\ref{SigmaB}) are square roots of the
covariance matrix $\, \mathcal{A} (x_1, x_2)\,$  in (\ref{Alpha}). The
general real square root of this covariance matrix is of the form  
\begin{multline}
  \label{SO2}
   \mathcal{A}^{1/2} (x_1, x_2) \,=\, \Sigma (x_{1}, x_{2}) \,= \,
   \begin{pmatrix} \rho & 0 \\ 0 & \sigma  \end{pmatrix}
    \begin{pmatrix} 1 & 0 \\ 0 & \varepsilon \end{pmatrix}
    \begin{pmatrix} 
      \cos \varphi & - \sin \varphi \\
      \sin \varphi & \cos \varphi 
    \end{pmatrix}\, \I{x_{1} > x_{2}} \\
    +\begin{pmatrix} \sigma & 0 \\ 0 & \rho   \end{pmatrix} 
    \begin{pmatrix}  1 & 0      \\ 0 & \delta \end{pmatrix}
    \begin{pmatrix} 
      \cos \vartheta & - \sin \vartheta \\
      \sin \vartheta & \cos \vartheta
    \end{pmatrix} \,\I{x_{1} \le x_{2}} \, , 
\end{multline}
where $\, \varepsilon \in \{ \, -1, 1\, \}\, $,  
$\,  \delta \in \{ \, -1, 1\, \} \, $, $\, 0 \le \varphi, \vartheta  \le 2
\pi\, $, or equivalently 
\begin{equation}
  \label{SO3}
  {\cA}^{1/2} (x_1, x_2) \,=\, \Sigma (x_{1}, x_{2}) \,= \,
  \Sigma_+\, \I{x_{1} > x_{2}} + \Sigma_-\, \I{x_{1} \le  x_{2}} 
\end{equation}
in the notation
of (\ref{Sigma1234}). We also write the drift term of the generator
\eqref{0.1} in vector  form
\begin{displaymath}
  \mathbf{G}(x_1,x_2)\,=\,
  \begin{pmatrix}g\\-h\end{pmatrix}\I{x_1\leq x_2}+
  \begin{pmatrix}-h\\g\end{pmatrix}\I{x_1>x_2}
\end{displaymath}
Each  choice   in (\ref{SO2}) (or equivalently, (\ref{SO3})) leads to a system  of stochastic differential equations with rank-based characteristics   
\begin{equation}
  \label{eq:xSigma}
  \d \bX(t)\,=\,\cA^{1/2} (\bX(t))\,\d\bU(t)+ \mathbf{G}(\bX(t))\,\d t\,,
\end{equation}
driven by a two-dimensional Brownian motion $\,\bU(\cdot) = (U_1(\cdot), \,U_2(\cdot))^\prime\,$. This system admits a weak solution, which is  unique in the sense of the probability distribution.

\begin{theorem}
  \label{Theorem 5} 
  With the notation of (\ref{Sigma1234}) and $\be_1 = (1, 0)'$,
  $\be_2 = (0, 1)'$,  the (unique in distribution) weak
  solution   of  the system of stochastic differential equations
  (\ref{eq:xSigma}) fails to be strong, if and only if
  \begin{equation}
    \label{eq:IP-1}
    ( \be_1 - \be_2)'\Sigma_-\,=\,- ( \be_1 - \be_2)'\Sigma_+\,.
  \end{equation}  
\end{theorem}

\begin{proof}
  We note that the equation satisfied by the difference
  $\,Y(\cdot)=X_1(\cdot) - X_2(\cdot) \,$ is 
  \begin{equation}
    \label{Y}
    \begin{aligned}[b]
      \d Y(t)\,&=\,\frac{(\be_1-\be_2)'}2 \,\Big(\big(\Sigma_++\Sigma_-\big)+
      \sign(Y(t)) \big(\Sigma_+-\Sigma_-\big)\Big)\d\bU(t)- \lambda\,
      \sign(Y(t))\,\d t\\
      \,&=\,- \lambda\, \sign(Y(t))\,\d t \,+  
      \Big( \big( \be_1 - \be_2 \big)'   \,\Sigma_+  \,
      \mathbf{ 1}_{\{ Y(t) >0\}}+ \big( \be_1 - \be_2 \big)' \, \Sigma_-  \, 
      \mathbf{ 1}_{\{ Y(t) \le0\}} \Big) \,\d\bU(t) \\
      \,&=\,\sign(Y(t))\,\d M(t)+\d N(t) \,,
    \end{aligned}
  \end{equation}
  where we have set 
  \begin{displaymath}
    M(t)\,:=\,\frac{(\be_1-\be_2)'}2\big(\Sigma_+-\Sigma_-\big)\bU(t)-\lambda
    t\, \quad\text{and}\quad 
    N(t)\,:=\,\frac{(\be_1-\be_2)'}2\big(\Sigma_++\Sigma_-\big)\bU(t)\,.
  \end{displaymath}
  Elementary calculation shows  that $\q{M,N}(\cdot)=0$, since
  $(\be_1-\be_2)'\Sigma_{\pm}$ are two unit vectors on the plane.
  
  \smallskip
  We observe next that, because the indicators  $\I{X_1(t) \ge
    X_2(t)}$ and $\I{X_1(t) < X_2(t)}$ are functions of
  the difference $Y( t) = X_1( t) - X_2( t)$,   the $\sigma-$algebra
  $\,{\mathfrak F}^{\,(X_1 , X_2)} (t)\,$ is contained in the
  $\sigma$--algebra   
  $\,{\mathfrak F}^{\,(U_1, \,U_2, Y)} (t)\,$, for every $t \in [0,
  \infty)$, by construction as in \eqref{eq:xSigma}.  
  When the    equation (\ref{Y}) for $\,Y(\cdot)\,$ is solvable strongly with respect to the planar Brownian  motion $\bU(\cdot) = (U_1(\cdot), \,U_2(\cdot))^\prime$, we have the
  filtration  comparisons  
  $$
  {\mathfrak F}^{\,X_1 - X_2} (t) = {\mathfrak F}^{\,Y} (t) 
  \subseteq {\mathfrak F}^{\,(U_1, \,U_2)} (t),  
  \quad \text{thus}   \quad
  {\mathfrak F}^{\,(X_1, X_2)} (t)\subseteq{\mathfrak F}^{\,(U_1, U_2,
    Y)}(t) = {\mathfrak F}^{\,(U_1,\, U_2)} (t): 
  $$ 
  the system \eqref{eq:xSigma} is then strongly solvable.   When the equation (\ref{Y}) for $\, Y(\cdot)\,$ admits no strong solution with
  respect to $\, \mathbf{ U} (\cdot)\, $,  the system (\ref{eq:xSigma}) admits no
  strong solution. 
  There are now three possibilities.

  \begin{enumerate}
  \item When $\,(\be_1-\be_2)'\Sigma_+=(\be_1-\be_2)'\Sigma_- \,$,
    the second line of the     equation (\ref{Y}) for $\,Y(\cdot)\,$ simplifies to \eqref{3.1} which,    as we have already remarked in subsection \ref{sec2.3},  has a pathwise    unique, strong solution;  so the system \eqref{eq:xSigma} is then strongly    solvable. 
    
        \smallskip
  \item If both $(\be_1-\be_2)'(\Sigma_+\pm\Sigma_-)$ are non-zero vectors,
    we can apply the \textsc{Cameron-Martin-Girsanov} theorem as in the proof
    of Theorem \ref{Theorem 3}. 
    To wit, we can restrict ourselves to a finite time horizon $[0,T]$ 
    and  define a suitable equivalent probability measure $\PQ$ under which
    $N(\cdot)$ and $M(\cdot)$ are continuous, strongly orthogonal
    martingales. In this case the 
    quadratic variations $\q{M}(\cdot) $ and $\q{N}(\cdot)$ are proportional,
    so the domination 
    condition of Theorem \ref{Theorem 0} also holds. Then $Y(\cdot)$ is
    adapted to 
    the filtration of $(M(\cdot),N(\cdot))$, and this proves that the solution
    to \eqref{eq:xSigma} is strong in this case as well.
    
    \smallskip
  \item Finally, it is possible that $\,(\be_1-\be_2)'(\Sigma_+ +\Sigma_-)\,$ is
    zero; this condition is formulated in the statement as \eqref{eq:IP-1}. 
    Then the perturbation $\, N(\cdot)\,$ vanishes  in (\ref{Y}), which
    becomes then an ordinary {\it drifted \textsc{Tanaka}   equation}
    \begin{equation}
     \label{TanakaDrift}
      \d Y (t)\,=\, \mathrm{sgn} \big( Y(t) \big)\, \d M(t)\,,\quad  M(t)\,=\,
      (\be_1-\be_2)'   \Sigma_+ \bU(t)-\lambda    t\,   
    \end{equation}
    driven by the Brownian motion $\, M(\cdot)\,$ with negative drift.  From   the \textsc{Tanaka} formula,  the theory of the \textsc{Skorokhod} reflection problem, and (\ref{4.3.a}), we note 
    $$
    M(t) = |Y(t)| - |y| - 2\, L^Y (t)\,, \qquad \mathrm{thus} \qquad
    \mathfrak{ F}^{\,M} (t) \, = \,\mathfrak{ F}^{\,|Y|} (t) \,\subsetneqq \,\mathfrak{F}^{Y}(t)\,, \quad 0    < t <\infty \,.
    $$    
    We consider also the independent, one-dimensional standard Brownian motions
$$
U_+ (\cdot) \,:=\, (\be_1-\be_2)'\Sigma_+ \mathbf{ U}(\cdot)\, , \qquad \mathcal{Z}(\cdot) \,:=\, \nub' \, \mathbf{ U}(\cdot)\,,
$$    
where the unit vectors $\, \Sigma_+' (\be_1-\be_2)\,$ and $\, \nub\,$ are orthogonal, and note 
$\,
\mathfrak{ F}^{\,\mathbf{ U}} (t) \,=\, \mathfrak{ F}^{\,  U_+} (t)   \vee \mathfrak{ F}^{\,\mathcal{ Z}} (t)  \,=\, \mathfrak{ F}^{\,M} (t) \vee \mathfrak{ F}^{\,\mathcal{ Z}} (t)\,, \quad 0 \le t < \infty\,$.

\smallskip

The     equation (\ref{TanakaDrift}) has no strong solution; that is, $\,
\mathfrak{ F}^{\,Y} (t)    \subseteq \mathfrak{ F}^{\,\mathbf{ U}} (t)  $
cannot possibly hold for $\, 0< t < \infty\,$. For if it did, then the
process $$\, W (\cdot) \,:= \, Y(\cdot) - y + \lambda \int_0^\cdot
\mathrm{sgn} (Y(t))\, \mathrm{d} t\,,$$ which generates exactly the same
filtration as the process $\, Y(\cdot)\,$, would be a standard
one-dimensional Brownian motion and  independent of $\,\mathcal{Z}(\cdot)\,$,
yet adapted to the filtration generated by $\, \mathbf{ U} (\cdot)\,$; but
this is   impossible, since $\, \mathfrak{ F}^{\,M} (t) = \mathfrak{
  F}^{\,|Y|} (t) \subsetneqq  \mathfrak{F}^{Y}(t)= \mathfrak{ F}^{\,W} (t)
\indep \mathfrak{ F}^{\,\mathcal{Z}} (t) \,$ for $\, 0 < t < \infty\,$ would
lead  then to  
$$ 
\mathfrak{ F}^{\,\mathbf{ U}} (t) \,=\, \mathfrak{ F}^{\,  U_+} (t)   \vee
\mathfrak{ F}^{\,\mathcal{ Z}} (t)  \,=\, \mathfrak{ F}^{\,M} (t) \vee
\mathfrak{ F}^{\,\mathcal{ Z}} (t)\, \subsetneqq\,  \mathfrak{ F}^{\,W} (t)
\vee \mathfrak{ F}^{\,\mathcal{ Z}} (t) \subseteq \,\mathfrak{ F}^{\,\mathbf{
    U}} (t)  \,. 
$$ 
  \end{enumerate}
  To summarize: the solution to the system \eqref{eq:xSigma} fails to be
  strong if, and only if, \eqref{eq:IP-1} holds. 
\end{proof}
\begin{remark}
  Under the condition (\ref{eq:IP-1}), the diffusion vectors $\, \mathbf{ s}_+
  := ( \be_1 
  - \be_2 )  ' \,\Sigma_+\,$ (for the half-plane $\{ x_1 >x_2\}$) and $\,
  \mathbf{ s}_- :=  ( \be_1 - \be_2  )'\,\Sigma_-\,$ (for  the half-plane
  $  \{ x_1 \le x_2\}$) in the equation (\ref{Y}) point in {\it exactly opposite directions}. As a
  result,   it becomes impossible for the 
  planar diffusion $\, \mathbf{ X} (\cdot) = (X_1 (\cdot), X_2(\cdot))\,$ of
  \eqref{eq:xSigma} to ``escape strongly from the diagonal
  $\,\{x_1=x_2\}\,$''.   
  
  Let us also note that $\,\mathbf{ s}_{\pm}\,$ are unit vectors, so the condition \eqref{eq:IP-1}
  holds if and only if their scalar product is $\,-1$. Working out this
  relation using the representation given in \eqref{Sigma1234}, we obtain the condition     \eqref{WEAK}. 
\end{remark}

\begin{remark}
  Since $\, (\sigma^{2} \varepsilon  - \rho^{2} \delta)^{2} +
  \rho^{2}\sigma^{2} (1 + \varepsilon \delta)^{2} = 1\, $, the condition of
  (\ref{WEAK})   is equivalent to  
  \begin{equation}
    \label{geom}
    \pi + \vartheta  \,  =\, \varphi + \psi
  \end{equation}
  (modulo $2\, \pi\,$), which involves only the angles $\, \vartheta  \,$
  and $\, \varphi \,$ of (\ref{SO2}), as well as the angle $\, \psi \in
  (-\pi,  \pi]\, $  
  determined from $\, \cos \psi = \rho \sigma ( 1+ \varepsilon \delta)\, $, 
  $\, \sin \psi = \sigma^{2} \varepsilon - \rho^{2} \delta\, $.  
\end{remark}

\begin{example}
  \label{52}
  The system  of (\ref{1.1}), (\ref{1.2}),   with diffusion matrix $\,
  \Sigma_B\,$ in (\ref{SigmaB}), corresponds to $\, \varepsilon = \delta
  =1\,$, $\, \varphi = \vartheta = 0\,$ and thus to $\, \psi\,$ that has to
  satisfy $\, \cos \psi = 2 \, \rho \, \sigma \ge 0\,$, $\, \sin \psi =
  \sigma^2 - \rho^2\,$; there is no way that (\ref{geom}) can hold,   so the
  system \eqref{eq:xSigma} is strongly solvable for all choices of $\, \rho\,$, $\sigma\,$ as
  in (\ref{2.0}).  
  
  Similarly, the system  of (\ref{2.14}), (\ref{2.15})  with diffusion
  matrix $\, \Sigma_W\,$ in (\ref{SigmaWV}) corresponds to $\, \varepsilon =
  -1\,$, $\delta =1\,$, $\, \varphi =0\,$, $\vartheta = - \pi / 2\,$, and
  thus $\,\psi = - \pi / 2\,$ as well. Once again there is no way for
  (\ref{geom}) to hold,   so the system \eqref{eq:xSigma} is always strongly
  solvable.  
  
  Whereas the system  of (\ref{2.12}), (\ref{2.13})  with diffusion matrix
  $\, \Sigma_V\,$ in (\ref{SigmaWV})  corresponds to $\, \varepsilon =
  1\,$, $\delta =-1\,$, $\, \varphi =0\,$, $\vartheta = - \pi / 2\,$, so
  $\,\psi =   \pi / 2\,$ for all choices of $\, \rho\,$, $\sigma\,$ as in
  (\ref{2.0}). In this case  (\ref{geom}) always holds,   and the system
  \eqref{eq:xSigma} is never strongly solvable. 
 \end{example}

\begin{example}
  \label{51}
  The  matrix $\, \mathcal{A} (x_1, x_2)\,$  in (\ref{Alpha}) has a total of
  64 square roots    of the form $\, \Sigma_1 \, \mathbf{ 1}_{\{ x_1 >x_2\}}
  + \Sigma_2 \, \mathbf{ 1}_{\{ x_1 \le x_2\}}\,$ with
  \[
  \Sigma_1 \,\in\, \left\{ 
    \begin{pmatrix}
      \pm \rho   &   0  \\
      0   &  \pm \sigma
    \end{pmatrix}\,,\,\, 
    \begin{pmatrix}
      0   &   \pm \rho   \\
      \pm \sigma   &  0
    \end{pmatrix}    
  \right\}
  \quad \text{and} \quad 
  \Sigma_2 \,\in\, \left\{ 
    \begin{pmatrix}
      \pm \sigma     &   0  \\
      0   &  \pm \rho
    \end{pmatrix}\,,\,\, 
    \begin{pmatrix}
      0   &   \pm \sigma   \\
      \pm \rho   &  0
    \end{pmatrix}
  \right\} \,.
  \]
  Among these,   48   lead to strongly solvable systems in the isotropic
  ($\,\rho = \sigma = 1 / \sqrt{2\,}\,$)  or degenerate ($\,\rho \, \sigma =
  0\,$) cases, whereas 56 choices  lead to strongly solvable systems in all
  other cases. 
\end{example}

\section{ \textsc{Joint Distributions} }
\label{sec5}

The representations (\ref{2.23}), (\ref{2.24}) involve the triple $\, ( Y^+
(t), Y^- (t), L^Y (t))\,$ as well as the independent random variable $\,
Q(t)\,$, where $\, Y(\cdot)\,$ is the diffusion  process in (\ref{3.1}), $\,
L^Y (\cdot)\,$ the local time of this process at the origin as in
(\ref{1.4}), and $\, Q(\cdot)\,$ an independent, standard Brownian
motion. Thus, in order to compute  the joint distribution of $\, (X_1 (t),
X_2 (t) )\,$ via (\ref{2.23}), (\ref{2.24}), 
we need first to find that of  
$\, ( Y^+ (t), Y^- (t), L^Y (t))\,$.  

\smallskip
In order to do this, we consider  the  ``reference probability measure" $\,
\mathbb{P}_\star\,$, under which the process $\, Y(\cdot)\,$ becomes
standard Brownian motion. According to (\ref{3.1}) and  the
\textsc{Girsanov} theorem (e.g., \cite{MR1121940}, section 3.5), we have for
every $\, t \in [0, \infty)\,$ the \textsc{Radon-Nikod\'ym} derivative 

\begin{multline} \label{5.1}
  \left.\frac{\d\mathbb{P}}{\d\mathbb{P}_\star}\right|_{\mathfrak{F}^{\,Y}(t)} = 
  \exp\cbr{-\lambda \int_0^t \sign \big(Y(s)\big) \d Y(s) - 
    \frac{\lambda^2}2\int_0^t \sign^2 \big( Y(s) \big) \d  s} \\
  =\exp\cbr{\lambda \Big( |y| - \big| Y(t) \big| + 2 L^Y (t) \Big) - \frac{\lambda^2}2 t},
\end{multline}
thanks to the \textsc{Tanaka} formulas once again.  Thus, for any
Borel subsets $\, A\,$, $\,B\,$ of $\, [0, \infty)\,$ and $\, \Theta \,$ of
$\, \R\,$, and with $\, \mathbb{E}_\star\,$ denoting expectation with
respect to the reference probability measure, we have  
\begin{multline}
  \label{5.2}
  \P{ Y^\pm (t) \in A\,, ~ Y^\mp (t) =0\,, ~ 2\, L^Y (t) \in B
  \,, ~Q (t) \in \Theta}=
  \exp\cbr{ \lambda  |y| - \frac{\lambda^2}2 t} \times %
  \\
  \E[\star]\br{\exp\cbr{\lambda\big( 2 L^Y (t) -  Y^\pm (t)\big)} 
    \I{Y^\pm (t) \in A\,, \, \,Y^\mp (t) =0\,, \, \,2\, L^Y (t) \in
      B}}\int_\Theta\frac{\exp\cbr{ - \frac{\vartheta^2}{2t}}}{\sqrt{\,2\, \pi\, t\,\,}} \d  \vartheta.
\end{multline}

\begin{remark}
  \label{Remark 3}
    We have the classical result $\,\, \mathbb{P}_\star \big(\,
    \mathrm{Leb} \{\, 0 \le t \le T : Y(t) =0 \} >0 \big) =0\,$ for the
    Lebesgue measure of the    Brownian zero-set, thus $\, \mathbb{P}
    \big(\, \mathrm{Leb} \{\, 0 \le t \le T : Y(t) =0 \} >0 \big) =0\,$
    holds from (\ref{5.1}) for every $\, T \in (0, \infty)\,$. We deduce
    $\, \mathbb{P}  \big(\, \mathrm{Leb} \{\, 0 \le t < \infty \,:\, X_1(t)
    =X_2 (t) \} =0 \big) \,=\,1\,$.

    Let us also recall that the transition probability density function $\,
    \mathbb{P} \big( Y(t) \in \d  \xi\, \big|\, Y(0) = y \big)\,=\, \mathfrak{p}_t (y, \xi)\,\d  \xi \,$ for the one-dimensional ``bang-bang" diffusion process $\, Y(\cdot)\,$ of (\ref{3.1}) is given by
    \begin{equation}
      \label{3.2}
      \mathfrak{p}_t (y, \xi)\,=\, { 1 \over \, \sqrt{2 \pi t\,} \,} \left( \exp \left\{ - { (y - \xi - \lambda\, t)^2\, \over 2\, t} \right\} + \lambda \, e^{\, -2 \lambda \xi} \int_{y+\xi}^\infty e^{\, - (u- \lambda t)^2/  (2 t) } \, \d  u \right) 
    \end{equation}
    for $\, \xi > 0\,$, and by 
    \begin{equation}
      \label{3.3}
      \mathfrak{p}_t (y, \xi)= \frac1{\sqrt{2 \pi t}} \left( \exp\cbr{2
          \lambda y - \frac{(y - \xi + \lambda\, t)^2}{2t}} + 
        \lambda \, e^{2 \lambda \xi} \int_{y-\xi}^\infty e^{ - (u- \lambda t)^2/  (2 t) } \, \d  u \right) 
    \end{equation}
    for $\, \xi  \le 0\,$. In particular, with $\, y=0\,$ the function 
    \begin{multline}
      \label{3.4}
      \xi \, \mapsto \, \mathfrak{p}_t (  \xi)\,\equiv \, \mathfrak{p}_t (0,
      \xi) =\\ \frac1{\sqrt{2 \pi t}} \left( \exp\cbr{ - { (| \xi | + \lambda\,
            t)^2\, \over 2\, t}} + 
        \lambda \, e^{\, -2 \lambda |\xi |} \int_{| \xi |}^\infty e^{\, - (u- \lambda t)^2/(2 t) } \, \d  u \right)\,
    \end{multline}
    is evenly symmetric about the origin. Similar formulas hold for $\, y
    <0\,$; for the details of these computations, see 
    \cite{MR744236}%
    .  
\end{remark}

\subsection{ \textsc{The Isotropic  Case with} $~\, y = x_1 -  x_2 \ge 0\,$}
\label{sec5.1}

This   equal variance   case $\, \rho = \sigma = 1 / \sqrt{2\,}\,$ affords
the most straightforward computation: from the representations of
(\ref{2.26})  and the independence of $\, Y(t)\,$ and $\, Q(t)\,$, we obtain  
\begin{equation}
  \label{5.3}
  \mathbb{P} \big( X_1 (t) \in \d  \xi_1\,, \, X_2 (t) \in \d  \xi_2 \big) \,=\,   { \,   \mathfrak{ p}_t \big(y,  \xi_1 - \xi_2\big)\,  \over\, 2\, \sqrt{\, 2\, \pi\, t\,}\,} \, \exp \left\{ - { \,(\xi_1 + \xi_2 - z - \nu \, t )^2\, \over 2 \, t} \right\} \,\d  \xi_1 \, \d  \xi_2 
\end{equation}
for $\, (\xi_1, \xi_2) \in \R^2\,$ in the notation of (\ref{2.1}),
(\ref{3.2}) and (\ref{3.4}). The resulting transition probability density
function  
$$
\mathfrak{P}_{\,t} (\xi_1, \xi_2) \,=\, { \,   \mathfrak{ p}_t \big(y,  \xi_1 - \xi_2\big)\,  \over\, 2\, \sqrt{\, 2\, \pi\, t\,}\,} \, \exp \left\{ - { \,(\xi_1 + \xi_2 - z - \nu \, t )^2\, \over 2 \, t} \right\}
$$ 

\medskip
\noindent
is continuous and strictly positive on all of $\, \R^2\,$, and of class $\,
{\cal C}^{\,\infty}\, $ on $\,  \, \R^2 \setminus \{ (\xi_1 , \xi_2):\xi_1 =
\xi_2 \}  \,$.  

\subsection{ \textsc{The Degenerate Case with} $\,\, y=x_1- x_2\ge 0$}
\label{sec5.2}

Let us focus now on the degenerate case with $\, \sigma =0\,$, thus $\,
\rho=1\,$, and with $\, x_1 \ge x_2\,$ as in Figure 1. From formulae
(6.5.9)-(6.5.11), page 440 in 
\cite{MR1121940}%
, we have the joint probability distribution computations 
\begin{multline}
  \label{5.4}
  \mathbb{P}_\star \big( \,   Y  (t)   \in \d  a \,,   ~ 2\, L^Y (t) \in \d  b
  \, \big)=\\ 
  \frac{( |a| + b + y)}{\sqrt{\, 2\, \pi\, t^3\,}\,} \cdot \exp \cbr{ - \frac{ (|a|+b+y)^2} {2 \, t}} \d  a\, \d  b\,,\quad a \in \R\,, ~~ b >0
\end{multline}
as well as 
\begin{multline}
  \label{5.5}
  \mathbb{P}_\star \big( \,   Y  (t)   \in \d  a \,,   ~ 2\, L^Y (t) =0 \,
  \big)= \\
  \frac1{\,\sqrt{\, 2\, \pi\, t \,}\,} \left(\exp\cbr{-\frac{( a-y)^2}{2 \,
        t}} -  
    \exp\cbr{ - \frac{( a+y)^2}{ 2 \, t}} \right)  \d  a\,  ,\quad a >0\,, 
\end{multline}
which are based on the theory of the so-called ``elastic Brownian
motion". Substituting into (\ref{5.2}) with $\, \Theta = \R\,$ we obtain
from these expressions  
\begin{multline}
  \label{5.6}
  \mathbb{P} \br{\,Y^+(t)\in\d a\,,~Y^-(t)=0\,,~2\,L^Y(t)\in\d b\,}= \\
  \shoveleft{\mathbb{P} \br{\,Y^-(t)\in\d a\,,~Y^+(t)=0\,,~2\,L^Y(t)\in\d b\,} =} \\
  = \exp\cbr{\lambda \big( y+ b-a \big) - \frac{ \lambda^2}2 \,
    t \,} \cdot \, \frac{ \,   ( a + b + y)\,}{\sqrt{\, 2\, \pi\, t^3\,}\,} 
  \cdot \exp\cbr{ - \frac{ (a+b+y)^2\,}{2 \, t}}\,\d a\,\d b\\
  =\, e^{\, - 2 \lambda a}  \cdot \,\frac{\,( a + b+y)\,}{\sqrt{\, 2\, \pi\,
      t^3\,}\,} \, 
  \exp \cbr{-\frac{ (a+b +y- \lambda \, t)^2\,}{2 \, t}} \d  a\, \d  b\,, \quad a>0\,, \, b>0\,.  
\end{multline}

\medskip
\noindent
In a similar fashion, we obtain 
\begin{multline}
  \label{5.7}
  \mathbb{P}\big(\,Y^\pm (t)\in\d a\,,~ Y^\mp(t)=0\,,~ 2\,L^Y(t)=0\,\big)\,= \\
  =\exp\cbr{ \lambda \big( y -a \big) - \frac{\lambda^2}2 \, t}\,  \cdot \, 
  \frac1{\sqrt{\, 2\, \pi\, t\,}}\,\left(  
    \exp\cbr{ - \frac{ ( a-y)^2\,}{2 \, t}} -  \exp\cbr{- \frac{( a+y)^2\,}{2
        \, t}} \right)\d  a\\
  =  \frac1{\sqrt{\, 2\, \pi\, t\,}} \left(\exp\cbr{ - \frac{ ( a-y +
        \lambda\, t)^2\,}{2 \, t}} - 
    e^{\,-\, 2  \lambda  a\,} \exp \cbr{- \frac{( a+y- \lambda\, t)^2\,}{2 \,t}} \right) \d  a
\end{multline}
for $\, a >0\,$. This expression vanishes, as it should, for $\, y=0\,$: in
this case the process $\, Y(\cdot)\,$ starts accumulating local time at the
origin straightaway, that is, $\, \mathbb{P}  ( L^Y (t) >0 ) =1\,$  holds
for every $\, t \in (0, \infty)\,$.   

\medskip
\noindent
$\bullet~$ We set out to compute the joint distribution of the random vector
$\, (X_1 (t), X_2(t))\,$ for given, fixed $\, t >0\,$. From Remark
\ref{Remark 3} and either (\ref{2.27}) or (\ref{4.7.a}), it is clear that
this distribution is supported in the planar region $\, \mathfrak{ B}_1
\cup\, \mathfrak{ B}_2\,$, the union of the two blunt (135$^\mathrm{o}\,$-)
wedges 
$$
\mathfrak{ B}_1\,:=\, \big\{ (\xi_1, \xi_2) : \xi_1 > \xi_2\,, ~ \xi_2 \le x_2 + g\, t \big\}\,, \quad    \mathfrak{ B}_2\,:=\, \big\{ (\xi_1, \xi_2) : \xi_1 < \xi_2\,, ~ \xi_1 < x_2 + g\, t \big\}\,.
$$

\medskip
\noindent
$\bullet~$ 
Let us consider the wedge $\, \mathfrak{ B}_1\,\,$ first. With given real
numbers $\, \xi_1\,$, $\,  \xi_2\,$ that satisfy $\,  \xi_1 > \xi_2\,, ~
\xi_2 \le x_2 + g\, t \,$,   setting 
$$
\mathfrak{C} \,:=\, \Big\{ (a,b) \in (0, \infty)^2\,:\,  x_2+ g\, t+  a -  \big(    b/ 2\big)    \ge \xi_1\,, \, ~x_2+ g\, t- \big(    b/ 2\big)   \le \xi_2 \Big\} \,,
$$ 
and with the help of the expressions (\ref{2.27}) and (\ref{5.6}), we obtain 
\begin{multline*}
\P{X_1 (t) \ge \xi_1, \, X_2 (t) \le \xi_2 }=\\
\int\int_{\mathfrak{C}}  ~ e^{\, - 2 \lambda a}  \cdot \,{ \,   ( a + b+y)\,
  \over \sqrt{\, 2\, \pi\, t^3\,}\,} \, \exp \left\{ - { (a+b+y - \lambda \,
    t)^2\, \over 2 \, t} \right\} \d  a\, \d  b = 
 \\\int_{2 (x_2 + g\, t - \xi_2)}^\infty \left( \int_{\frac b2 +\xi_1 -x_2
    - g\,t}^{\infty} e^{\, - 2 \lambda a}  \cdot \frac{ \,   ( a + b+y)\,}{
    \sqrt{\, 2\, \pi\, t^3\,}}  
  \exp \cbr{ - { (a+b+y - \lambda \, t)^2\, \over 2 \, t}} \d  a \right) \d  b \,.  
\end{multline*}
We differentiate this expression, first with respect to $\, \xi_2\,$, then
with respect to $\, \xi_1\,$; recalling the notation of (\ref{2.1}), we
obtain 
\begin{multline}
  \label{5.8}
  \mathbb{P} \big( X_1 (t) \in \d  \xi_1, \, X_2 (t) \in \d \xi_2 \big) =\\
  2\, e^{\, - 2\, \lambda\, (\xi_1 - \xi_2 )\,} \cdot \frac{\, \, \xi_1 - 3\,
    \xi_2  +z +  2 \, g\,t   \, \,}{\sqrt{\, 2\, \pi\, t^3\,}} 
  \cdot\, \exp \left\{ - { 1 \over \, 2\, t\,} \,\Big( \xi_1 - 3\, \xi_2  +z  +
    \nu\, t \Big)^2 \right\} \, \d  \xi_1  \, \d  \xi_2\,;\\ 
  \xi_1 > \xi_2\,, ~ ~~\xi_2 < x_2 + g\, t \,.
\end{multline}
On the other hand, with $\, \xi_1 > \xi_2= x_2 + g\, t \,$ there is no
accumulation of local time at the origin over the interval $\, [0, t]\,$, so
the expressions of (\ref{2.27}) and (\ref{5.7}) give  
\begin{multline}
  \label{5.9}
    \mathbb{P} \big( X_1 (t) \in \d  \xi_1, \, X_2 (t)=  \xi_2 \big)\,=\, {\,
      1\, \over \,\sqrt{\, 2\, \pi\, t\,}\,} \cdot \bigg( \exp \left\{ - \,{ (
        \xi_1-x_1 + h\, t)^2\, \over 2 \, t} \right\}\\
    - \,e^{\, - \,2 \lambda ( \xi_1 - \xi_2)}\, \cdot \, \exp \left\{ - \,{ ( \xi_1 - 2 \xi_2+x_1  - h  \, t)^2\, \over 2 \, t} \right\} \bigg) \,\d  \xi_1\,\,; \quad \xi_1 > \xi_2 = x_2 + g\, t\,.
\end{multline}
\noindent
$\bullet~$ 
Next, we consider the wedge $\, \mathfrak{ B}_2\,\,$; equivalently, we work
on the event $\, \big\{ X_1 (t) < X_2 (t)\big\}\,$, on which the initial
order $\, x_1 \ge x_2\,$  stands reversed at time $\,t\,$ and     $\,   L^Y
(t) >0  \,$ holds a.e. In particular, the joint distribution of the random
vector $\, (X_1 (t), X_2(t))\,$ assigns zero mass to the region $\, \{
(\xi_1, \xi_2) : \xi_2 > \xi_1 = x_2 + g\, t \}\,$, as we have already
observed.  

With given real numbers $\, \xi_1\,$, $\, \xi_2\,$ that satisfy $\,  \xi_1 <
\xi_2\,, ~ \xi_1 < x_2 + g\, t \,$,   denoting 
$$
\mathfrak{D} \,:=\, \Big\{ (a,b) \in (0, \infty)^2\,:\,  x_2+ g\, t  -  \big(    b/ 2\big)    \le \xi_1\,, \, ~x_2+ g\, t+a- \big(    b/ 2\big)   \ge \xi_2 \Big\} \,,
$$ 
and with the help of the expressions (\ref{2.27}) and (\ref{5.6}), we obtain
then  
\begin{multline*}
\mathbb{P} \big( X_1 (t) \le \xi_1, \, X_2 (t) \ge \xi_2 \big)=\\
 \int\int_{\mathfrak{D}}   e^{\, - 2 \lambda a}  \cdot \,{ \,   ( a + b+y)\, \over \sqrt{\, 2\, \pi\, t^3\,}\,} \, \exp \left\{ - { (a+b+y - \lambda \, t)^2\, \over 2 \, t} \right\} \d  a\, \d  b= \\
\int_{2 (x_2 + g\, t - \xi_1)}^\infty \left( \int_{\frac b2 %
    +\xi_2 -x_2 - g\,t}^{\infty} ~e^{\, - 2 \lambda a}  \cdot \,{ \,   ( a + b+y)\, \over \sqrt{\, 2\, \pi\, t^3\,}\,} \, \exp \left\{ - { (a+b+y - \lambda \, t)^2\, \over 2 \, t} \right\} \d  a \right) \d  b \,.  
\end{multline*}
Differentiating this expression, first with respect to $\, \xi_1\,$ and then
with respect to $\, \xi_2\,$,  we obtain the following expression for the
probability density function:  
\begin{multline}
  \label{5.10}
  \mathbb{P} \big( X_1 (t) \in \d  \xi_1, \, X_2 (t) \in \d \xi_2 \big)\,=\,\\
  2\, e^{\, - 2\, \lambda\, (\xi_2 - \xi_1 )\,} \cdot \frac{\, \, \xi_2 - 3\,
    \xi_1  +z +  2 \, g\,t   \, \,}{\sqrt{\, 2\, \pi\, t^3\,}} %
  \cdot\, \exp \left\{ - { 1 \over \, 2\, t\,} \,\Big( \xi_2 - 3\, \xi_1  +z
    + \nu\, t \Big)^2 \right\} \, \d  \xi_1  \, \d  \xi_2\,;\\ %
  \xi_1< \xi_2\,, ~ ~~\xi_1 < x_2 + g\, t \,. 
\end{multline}

This is the same as the expression on the right-hand side of (\ref{5.8}),
except   $\, \xi_1\,$, $\, \xi_2\,$ have now traded places.

\begin{figure}[H]
\centerline{ $$\boxed{\hbox{~~The following Figure plots the joint
      Probability Density Function of $\, (X_1 (t), \, X_2 (t))\,$.~\,}
  } $$ } 
  \vspace{4ex}
   \begin{center}
        \includegraphics[width=39mm]{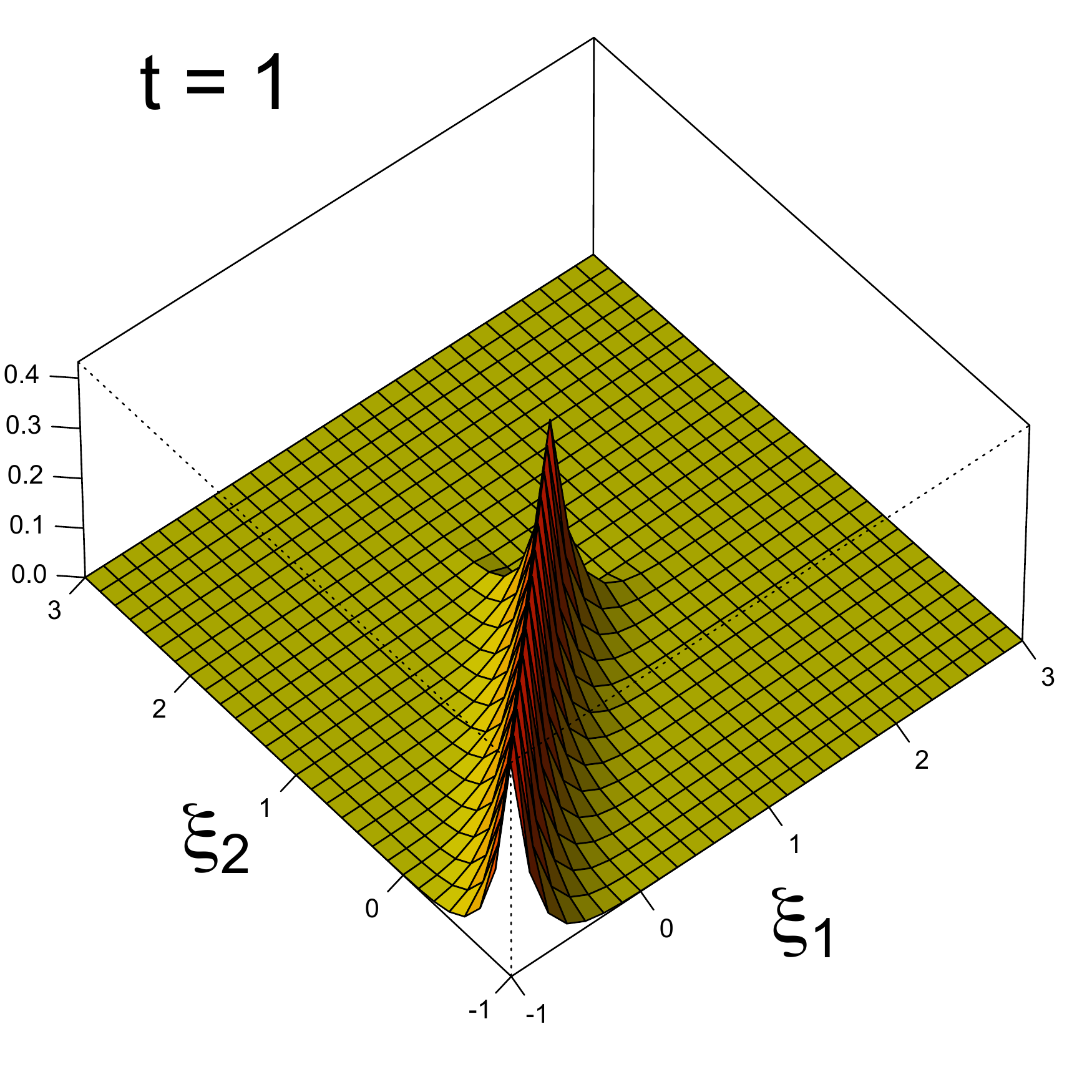}
        \includegraphics[width=39mm]{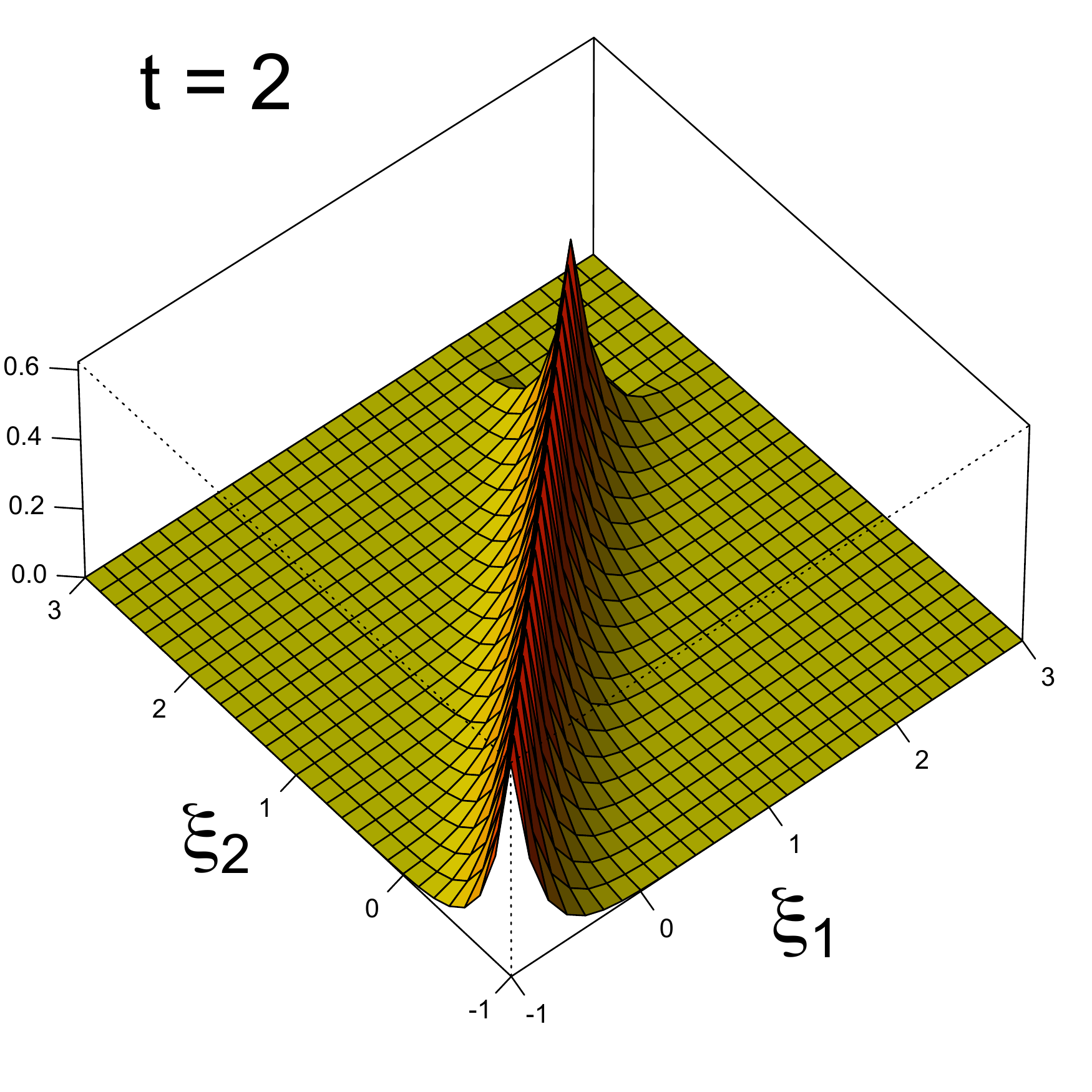}
    \caption{Joint density: $\, g = h = 1\, $, $\, \sigma = 0\, $, $\, \rho = 1\, $, $\, x_{1} = x_{2} = 0\, $, $\, t = 1\, $ (left), $\, t = 2\, $ (right).}
  \end{center}
\end{figure}

\begin{remark}
  \label{density}
  The joint distribution for the ranks $\,R_1 ( t) = \max ( X_1
  (\cdot), X_2 ( t) )\, , ~ \, R_2 ( t)\, = \,\min ( X_1 ( t), X_2 ( t)
  )\,$ is supported by the planar region $\, \{ (\rho_1, \rho_2) : \rho_1
  > \rho_2\, , ~ \rho_2 \le x_2 + g\, t\}\,$; from
  (\ref{5.8})-(\ref{5.10}), or directly from $\, R_2 (t)= r_2 + g\, t -
  L^Y (t)\,$ and $\, R_1 (t)= r_2+ g\, t + | Y (t)|- L^Y (t)\,$ and
  (\ref{5.4})-(\ref{5.5}),    it is computed as     
  \begin{multline}
    \label{5.11}
      \mathbb{P} \big( R_1 (t) \in \d  \rho_1,\, R_2 (t) \in \d \rho_2\big) = \\
      \frac{\,4 \,\big( 2\,    g \,t +z + \rho_1 -  3\,
        \rho_2 \big) \,  }{\,\sqrt{\, 2\, \pi\, t^3\,}\,}\, %
      \cdot \, \exp \cbr{ - 2 \, \lambda \, \big( \rho_1 - \rho_2\big) - { 1
          \over \, 2\, t\,} \,\Big( z  + \nu \, t + \rho_1 - 3 \,\rho_2
        \Big)^2  \, } \, \d  \rho_1\, \d  \rho_2 \,;%
      \\
      \rho_2 <   x_2+ g\,t\,, ~  ~\rho_1 > \rho_2  
  \end{multline}
  and
  \begin{multline}
    \label{5.12}
    \mathbb{P} \big( R_1 (t) \in \d  \rho_1, \, R_2 (t)=  \rho_2 \big) = \\
    \frac1{\sqrt{2 \pi t}} \cdot \bigg( \exp\cbr{ -\frac{(\rho_1-x_2 -y+ h\, t)^2}{2\,t}}%
    - e^{ -2 \lambda (\rho_1 - \rho_2) } \cdot  \exp \cbr{ - \frac{(\rho_1 - 2\rho_2 +x_1 - h  \, t)^2}{ 2 \, t}} \bigg)
    \d  \xi_1;\\ %
    \rho_1 > \rho_2 = x_2 + g\, t\,. 
  \end{multline}
 \end{remark}

\begin{remark}
  \label{gap at front}
  \noindent
  Once again, the probabilities in (\ref{5.9}) and (\ref{5.12}) vanish
  for $\, y=0\,$, i.e.,  when the two particles start at the same point.
  In this case  the distribution of $\, (X_1 (t), X_2 (t))\,$ is
  absolutely continuous with respect to Lebesgue measure on the plane,
  with probability density function $\, \mathfrak{ P}_t (\xi_1, \xi_2)\,$
  given by (\ref{5.8}) on $\, \xi_1 > \xi_2\,,  ~\xi_2 < x_2 + g\, t \,$
  (in the wedge $\, \mathfrak{ B}_1$), and    by (\ref{5.10}) on $\, \xi_1
  < \xi_2\,,  ~\xi_1 < x_2 + g\, t \,$ (in the wedge $\, \mathfrak{
    B}_2$).   

  Even in this case, though, there is a discontinuity along the front $\,
  \{ (\xi_1, \xi_2) \in \R^2 \,:\,  \xi_1 \wedge \xi_2 = x + g \,t \,\}\,$
  (cf. Figure 2), of size  
  $$
  \frac{\,2 \, \big| \xi_1 - \xi_2\big|\,}{\,\sqrt{\, 2\, \pi\, t^3\,}\,}\, \exp \left\{ - 2 \, \lambda \, \big| \xi_1 - \xi_2\big| - { 1 \over \, 2\, t\,} \,\Big( \big| \xi_1 - \xi_2 \big|  - \lambda \, t \Big)^2  \,  \right\} \,.
  $$
\end{remark}

\subsection{ \textsc{The Non-Degenerate Case with Unequal Variances} }
\label{sec6.3}

More generally, and by virtue of (\ref{4.2}), (\ref{5.6}) and (\ref{5.7}),
the quadrivariate joint distribution of $\, (Y^+ (t),\, Y^- (t),$ $ L^Y
(t),\, Q(t))\,$ is given as  
\begin{equation}
  \label{5.13}
  \mathbb{P} \big[ \, Y^\pm (t) \in \d  a \,, ~ Y^\mp (t) =0\,, ~ 2\, L^Y (t) \in \d  b \,, ~ Q (t) \in \d  \vartheta \, \big]\,=\,f_1 \big( a, b, \vartheta \big)\, \d  a \, \d  b \,  \d  \vartheta
\end{equation}
for  $\, a>0\,$, $\, b>0\,$,  $\, \vartheta \in \R\,$,  and 
\begin{equation}
  \label{5.14}
  \mathbb{P} \big[ \, Y^\pm (t) \in \d  a \,, ~ Y^\mp (t) =0\,, ~ 2\, L^Y (t) =0 \,, ~Q (t) \in \d  \vartheta \, \big]\,=\,f_2 \big( a,   \vartheta \big)\, \d  a \,    \d  \vartheta
\end{equation}
for  $\, a>0\,$,    $\, \vartheta \in \R\,$, with 
\begin{equation}
  \label{5.15}
  f_1 \big( a, b, \vartheta \big)\, 
  =\, e^{\, - 2 \lambda a}  \cdot \,{ \,   ( a + b+y)\, \over  \, 2\, \pi\, t^2\, } \, \exp \left\{ - { \, \vartheta^2 + (a+b +y- \lambda \, t)^2\, \over 2 \, t} \right\}\,,
\end{equation}
\begin{equation}
  \label{5.16}
  f_2 \big( a,   \vartheta \big)\,
  =\,  {\, e^{\, - \,\vartheta^2 / (2t)}\, \over  \, 2\, \pi\, t\, } \left(  \exp \left\{ - { ( a-y + \lambda\, t)^2\, \over 2 \, t} \right\} -   \exp \left\{ \, -\,2  \lambda  \,a- { ( a+y- \lambda\, t)^2\, \over 2 \,  t} \right\} \right)\,.
\end{equation}

\bigskip
\noindent
In conjunction with the skew representations of (\ref{2.23}) and
(\ref{2.24}), now written in the form 
$$
X_1 (t) \,=\, x_1 + \mu\, t + \Psi_1 (t)\,, \qquad X_2 (t) \,=\, x_2 + \mu\, t+ \Psi_2 (t)
$$
with 
\begin{align}
  \label{5.17}
  \Psi_1 (t) \,&:=\,  \rho^2 \,\big( \,Y^+ (t) - y^+ \big) - \sigma^2 \,\big(
  \,Y^- (t) - y^- \big)- \gamma\, L^Y (t) \,+  \, \rho\, \sigma\, Q (t)\,,\\
  \label{5.18}
  \Psi_2 (t) \,&:=\,  - \sigma^2 \,\big( \,Y^+ (t) - y^+ \big) + \rho^2 \,\big( \,Y^- (t) - y^- \big) - \gamma\, L^Y (t) \,+  \, \rho\, \sigma\, Q (t)\,,
\end{align}
it clearly suffices to compute the joint distribution of $\, ( \Psi_1 (t),\,
\Psi_2 (t))\,$. This is facilitated by the observation that the system of  
(\ref{5.17}), (\ref{5.18}) can be ``inverted", in the sense 
$$
Y^+ (t) - y^+ \,=\, { 1 \over \, \gamma\,}\, \Big( \rho^2\, \Psi_1 (t) + \sigma^2\, \Psi_2 (t) \Big) + L^Y (t)  - { \, \rho \, \sigma \,  \over \, \gamma\,}\,Q(t)\,,
$$
$$
Y^- (t) - y^- \,=\, { 1 \over \, \gamma\,}\, \Big( \sigma^2\,   \Psi_1 (t) + \rho^2\,  \Psi_2 (t) \Big) + L^Y (t)  - { \, \rho \, \sigma \,  \over \, \gamma\,}\,Q(t)\,.
$$

To proceed further, the cases $\, y \ge 0\,$, $\,y <0\,$ and $\, \gamma
>0\,$, $ \,\gamma<0\,$ have to be considered separately. We shall discuss
briefly only the case  
$\,
y=x_1- x_2\ge 0 \,,~~ \gamma >0 ~~~(\mathrm{i.e.,} ~~\rho >\sigma >0)\, 
$
In this case, the  joint probability density function 
$$
\mathbb{P} \big( \Psi_1 (t) \in \d  \psi_1\,, ~ \Psi_2 (t) \in \d  \psi_2 \big)\,=\, \mathfrak{ P}_t \big( \psi_1, \, \psi_2 \big)\, \d  \psi_1\, \d  \psi_2
$$ 
of $\, ( \Psi_1 (t),\, \Psi_2 (t))\,$ in (\ref{5.17}), (\ref{5.18}) is given by 
\begin{multline}
  \label{5.19}
  \gamma \, \,\mathfrak{ P}_t \big( \psi_1, \, \psi_2 \big)\,=\,\\
 = \int_{-\infty}^{\,\thetab (\psi_1, \psi_2)} \left\{  \int_0^\infty f_1 \big(
    b - \mathfrak{b}^*  (\psi_1, \psi_2,   \vartheta) \big)\, \d  b + f_2
    \big( - \mathfrak{ b}^*  (\psi_1, \psi_2,   \vartheta) \big) \right\} \,
  \d  \vartheta\\ 
  + \int^{ \infty}_{\,\thetab (\psi_1, \psi_2)} \left( \, \int^{
      \infty}_{\,\mathfrak{ b}^* (\psi_1, \psi_2, \vartheta)} \, f_1 \big( b-
    \mathfrak{ b}^*  (\psi_1, \psi_2,   \vartheta) \big)\, \d  b \right) \d\vartheta\,, 
\end{multline}
where we have set 
$$
\thetab (\psi_1, \psi_2) := {\, \psi_1 + \psi_2 + \gamma \, y \, \over 2\, \rho\, \sigma}\,, \qquad \mathfrak{b}^*  (\psi_1, \psi_2,   \vartheta) :=   { \, 2  \rho  \sigma \over \gamma} \big( \vartheta - \thetab (\psi_1, \psi_2) \big)\,.
$$

\noindent
After some calculations, (\ref{5.19}) reduces to 
\begin{multline}
  \label{5.20}
  \mathfrak{ P}_t \big( \psi_1, \, \psi_2 \big) = 
  \frac{1}{\sqrt{2 \pi t}} \Big [   ( 1+ \Phi \big(A_{+}(\psi_{1},
  \psi_{2})\big) 
  \exp\cbr{- \frac{1}{2t} \big( \psi_{1} + \psi_{2} + \lambda \gamma t \big)^{2}}\\
  - \,e^{\,2 \lambda \gamma^{2} y} \,\Phi\big(A_{-}(\psi_{1}, \psi_{2})\big) \exp \Big \{ - \frac{1}{2t} 
  \big( \psi_{1} + \psi_{2} + 2\gamma y + \lambda \gamma t \big)^{2}\, \Big\}
  \Big]\, ,
\end{multline}
where $\,\Phi(\cdot):=\big(1/\sqrt{2\pi\,}\,\big)\int^{\cdot}_{-\infty}
e^{\,-z^{2}/2} \,\d  z\,$ is the cumulative
standard normal distribution function and  
\[
\sqrt{t\,\,}\,A_{\pm}(\psi_{1}, \psi_{2}) \,:=\, \frac{\gamma}{2 \rho \sigma} \big(\psi_{1} + \psi_{2} + \gamma y\big) 
\pm 2 \, \rho \, \sigma \big( y \mp \lambda t \big) \, , \qquad 
(\psi_{1}, \psi_{2}) \in \R^{2}\, . 
\]

\section{ \textsc{Time Reversal} }
\label{sec6}

Let us consider now the time-reversed versions (``free'' and ``anchored'',
respectively)  
\begin{equation}
  \label{6.1}
  \widehat{Y} (t) \,:=\, Y (T-t)\,, \qquad \widetilde{W} (t) \,:=\, W (T-t) - W(T)\,,\qquad 0 \le t \le T\,,
\end{equation}
over a given time horizon $\, [0, T]\,$ of finite length, of the diffusion
process $\, Y(\cdot)\,$ in (\ref{3.1}) and of the Brownian motion $\,
W(\cdot)\,$ which drives that equation. Both of the  processes introduced in
(\ref{6.1}) are adapted to the {\it backwards filtration} $\,
\widehat{\mathbb{F}} = \{ \widehat{\F} (t) \}_{0 \le t \le T}\,$ defined as 
\begin{equation}
  \label{6.2}
  \widehat{\F} (t) \,:=\, \sigmab \big( Y(T) \big)
  \vee\mathfrak{F}^{\,\widetilde{W}} (t)\,,\qquad
  \mathfrak{F}^{\,\widetilde{W}} (t) \,:=\, \sigmab \,\big(\, \widetilde{W}
  (\theta)    ;\, 0 \le \theta \le t\big)\,. 
\end{equation}
Note that the process $\, \widetilde{W} (\cdot)\,$ is Brownian motion with
respect to the filtration $\, \mathbf{F}^{\,\widetilde{W}} =
\{\mathfrak{F}^{\,\widetilde{W}} (t) \}_{0 \le t \le T}\,$; this is because
$\, \widetilde{W}(\cdot)\,$  has continuous paths, and is easily checked to
be an $\, \mathbf{F}^{\,\widetilde{W}}-$martingale with the same quadratic
variation  as   Brownian motion. However, $\, \widetilde{W}(\cdot)\,$ is
only a semimartingale with respect to the larger backwards filtration $\,
\widehat{\mathbb{F}}\,$; its semimartingale decomposition is provided by the
fact that {\it the process 
  \begin{equation}
    \label{6.3}
    W^\# (t) \,:=\, \widetilde{W}(t) - \int_0^t \mathfrak{q} \big( T-s, \widehat{Y} (s) \big)\, \d  s  
    \,=\, W(T-t) - W(T) - \int_{T-t}^T \mathfrak{q} \big(  s,   Y  (s) \big)\, \d  s 
  \end{equation}
  for $\,0 \le t \le T\,$, is a standard $\, \widehat{\mathbb{F}}-$Brownian motion}. Here 
\begin{equation}
  \label{6.4}
  \mathfrak{q} \big( \tau , \xi\big) \,:=\, 
  { \partial  \over \, \partial \xi\,} 
  \log  \mathfrak{p}_\tau \big( y , \xi\big)\,, \qquad (\tau, \xi)
  \in (0, \infty) \times \R 
\end{equation}
is the logarithmic derivative of the transition probability density function
in (\ref{3.2}), (\ref{3.3}).  

\medskip
\begin{remark}
  This result is proved as in 
  \cite{MR1329103}%
  ,  who notes also the following  corollary: the (once-more-time-reversed)
  process     
  \begin{multline}
    \label{6.7.a}
    \etab (t) \coloneq 
    W^\# (T-t) - W^\#(T) \,
    =\, Y(t) + \int_0^t \Big( \mathfrak{q} \big(  s, Y (s) \big) + \lambda\,   \sign\big(    Y (s) \big) \Big)  \,\d  s\\
    = \,W(t) + \int_0^t \mathfrak{q} \big(  s, Y (s) \big)\, \d  s \,=\,
    \int_0^t \Big( \sign \big(Y(s) \big)\, \d  V^\flat(s) + \mathfrak{q} \big(
    s, Y (s) \big)  \, \d s\Big) \,, \quad 0 \le t \le T
  \end{multline}
   is also   Brownian motion {\it with respect to its own filtration,}  and
  is adapted to the filtration ${\bf F}^Y \equiv {\bf F}^W\,$ (though not to
  the filtration $\,{\bf F}^{V^\flat}\,$; recall %
 the filtration identity (\ref{4.3.a}) in this regard).  

  The process $\, \etab (\cdot)\,$ of (\ref{6.7.a})  is independent of the
  random variable $\, Y(T)$. In fact, for every given $\, t \in [0,T)\,$,
  the $\, \sigma-$algebra $\, \mathfrak{ F}^{\,\etab} (t)= \sigma ( \etab
  (s)\,, \, 0 \le s \le t ) \,$  is independent of the $\, \sigma-$algebra
  $\,   \sigma ( Y (\theta)\,, \, t \le \theta \le T ) \,$; in particular,
  of the random variable $Y(t)$. This shows that the inclusion $\,
  \mathfrak{ F}^{\,\etab} (t) \subset \mathfrak{ F}^{\,Y} (t)\,$, quite
  obvious from (\ref{6.7.a}), is strict; or, put another way, that the
  stochastic integral equation 
  \begin{equation}
    \label{6.7.b}
    Y(t) \,=\,y  - \int_0^t \Big(\lambda \,\mathrm{ sgn} \big(    Y (s)
    \big) + \mathfrak{q} \big(  s, Y (s) \big) \Big) \d  s \,+\,
    \etab (t)\,, \qquad  0 \le t \le T 
  \end{equation}
  cannot possibly have a strong solution.\hfill$\Box$
\end{remark}

\smallskip
A   result by  %
\cite{MR866342} %
 (see also  the expository papers by 
\cite{MR942014} and  %
\cite{MR1329103}, among many other works on this subject) is that, with
respect to the 
backwards filtration $\, \widehat{\mathbb{F}} = \{ \widehat{\F} (t) \}_{0
  \le t \le T}\,$,    the time-reversed process $\, \widehat{Y}(\cdot)\,$ in
(\ref{6.1}) is not just a semimartingale   but also a {\it diffusion
  process} driven by the $\, \widehat{\mathbb{F}}-$ Brownian motion in
(\ref{6.3}), of the form 
\begin{equation}
  \label{6.8}
  \widehat{Y} (t) \,=\, \widehat{Y} (0) + \int_0^t \,\widehat{\mathfrak{b}} \,\big( T-s, \widehat{Y} (s) \big) \, \d  s \,+\, W^\# (t) \,, \qquad 0 \le t \le T\,.
\end{equation}
Here, the new (backward) drift function $\, \widehat{\mathfrak{b}} (\cdot\,,
\cdot)\,$ is given in the notation of (\ref{6.4}) by the {\it generalized
  \textsc{Nelson} equation}  
\begin{equation}
  \label{6.9}
  \widehat{\mathfrak{b}} (\tau, \xi)  \,=\,  \lambda \, \sign \big( \xi \big) + \mathfrak{q} (\tau, \xi)\,, \qquad (\tau, \xi) \in (0, T] \times \R\,.
\end{equation}

\medskip
\begin{remark}
  In the special case $\, y=0\,$, the function of (\ref{6.4}) takes the form 
  \begin{multline}
    \label{6.6}
    \mathfrak{q} \big( \tau, \xi\big)\, =\, \left( \varphi^{(\lambda)} (\tau,
      -\xi) + \lambda\, e^{\, 2 \lambda \xi} \int_{-\xi }^\infty
      \varphi^{(\lambda)} (\tau, -u) \, \d  u \right)^{-1} \cdot\\
    \cdot  \left( \Big( 2 \,\lambda - {\, \xi \, \over \tau } \Big)
      \,\varphi^{(\lambda)} (\tau, -\xi) + 2\, \lambda^2\, e^{\,  2 \lambda
        \xi} \int_{-\xi }^\infty  \varphi^{(\lambda)} (\tau, -u) \, \d  u
    \right)\,, \qquad \xi \le 0\,, 
  \end{multline}
   and $\, \mathfrak{q}  ( \tau,  \xi    )     = -\, \mathfrak{q}  ( \tau, - \xi  )\,$ for $\, \xi >0\,$,    $\, \tau >0\,$ with 
  \begin{equation}
    \label{6.7}
    \varphi^{(\lambda)} (\tau, \xi) :=
    \frac1{\sqrt{\,2\, \pi\, \tau \,}}\, \exp\cbr{ - \frac{ \, ( \xi + \lambda\, \tau )^2\,}{ 2\,
        \tau}}\,. 
  \end{equation}  
   In particular,  we obtain from the equations   (\ref{6.6})-(\ref{6.7}) the
  explicit expression  
  \begin{multline}
    \label{6.10}
    \widehat{\mathfrak{b}}  \big( \tau, \xi\big)\, =\,\lambda \, \sign (\xi)
    -\left( \varphi^{(\lambda)} (\tau, |\xi|) + 
      \lambda\, e^{\,- 2 \lambda |\xi|} \int_{| \xi | }^\infty
      \varphi^{(\lambda)} (\tau, -u) \, \d  u \right)^{-1} \cdot \\
    \cdot  \left( \Big( 2 \,\lambda + {\, | \xi | \, \over \tau} \Big)
      \,\varphi^{(\lambda)} (\tau, | \xi |) + 2\, \lambda^2\, e^{\,- 2 \lambda
        | \xi |} \int_{| \xi |}^\infty  \varphi^{(\lambda)} (\tau, -u) \, \d
      u \right)\,, \qquad \xi \in \R\,.  
  \end{multline}
   The singularity at $\,\tau=0\,$ of the drift in (\ref{6.10})  
  is of the ``bridge" type: it ensures that, as the time-to-go $\,\tau
  \downarrow 0\,$ decreases to zero, the backward diffusion $\, \widehat{Y}
  (\cdot)\,$  zooms into the prescribed terminal condition $\, \widehat{Y}
  (T) =0\,$ (the initial condition $\, Y(0)=0\,$ of the forward process), as
  it must. \hfill$\Box$
\end{remark}

Let us also note that the semimartigale local time $\, L^{\widehat{Y}}
(\cdot)\,$ of the backward diffusion $\, \widehat{Y} (\cdot)\,$, and the
semimartigale local time $\, L^{Y} (\cdot)\,$ of the forward  diffusion $\,
Y (\cdot)\,$, are linked via 
\begin{equation}
  \label{6.12}
  L^{\widehat{Y}} ( t)\,=\, L^{Y} (T) - L^{Y} (T-t)\,, \qquad 0 \le t \le T\,.
\end{equation}

\subsection{ \textsc{A Time Reversal for Names} }
\label{sec6.1}

We consider now  the ``free" time-reversals
\begin{equation}
  \label{6.13}
  \widehat{X}_1 (t) \,:=\, X_1(T-t)\,, \quad \widehat{X}_2 (t) \,:=\, X_2(T-t)\,,\qquad 0 \le t \le T
\end{equation}
of the components of the vector process $\, (X_1 (\cdot), X_2 (\cdot))\,$ constructed in section \ref{sec3}, as well as their ``anchored" versions 
\begin{equation}
  \label{6.14}
  \widetilde{X}_j (t) \,:=\,  \widehat{X}_j (t) - \widehat{X}_j (0) \,=\, X_j(T-t)- X_j (T)\,;  \qquad 0 \le t \le T\,, ~~ j=1,\, 2\,.
\end{equation}
By analogy with (\ref{6.2}), we also look at the new backwards filtration $\, \widetilde{\mathbf{F}} = \{ \widetilde{\mathfrak{F}} (t) \}_{0 \le t \le T}\,$   given by
$$
\widetilde{\mathfrak{F}} (t) \,:=\, \sigmab \big( Y (T) \big) \vee \mathfrak{F}^{\,(\widetilde{  Q}, \widetilde{  W})} (t)\,,\quad \text{with} \quad \mathfrak{F}^{\,(\widetilde{  Q}, \widetilde{  W})} (t) = \sigmab \big( \,\widetilde{Q} ( \theta), \, \widetilde{W}(\theta)\,  ;\, 0 \le \theta \le t\big) 
$$
generated by the random variable $\, Y(T)\,$ and by the time-reversed versions 
\begin{equation}
  \label{6.15}
  \widetilde{Q} (t) = Q (T-t) - Q(T)\,, \qquad  \widetilde{W}  (t)=W (T-t) - W(T)\,, \qquad 0 \le t \le T
\end{equation}
of the independent Brownian motions $\, W(\cdot)\,$, $\, Q(\cdot)\,$ of
(\ref{2.20}), (\ref{2.22}) in the manner of (\ref{6.1}). In particular, $\,
\widetilde{Q}(\cdot)\,$ is independent of both $\, \widetilde{W}(\cdot)\,$
and $\, Y(T)\,$, thus also of the Brownian motion $\, W^{\#} (\cdot)\,$ in
(\ref{6.3}).  

\smallskip
Then the skew representations (\ref{2.23}) and (\ref{2.24}), along with the
notation of (\ref{6.1}) and the local time identity  (\ref{6.12}), imply
that  the ``anchored  time-reversals" of    (\ref{6.14})    are  $\,
\widetilde{\mathbf{F}}-$adapted and are given by  
$$
\widetilde{X}_1 (t) \,=\, - \mu\, t + \rho^2 \,\big( \, \widehat{Y}^+ (t) -
\widehat{Y}^+ (0) \big) - \sigma^2 \,\big( \, \widehat{Y}^- (t) -
\widehat{Y}^- (0) \big)+ \gamma\, L^{\widehat{Y}} (t) \,+  \, \rho\, \sigma
\, \widetilde{Q} (t)  
$$
\begin{equation}
  \label{6.15.a}
  ~~~=\, - \mu\, t + \int_0^t \big( \, \rho^2 \, \mathbf{ 1}_{ \{
    \widehat{Y} (s) > 0 \} }    + \sigma^2 \,  \mathbf{ 1}_{ \{ \widehat{Y}
    (s) \le 0 \} }  \big) \,\d  \widehat{Y}(s)+ 2\,\gamma\,
  L^{\widehat{Y}} (t) \,+  \, \rho\, \sigma \, \widetilde{Q} (t)\,, 
\end{equation}

\noindent
and 
$$
\widetilde{X}_2 (t) \,=\, - \mu\, t - \sigma^2 \,\big( \,\widehat{Y}^+ (t) -
\widehat{Y}^+ (0) \big) + \rho^2 \,\big( \,\widehat{Y}^- (t) - \widehat{Y}^-
(0)  \big) + \gamma\, L^{\widehat{Y}} (t) \,+  \, \rho\, \sigma \,
\widetilde{Q} (t)  
$$
\begin{equation}
  \label{6.15.b}
  ~~~~=\, - \mu\, t - \int_0^t \big( \, \rho^2 \, \mathbf{ 1}_{ \{
    \widehat{Y} (s) \le 0 \} }    + \sigma^2 \,  \mathbf{ 1}_{ \{
    \widehat{Y} (s) > 0 \} }  \big) \,\d  \widehat{Y}(s)+ 2\,\gamma
  \, L^{\widehat{Y}} (t) \,+  \, \rho\, \sigma \, \widetilde{Q} (t)\,, 
\end{equation}
respectively, thanks once again to the \textsc{Tanaka} formulas. 

\smallskip
We recall now the ``backwards dynamics" of (\ref{6.8})  as well as the
notation of (\ref{6.4}), (\ref{6.9}), and write these equations in the {\it
  time-reversed skew representation}  form 
\begin{multline}
  \label{7.33}
\widetilde{X}_1 (t)=\widehat{X}_1 (t) - \widehat{X}_1 (0)=\\
=\int_0^t  \left(
  h\,\mathbf{ 1}_{ \{ \widehat{X}_1 (s) >  \widehat{X}_2 (s) \} } -g\,
  \mathbf{ 1}_{ \{ \widehat{X}_1 (s) \le  \widehat{X}_2 (s) \} } \right)
\d  s \,+  \, ( \rho^2 - \sigma^2)\,    L^{|\widehat{X}_1 -
  \widehat{X}_2|} (t)   
  +  \, \rho\, \sigma  \,   \widetilde{Q} (t)\,+\\
  \int_0^t \big( \,
  \rho^2 \, \mathbf{ 1}_{ \{ \widehat{X}_1 (s) >  \widehat{X}_2 (s) \} }
  + \sigma^2 \, \mathbf{ 1}_{ \{ \widehat{X}_1 (s) \le  \widehat{X}_2 (s) \}
  }   \big) \,\big[ \, \mathfrak{ q} \big( T-s, \widehat{X}_1 (s) -
  \widehat{X}_2 (s) \big) \,   \d  s +\,  \d  W^\#  (s)\,
  \big] 
\end{multline}
and 
\begin{multline}
  \label{7.34}
  \widetilde{X}_2 (t)=\widehat{X}_2 (t) - \widehat{X}_2 (0)=\\
=  \int_0^t
  \left(h\, \mathbf{ 1}_{ \{ \widehat{X}_1 (t) \le  \widehat{X}_2 (t) \} }
    -g\, \mathbf{ 1}_{ \{ \widehat{X}_1 (t) >  \widehat{X}_2 (t) \} } \right)
  \d  t\, +\,   ( \rho^2 - \sigma^2)\,    L^{|\widehat{X}_1 -
    \widehat{X}_2|} (t) \,
  +  \, \rho\, \sigma  \,   \widetilde{Q} (t)\,+\\
  \int_0^t  \big( \,
  \rho^2 \, \mathbf{ 1}_{ \{ \widehat{X}_1 (s) \le  \widehat{X}_2 (s) \} }
  + \sigma^2 \, \mathbf{ 1}_{ \{ \widehat{X}_1 (s) >  \widehat{X}_2 (s) \} }
  \big) \,\big[ \, \mathfrak{ q} \big( T-s, \widehat{X}_1 (s) -
  \widehat{X}_2 (t) \big) \,   \d  s +\,  \d  W^\#  (s)\,
  \big] .
\end{multline}

\begin{remark}
  \label{Werner}
  A somewhat interesting dichotomy emerges. In the case of equal variances
  ($\,\rho^2 = \sigma^2 = 1/2 \,$) these anchored time-reversals are $\,
  \widetilde{\mathbf{F}}-$adapted \textsc{It\^o}  processes: the  bounded
  variation terms in their semimartingale decompositions are absolutely
  continuous with respect to \textsc{Lebesgue} measure.  In the case of
  unequal variances $\, \rho^2 \neq \sigma^2\,$, terms which are singular
  with respect to \textsc{Lebesgue} measure  appear, and are governed by
  local time.   

  To the best of our knowledge, this is the first time such a  feature is
  observed in the context of time-reversal of a ``purely forward" stochastic
  differential equation without reflection; its occurrence and significance
  need to be understood further. For a similar but different phenomenon, in
  the context of time reversal of Brownian motion reflected on an
  independent time-reversed Brownian motion, see  \cite{MR1917545} %
  (as well as  \cite{MR1785393}, \cite{MR1902187}).  
\end{remark}

\subsection{ \textsc{A Time Reversal for Ranks} }
\label{sec6.2}

By analogy with (\ref{6.13}), we introduce the time-reversed versions
\begin{equation}
  \label{6.20}
  \widehat{R}_1 (t) \,:=\,R_1(T-t) 
  \,=\, \max ( \widehat{X}_1 ( t), \widehat{X}_2 ( t) )\,, \quad
  \widehat{R}_2 (t) \,:=\,R_2(T-t)
  \,=\, \min ( \widehat{X}_1 ( t), \widehat{X}_2 ( t) ) 
\end{equation}
of the ranked processes in (\ref{4.1}) for $\, 0 \le t \le T\,$. We have $\,
\widehat{R}_1 (t) + \widehat{R}_2 (t)= \widehat{X}_1 ( t)+ \widehat{X}_2 (
t) \,$ and $\,  \widehat{R}_1 (t) - \widehat{R}_2 (t)= |\widehat{X}_1 ( t)-
\widehat{X}_2 ( t)| = | \widehat{Y} (t)|  \,$, so the time-reversed skew
representations of (\ref{6.15.a}), (\ref{6.15.b}) cast the ``anchored"
versions of the processes of (\ref{6.20}) in the form 
\begin{equation}
  \label{6.21}
  \widetilde{R}_1 (t) \,:=\, \widehat{R}_1 (t) - \widehat{R}_1 (0)\,=\, - \mu \, t + \rho^2 \big( \big| \widehat{Y} (t) \big| - \big| \widehat{Y} (0) \big| \big) + \gamma\, L^{\widehat{Y}} (t) + \rho \, \sigma\, \widetilde{Q} (t)\,,
\end{equation}
\begin{equation}
  \label{6.22}
  \widetilde{R}_2 (t) \,:=\, \widehat{R}_2 (t) - \widehat{R}_2 (0)\,=\, - \mu \, t - \sigma^2 \big( \big| \widehat{Y} (t) \big| - \big| \widehat{Y} (0) \big| \big) + \gamma\, L^{\widehat{Y}} (t) + \rho \, \sigma\, \widetilde{Q} (t)\,.
\end{equation}

\medskip
\noindent
In conjunction with the reverse-time dynamics of (\ref{6.8}), (\ref{6.9}),
the \textsc{Tanaka-Meyer} formulas give now 
\begin{equation}
  \label{6.23}
  \big| \widehat{Y} (t) \big| - \big| \widehat{Y} (0) \big|\,=\, \lambda\,t
  + \int_0^t \sign \big( \,\widehat{Y} (s) \big)\, \mathfrak{ q} \big( T-s,
  \widehat{Y} (s) \big)\, \d  s + V^{\#} (t)+ 2\, L^{\widehat{Y}}
  (t)\,, 
\end{equation}
where
\begin{equation}
  \label{6.23b}
  V^{\#} (t) \,:=\, \int_0^t \sign \big( \,\widehat{Y} (s) \big)\,
  \d  W^{\#} (s)\,, 
  \qquad 0 \le t \le T 
\end{equation}
is a Brownian motion of the backwards filtration $\, \widetilde{\mathbf{F}}
\,$, and is independent of the Brownian motion $\, \widetilde{Q} (\cdot)\,$.

Substituting the expression of (\ref{6.23}) back into (\ref{6.21}) and
(\ref{6.22}), and recalling (\ref{2.1}), (\ref{2.25}), we obtain the
dynamics 
\begin{equation}
  \label{6.25}
  \widehat{R}_1 (t) - \widehat{R}_1 (0)\,=\, h \, t + \rho^2 \int_0^t \sign
  \big( \,\widehat{Y} (s) \big)\, \mathfrak{ q} \big( T-s, \widehat{Y} (s)
  \big)\, \d  s + \rho\, V^{\#}_1 (t)+ \big( 4 \,\rho^2 - 1 \big) \,
  L^{\widehat{Y}} (t)\,,  
\end{equation}
\begin{equation}
  \label{6.26}
  \widehat{R}_2 (t) - \widehat{R}_2 (0)\,=\, - g \, t - \sigma^2 \int_0^t
  \sign \big( \,\widehat{Y} (s) \big)\, \mathfrak{ q} \big( T-s, \widehat{Y}
  (s) \big)\, \d  s + \sigma\, V^{\#}_2 (t)- \big( 4 \,\sigma^2 - 1
  \big) \, L^{\widehat{Y}} (t)\,. 
\end{equation}
Here the processes
\begin{equation}
  \label{6.27}
  V^{\#}_1 (\cdot)\,:=\, \rho\, V^{\#}  (\cdot) + \sigma \, \widetilde{Q}
  (\cdot)\,,\qquad V^{\#}_2 (\cdot)\,:=\, \rho\, \widetilde{Q} (\cdot) -
  \sigma \, V^{\#}  (\cdot)  
\end{equation}
are independent, standard Brownian motions of the  backwards filtration $\,
\widetilde{\mathbf{F}} \,$.  

Comparing the equations of (\ref{6.25}), (\ref{6.26}) with those of
(\ref{4.3}), (\ref{4.4}), we see that $\, V^{\#}_1 (\cdot)\,$, $\, V^{\#}_2
(\cdot)\,$ play in the context of time reversal the same r\^oles that the
Brownian motions $\, V_1 (\cdot)\,$, $\, V_2 (\cdot)\,$ play on the forward
context: to wit, that of independent, standard Brownian motions associated
with individual {\it ranks.} 

\subsection{ \textsc{Steady State}  }
\label{sec7.3}

Finally, let us note that the diffusion process $\, Y(\cdot)\,$ of
(\ref{3.1}) has invariant distribution with double exponential probability
density function  
$$\, 
\mathfrak{ p} (\xi) \,=\,   \lambda \,  \,e^{\, - \,2 \, \lambda \,
  |\xi|}\,\,,~ ~~~~~\,\xi \in \R\,. 
$$ 
For this function, the analogue $\, \mathfrak{q} (\xi) = \big( \partial
/ \partial \xi \big) \, \log\, \mathfrak{ p} (\xi)\,$ of the logarithmic
derivative in (\ref{6.4}) becomes $\, \mathfrak{q} (\xi) = -2\, \lambda\,
\sign (\xi)\,$ in the notation of (\ref{5.1}), and from the generalized
\textsc{Nelson} equation (\ref{6.9}) the drift in the backward equation
(\ref{6.8}) becomes $$\, \,\widehat{\mathfrak{b}} (\xi) \,=\,- \lambda\,
\sign  (\xi)\,.$$ This reflects the fact that the diffusion process $\,
Y(\cdot)\,$ of (\ref{3.1}) is strictly time-reversible  when started at its
invariant distribution, that is, the processes $\, Y(\cdot)\,$ and $\,
\widehat{Y} (\cdot)\,$ are then identically distributed; and then the
Brownian motion in (\ref{6.7.a}) takes the form  
$\,
\etab (t)= \int_0^t \,\sign \big( Y(s)\big)\, \big( {\rm d} V^\flat(s) - 2\,
\lambda\, {\rm d} s \big)\,$, $\,\,0 \le t \le T\,$. 

\smallskip
In such a  setting, the equations of  (\ref{6.7.b}) and (\ref{7.33}),
(\ref{7.34}) continue to hold, now with $\, \mathfrak{q} (\cdot) = -2\,
\lambda\, \sign (\cdot)$; whereas the equations of (\ref{6.25}),
(\ref{6.26}) assume the rather concrete form 
$$
\widehat{R}_1 (t) - \widehat{R}_1 (0)\,=\, \Big( h - 2\, (g+h)\, \rho^2
\,\Big) \, t   \,+\, \rho\, V^{\#}_1 (t)\,+\, { \,   4 \,\rho^2 - 1   \,
  \over 2}\,\, L^{\widehat{R}_1- \widehat{R}_2} (t)\,,  
$$
$$
\widehat{R}_2 (t) - \widehat{R}_2 (0)\,=\, \Big(   2 \,(g+h)\, \sigma^2 - g\Big) \, t   \, + \,\sigma\, V^{\#}_2 (t)\,-\, { \,      4 \,\sigma^2 -1   \, \over 2}\, \, L^{\widehat{R}_1- \widehat{R}_2} (t)\,,
$$
where $\, L^{\widehat{R}_1- \widehat{R}_2} (\cdot) 
= 2 \, L^{\widehat{Y}}(\cdot) \,$ is the local time accumulated 
at the origin by $\, \widehat{R}_1(\cdot) - \widehat{R}_2(\cdot) = 
\lvert \widehat{Y}(\cdot) \rvert \, $. 

\medskip

\section{  \textsc{A Generalization of the perturbed Tanaka  equation} }
\label{App}

We have the following generalization of Theorem \ref{Theorem 0}, which  partially   answers a question posed  by Professor Marc \textsc{Yor}.   In what follows, we shall agree to denote by $\, \langle \mathcal{Z} \rangle  (\cdot)\,$ the quadratic variation $\, \langle \mathcal{V} \rangle  (\cdot)\,$ of the continuous local martingale $\,   \mathcal{V}   (\cdot)\,$ in the decomposition of a continuous semimartingale $\,   \mathcal{Z}    (\cdot) = \mathcal{Z}    (0) + \mathcal{V}    (\cdot) + \mathcal{B}    (\cdot)\,$, where $\, \mathcal{B}    (\cdot)\,$ is a continuous, adapted process of finite variation on compact intervals. 

\begin{theorem}  
  \label{Theorem 4} 
  Let $\, f : \R \ra \R\,$ be a     function of finite variation, that
  is, $\, f = f_+ - f_-\,$ where $\, f_\pm : \R \ra \R\,$ are    increasing. On some filtered probability space $\, (\Omega, \F,
  \mathbb{P})$, $\mathbb{F}= \{  \mathfrak{ F}(t)\}_{t \ge 0}\,$, consider
  two  continuous local martingales   $\, M(\cdot)\,$, $  N(\cdot)\,$ which
  satisfy the conditions of Theorem \ref{Theorem 0}, and a
  continuous, adapted process $\, A(\cdot)\,$ with $A(0)=0$ and finite total 
  variation $\, \breve{A} (t)\,$ on compact intervals of the form $\, [0,t]\,$.
  
  Then pathwise
  uniqueness holds for the stochastic differential equation  
  \begin{equation} 
    \label{eq: prokaj 2} 
    Y(t)\, =\,y+  \int_0^t f \big(Y(s)\big)\, {\mathrm d} M(s) + A(t) +    N(t) \, , \quad 0 \le t < \infty\,.
  \end{equation}
\end{theorem}

\begin{proof} We shall prove the statement for $\,f=f_{+}-f_{-}\,$ with bounded, 
  increasing $f_{\pm}$. The  general case will then follow  by considering 
  $\,f^{(n)}= f_{+}^{(n)}-f_{-}^{(n)}$, where $\,f_{\pm}^{(n)}=\min \big(n, \max(-n, f_{\pm})\big) \,$, $n \in \mathbb{N}\,$ are       truncated versions of $f$. If     pathwise uniqueness holds in the bounded    case, then the solution is pathwise unique up to the stopping time    $\,\sup  \vartheta_n\,$, where     $\, \vartheta_n=\inf\set{t>0}{\abs{f_{\pm}(Y(t))}>n}$.     That is, if we have two    solutions, $X(\cdot)$, $Y(\cdot)$ to \eqref{eq: prokaj 2} on the same probability space, issued from the same initial value and driven by the processes $M(\cdot)$, $N(\cdot)$ and $A(\cdot)$, then     $X(t)=Y(t)$ holds for    $\,t<\zeta(X)\wedge\zeta(Y)$, where $\,\zeta(X),\,\zeta(Y)\,$ denote the lifetimes    of the processes $X(\cdot)$ and $Y(\cdot)$, respectively. 

  \smallskip
  
  Suppose then that the processes $\, X(\cdot)\,$, $\, Y(\cdot)\,$
  both solve the equation (\ref{eq: prokaj 2}), so their difference 
  \begin{equation}
    \label{D}
    D (t) \,:=\, X(t) - Y(t) \,=\, \int_0^t \big( f (X(s) ) - f (Y(s) )
    \big)\, \d  M(s)\,, \qquad 0 \le t < \infty 
  \end{equation}
  is a continuous local martingale with quadratic variation
  \begin{equation}
    \label{<D>}
    \q{D}(t) \,=\, \int_0^t \big( f (X(s) ) - f (Y(s) )  \big)^2\,
    \d  \q{ M } (s) \,\le\, \big\| f \big\|_{TV}
    \int_0^t \big| f (X(s) ) - f (Y(s) )  \big|\, \d  \q{M}(s)\,. 
  \end{equation}
  Here $\, \big\| f \big\|_{TV}\,$ is the total variation of $\,f$
  over the real line, and we have used the elementary comparison  $\,  ( f (x) -
  f(y)  )^2 \le \| f \|_{TV}\, |f(x)-f(y)|\,$, $~\forall~ (x, y) \in   \R^2\,$.  
  Following 
  \cite{0527.60062}%
  ,  we shall establish the  estimate
  \begin{equation}
    \label{eq:D L}
    \E \int_0^T \frac{\,\d \q{D}(t) \,}{\,D(t)\,} \, 
    \I{  D(t) > 0  }\, \le\, c \cdot \big\|f \big\|_{TV}^2
    \cdot \sup_{a \in \R \atop  u \in [0,1] } \E{2\,L^{(u)} (T,a)}    
  \end{equation}
  where, for each $\, u \in [0,1]\,$, the quantity  $\, L^{(u)} (T,a) \equiv \Lambda^{Z^{(u)}} (T,a)\,$ is the  local time accumulated  at the site $\, a \in \R\,$  during the   time interval $\, [0,T]\,$  by  the  continuous semimartingale  
  \begin{multline}
    \label{Zu}
    Z^{(u)}(\cdot) \,:=\, (1-u)\, X(\cdot) +\, u\, Y(\cdot)\,=\,\\ 
    y + \int_0^\cdot \Big((1-u)\, f \big( X(t) \big) + \,u\, f \big( Y(t) \big)\Big)\, \d  M(t) + A (\cdot) + N (\cdot)\,.~~~~~~~
  \end{multline}
  
  To this end, we introduce a sequence  $\, \{ f_k \}_{k \in \N} \subset
  \mathcal{C}^1(\R)\,$  of continuous and continuously differentiable
  functions that converge to $\,f\,$ pointwise, and are bounded in total
  variation norm, that is $\sup_{k\in\N}\norm{f_k}_{TV}<\infty$. Since $\limsup_{k}
  \norm{f_k}_{TV}\leq \norm{f}_{TV}$ obviously holds, it is only possible if
  $f$ is of bounded variation and in this case an approximating sequence is
  easily  obtained, e.g., by mollifiers.
  We note the identity 
  $\,\,
  f_k\big( X(t) \big)- f_k \big( Y(t) \big)\,=\, 
  \big( X(t) - Y(t) \big) 
  \int_0^1 f_k^{\prime}\big( Z^{(u)}(t) \big)\,\d u\,$, as well as the comparison 
  \begin{multline*}
      \E\int_0^T \frac{\,\big| f_k \big( X(t) \big) - f_k \big( Y(t) \big)
        \big| \,}{X(t) - Y(t)} \, \I{ X(t) - Y(t) > \delta}\,\d\q{M}(t)\\
      \le \int_0^1  \left( \E \int_0^T  \big| f_k^{\prime} \big( Z^{(u)}(t) \big) \big| \, \I{ X(t) - Y(t) > \delta}\, \d\q{M}(t) \right)\d u
  \end{multline*}  
  for $\delta> 0$.   We also note
  that the orthogonality of $\, M(\cdot)\,$ and $\, N(\cdot)\,$  implies 
  $\,\q{Z^{(u)}}(\cdot)\ge \q{N}(\cdot) \,$ in conjunction with
  (\ref{Zu}). Let us recall now the domination condition (\ref{dom}), which  gives   
  $$
  \q{M}(\cdot)\,\le\,c\,\q{N}(\cdot)\,
  \le\,c\,\q{Z^{(u)}}(\cdot)\,\le\,C\,\q{N}(\cdot) 
  $$
  for a suitable real constant $\, C >c\,$. We deduce from this
  \begin{align*}
    \int_0^1  &\left( \E \int_0^T  \big| f_k^{\prime} \big( Z^{(u)}(t) \big)
      \big| \, \I{ D(t) > \delta }\, \d  \q{M}(t)\right)\d u
    \,\\
    &\le c \, \int_0^1  \left( \E \int_0^T  \big| f_k^{\prime} \big(
      Z^{(u)}(t) \big) \big|\,\d\q{Z^{(u)}}(t)\right)\d u \\&
    =\, c \, \int_0^1  \left( \E \int_\R  \,\big| f_k^{\prime} ( a)   \big|
      \,2\,L^{(u)} (T,a)\,\d a\right)\d u \,\\ 
    &\le\, c \cdot \sup_{a \in \R \atop  u \in [0,1]} \E{2\,L^{(u)} (T,a)} 
    \cdot \int_\R  \,\big| f_k^{\prime} ( a)\big|\,\d  a  
    \le\, c \cdot \norm{f_k}_{TV}
    \cdot \sup_{a \in \R \atop u\in [0,1]} \E{ 2\,L^{(u)} (T,a)}\,.
  \end{align*}
   Letting $k\uparrow\infty$ and then $\delta\downarrow0\,$, 
  these estimates give  
  $$
  \E \int_0^T \frac{\,\big| f  \big( X(t) \big) - f  \big( Y(t) \big) \big|
    \,}{X(t) - Y(t)} \, \mathbf{ 1}_{\{ X(t) - Y(t) > 0 \}}\, \d\q{M}(t)\,
  \le\, c\cdot\norm{f}_{TV}\cdot\sup_{a \in \R \atop  u\in [0,1]}\E{2\,L^{(u)} (T,a)}  
  $$   
  and, in conjunction with (\ref{<D>}), the claim  
  \eqref{eq:D L} as well.

  \smallskip
  Suppose now we can show
  \begin{equation}
    \label{EL}
    \sup_{a \in \R \atop  u \in [0,1] }\E{2\,L^{(u)} (T,a)} \,<\, \infty\,.
  \end{equation}
  On the strength of \eqref{eq:D L}, this   then implies  
  $\,
  \mathbb{E} \int_0^T (D(t))^{-1} \, \mathbf{ 1}_{\{ D(t) > 0 \}}\,\d  \langle
  D \rangle (t)<  \infty \,$ for each $\,\, T \in (0, \infty)\,$;  and arguing
  as in  
  \cite{0527.60062}%
  ,  Lemma 1.0 (see also Exercise 3.7.12, pages 225-226 in
  \cite{MR1121940}%
  ),  we deduce that the local time $\,L
  (\cdot) \equiv  \Lambda^D (\cdot\,, 0)\,$, accumulated at the origin by the
  continuous 
  local martingale $\, D(\cdot)\,$ of (\ref{D}), is identically equal to
  zero. But then, by the \textsc{Tanaka} formula once again, we obtain that
  $\, | D(\cdot)|\,$ is a local martingale, thus also a (nonnegative)
  continuous    supermartingale   with $\, | D(0)|=0\,$, and consequently $\,    X(\cdot) -
  Y(\cdot)\equiv  D(\cdot)  \equiv 0\,$; that is, pathwise uniqueness holds.   
  
  \smallskip
  \noindent
  $ \bullet~$  The property (\ref{EL})  is checked by standard methods:   with
  the help of the \textsc{Tanaka} formula 
  $$
  \big| Z^{(u)}(T) - a \big| \,=\,  \big| Z^{(u)}(0) - a \big| + \int_0^T \sign \big( Z^{(u)}(t) - a \big)\, \d  Z^{(u)}(t) \,+\, 2\,L^{(u)} (T,a)
  $$
  and  the It\^o isometry (e.g., 
  \cite{MR1121940}%
  , p. 144)  we get
  \begin{align*}
    \E{2\,L^{(u)} (T,a)} &\leq 
    \E\abs{Z^{(u)}(T)-Z^{(u)}(0)}+\E^{1/2}\zjel{\q{Z^{(u)}}(T)} + 
    \E \,\big(    \breve{A} (T)  \big)  \\ &
    \leq 2\,\E^{1/2}\zjel{\q{Z^{(u)}}(T)} + 2\,\E \,\big( \breve{A} (T)
    \big)\, \\ &
    \leq 2\, \sqrt{C\,}\,\,\E^{1/2}\zjel{\q{N}(T)}+ 2 \,\E \,\big( \breve{A} (T)    \big)\,.
  \end{align*}
  This last quantity does not depend on $\, a \in \R\,$ or $\, u \in [0,1]\,$;
  if it is also finite, we are done.  
  
  If not, we deploy standard localization arguments: to wit, we consider the
  stopping times 
  $$\, \tau_m = \inf \, \big\{ \, t \ge 0\,:\, \max \big( \breve{A}(t),\q{N} (t) \big) \ge
  m \,\big \} \,,$$ 
  and deduce $\, D ( \,\cdot \,\wedge \tau_m) \equiv 0\,$ from the above
  analysis, for every $\, m \in \N\,$, $\,\mathbb{P}-$a.s. But then $\,
  \lim_{m \to \infty} \tau_m = \infty\,$ also holds
  $\,\P$--a.s., and this leads to $\, D ( \,\cdot  \,) \equiv 0\,$ once
  again. 
\end{proof} 

Our next result    covers
cases discussed in Theorems \ref{Theorem 3} and \ref{Theorem 5}. More
importantly, it  generalizes Theorem \ref{Theorem 4}, by replacing in the
equation (\ref{Prokaj}) (respectively, in the equation (\ref{eq: prokaj 2})) 
both driving local martingales $\,M(\cdot)\,$ and $\,
N(\cdot)\,$  by semimartingales. This generalization can be construed as an
analogue  of the results in 
\cite{MR0336813} %
and %
\cite{MR532447}%
. With bounded, measurable $\, f : \R^m
\rightarrow \R^m\,$, these authors show (for $\,m=1\,$ and for general $\,m
\in \N\,$, respectively) 
that, even in situations in which the ordinary differential equation $\,
\d  Y (t) = f (Y(t))\, \d  t\,$ might not be solvable, the
addition of a Brownian perturbation $\, W(\cdot)\,$ as in $$\, \d  Y
(t)\, =\, f \big(Y(t)\big)\, \d  t +  \d  W(t)\,$$ restores to
the differential equation a pathwise unique, strong solution.    

\medskip

\begin{proposition}
  \label{Prop}
  With the same assumptions and notation as in Theorem \ref{Theorem 4}, and with $\, \Gamma (\cdot)\,$ a continuous, adapted process of finite first variation on compact intervals, pathwise uniqueness holds for the stochastic differential equation 
  \begin{equation} 
    \label{eq: prokaj 3} 
    Y(t)\, =\,y+  \int_0^t f \big(Y(s)\big)\, {\mathrm d} \big( M(s) + \Gamma (s) \big) + A(t) +    N(t) \, , \quad 0 \le t < \infty \,,
  \end{equation}
  provided that either 
  \begin{rlist}
  \item \label{prop:8.1:it:1}
    the process $\, \Gamma (\cdot)\,$ is increasing, and the function $\, f(\cdot)\,$ is decreasing;  or
  \item \label{prop:8.1:it:2}
    the process $\, \Gamma (\cdot)\,$ is decreasing, and the function $\, f(\cdot)\,$ is increasing;  or
  \item \label{prop:8.1:it:3}
    the process $\, \Gamma (\cdot)\,$ is of the form $$\, \,
    \Gamma (\cdot)  =   \int_0^\cdot \varphi (t) \,    \d  \langle M  \rangle (t)  \, $$ for some progressively measurable process $\, \varphi (\cdot)\,$    which   satisfies the integrability condition $ \, \, \int_0^T   | \varphi (t)|^2 \, \d  \langle M  \rangle (t) < \infty\,$,    $\,~\forall~ ~T \in (0, \infty)\,$.
  \end{rlist}
\end{proposition}
\begin{proof} We show exactly as before that the continuous   semimartingale $\, D(\cdot) = X(\cdot) - Y(\cdot)\,$ accumulates zero local time $\,L  (\cdot) \equiv  \Lambda^D (\cdot\,, 0) \equiv 0\,$ at the origin, so we have 
  \begin{align}
    \label{|D|}
    | D (\cdot)| &=  \int_0^\cdot  \eta (t) \,\big(  \d  M(t) + \d  \Gamma (t)
    \big), 
    \intertext{with}
    \nonumber
    \eta (t) &:= \sign \big(  X(t) - Y(t)\big) \big( f (X(t) ) - f (Y(t)
    )\big)\,. 
  \end{align}
  Under either of the conditions {\it(\ref{prop:8.1:it:1})} or
  {\it(\ref{prop:8.1:it:2})}, the process  $\, |D  (\cdot) |\,$  is now a
  continuous   local supermartingale thanks to the representation
  (\ref{|D|}), thus a true supermartingale (by \textsc{Fatou}'s lemma) since
  it is nonnegative. As we 
  have $\, D(0) =0\,$, we conclude $\, D(\cdot) \equiv 0\,$ just as before.   

  \smallskip
  Under the condition {\it(\ref{prop:8.1:it:3})}, we note first that the
  process $\,\varphi (\cdot)\,$  satisfies also the integrability
  condition $\,\int_0^T|\varphi (t)|\,\d\q{M}(t)<\infty\,$ for every     
  $\,T\in (0, \infty)\,$,  by the  \textsc{Kunita-Watanabe} inequality;
  cf.$\,$%
  \cite{MR1121940}%
  ,   p.$\,$142. We introduce   the
  continuous local martingale  
  \begin{equation}
    \label{||D||}
    K (\cdot) := \int_0^\cdot \varphi (t) \, \d    M   (t)\,,
    \qquad \hbox{note} \qquad 
    \Gamma (\cdot) \,=\,  \langle M, K  \rangle (\cdot)\,,
  \end{equation}
  and  summon  the ``stochastic  exponential"
  $$
  \mathfrak{Z} (\cdot)\,:=\, 
  \exp\cbr{ - K ( \cdot) - \frac12 \, \q{K} (\cdot)},
  $$
  the unique solution of the ``simplest stochastic  integral equation" (in the
  terminology of %
  \cite{MR0247684}%
  ) $\,\mathfrak{Z}(\cdot) =1 - \int_0^\cdot\mathfrak{Z}(t) \,\d  K(t)\,$.   

  \smallskip
  We want to show $\, D(\cdot) \equiv 0$.    Observe that $\,\Pi (\cdot)
  :=\mathfrak{Z} (\cdot) \abs{D} (\cdot)\,$ is a local
  martingale by the product rule of the stochastic calculus: 
  \begin{align*}
    \Pi(\cdot)\,&=\,\int_0^\cdot \mathfrak{Z} ( t)\,  \d |D(t)| + \int_0^\cdot|D(t)|\, \d  \mathfrak{Z}(t)
    - \int_0^\cdot\eta (t) \, \mathfrak{Z} (t)\, \d \q{M, K} (t)\\
    \,&= \, \int_0^\cdot\mathfrak{Z} (t)\eta (t)\,\d M(t) + \int_0^\cdot|D( t)| \,\d\mathfrak{Z}(t)\,.
  \end{align*}
  Here $M(\cdot) $ and $\mathfrak{Z}(\cdot)$ are continuous local martingales and the
  integrands are locally bounded, hence both terms on the right are continuous local
  martingales. Now, $\Pi (\cdot)$ is a nonnegative,  continuous local martingale, starting
  from the origin at time zero, so it  must stay at the origin at all times;
  and since $\, \mathfrak{Z} (\cdot)\,$ is strictly positive, this leads to
  the desired conclusion $\, D(\cdot) \equiv 0\,$.   
\end{proof}

\subsection{  \textsc{Maximality}} %
\label{Max}

Let us place ourselves in the non-degenerate case $\, \rho\, \sigma >0\,$ and in the constructive setup (synthesis) of section \ref{sec3}, where we work with the filtration $\,\mathbf{ F}^{\, (W_1, W_2)}\,$ generated by the planar Brownian motion $\, (W_1 (\cdot), W_2 (\cdot))\,$. 

In  the terminology of \cite{MR2483736}, we shall say that a one-dimensional
Brownian motion $\, \betab (\cdot)\,$ in the  filtration $\,\mathbf{ F}^{\, (W_1, W_2)}\,$, is
\begin{rlist}
\item\emph{complenentable in} $\,\mathbf{ F}^{\, (W_1,W_2)}\,$, if there
  exists an 
independent one-dimensional Brownian motion  $\, \etab (\cdot)\,$ %
 in the filtration   $\,\mathbf{F}^{\, (W_1, W_2)}\,$, and for which
$\,\mathbf{ F}^{\,(\betab, \etab)}=\mathbf{ F}^{\, (W_1, W_2)}\,$;
\item\emph{maximal in} $\,\mathbf{ F}^{\, (W_1, W_2)}\,$, if
  for any one-dimensional Brownian motion $\, \thetab (\cdot)\,$  %
 in the filtration    $\,\mathbf{ F}^{\, (W_1, W_2)}\,$ and for which 
  $\,\mathfrak{ F}^{\, \betab} (t)  \subseteq \mathfrak{ F}^{\, \thetab} (t) \,$
  holds for all  $ \, 0 \le t < \infty\,$, we have $\, \mathfrak{ F}^{\,
    \thetab} (t) = \mathfrak{ F}^{\, \betab} (t)  \,$, $ \, 0 \le t <
  \infty\,$. 
\end{rlist}
\noindent
\cite{MR2483736} show that complementability implies maximality; it is still
an open question whether the reverse is true.  

In section \ref{sec3} we start with the  planar Brownian motion $\, (W_1
(\cdot), W_2 (\cdot))\,$ and construct 12 one-dimensional Brownian motions in the filtration    $\,\mathbf{ F}^{\, (W_1,  W_2)}\,$. We shall examine these properties for each of them.  

\smallskip
Let us start with the independent one-dimensional   Brownian motions $\, B_1
(\cdot)\,$ of (\ref{3.12}) and $\, B_2 (\cdot)\,$ of (\ref{3.13}). In view of
(\ref{4.12}), (\ref{4.12.a}), Theorem \ref{Theorem 3} and Proposition
\ref{Proposition 2}, we have then $\, \mathfrak{ F}^{\, (B_1, B_2)} (t)=
\mathfrak{ F}^{\, (X_1, X_2)} (t) =  \mathfrak{ F}^{\, (W_1, W_2)} (t)\,$, $
\, 0 \le t < \infty\,$.  Thus, each of   $\, B_1 (\cdot)\,$ of (\ref{3.12})
and $\, B_2 (\cdot)\,$   is  complementable (by the other one), therefore also
maximal, in the filtration $\,\mathbf{ F}^{\, (W_1, W_2)}\,$.

It follows also fairly directly  from   (\ref{3.5.a}) and (\ref{3.5.ab}), that
the independent Brownian motions $\,W(\cdot)\,$ and $\,U^\flat (\cdot)\,$
complement each other in $\,\mathbf{ F}^{\, (W_1, W_2)}\,$; the same is true
of $\,U(\cdot)\,$ and $\,W^\flat (\cdot)\,$. Thus, all these   one-dimensional
Brownian motions are maximal  in   $\,\mathbf{ F}^{\, (W_1, W_2)}\,$. 

On the other hand, let us consider the one-dimensional Brownian motion $\,
V^\flat (\cdot)\,$ of (\ref{3.5.bb}). This process is adapted to the
filtration $\,\mathbf{ F}^{\,  W }\,$ generated by the Brownian motion $\,
W(\cdot)\,$ of (\ref{3.5.a}); this is a consequence of the representation  
$\, V^{\,\flat} (\cdot) \,=\, \int_0^\cdot\,\sign \big( Y(t) \big)\,\d  W
(t)\,$ in (\ref{2.18}), and of the strong solvability of the stochastic
equation (\ref{1.3}). But we also have $\, \mathfrak{ F}^{\,V^\flat} (t)\,=\,
\mathfrak{ F}^{\,|Y|} (t) \,  \subsetneqq \,\mathfrak{F}^Y(t) =
\mathfrak{F}^W(t)  \,$,  
$\, \,0 < t < \infty\,$ from (\ref{4.3.a}), so $\, V^\flat (\cdot)\,$ cannot
possibly be maximal in the  two-dimensional filtration $\,\mathbf{ F}^{\,
  (W_1, W_2)}\,$. 
  
As a corollary of these considerations, we observe  also that a linear
combination of independent Brownian motions which are maximal, such as $\,
V^\flat (\cdot) = \rho\, V_1 (\cdot) - \sigma \, V_2 (\cdot)\,$ in
(\ref{3.5.bb}), can fail to be maximal. 
  
\smallskip
As for the Brownian motion $\, Q^\flat (\cdot)\,$ of (\ref{3.5.b}), it can
fail  to be maximal. Indeed, it is clear from (\ref{3.5.b}), (\ref{3.5.bb})
and the observations in the previous paragraph   that    in the isotropic case
$\, \rho = \sigma\,$  this Brownian motion $\,   Q^{\,\flat} (\cdot) \equiv
V^\flat (\cdot)\,$ is  not maximal in   $\,\mathbf{ F}^{\, (W_1, W_2)}\,$.  It
turns out that  $\, Q^\flat (\cdot)\,$ {\it is}    maximal in $\,\mathbf{
  F}^{\,    (W_1, W_2)}\,$, however, if   $\, \rho \neq \sigma\,$.    To see
this, recall the notation $\, \delta=2\rho\sigma\,$ and $\,
\gamma=\rho^2-\sigma^2\,$ from (\ref{2.19}), and observe   
\begin{multline*}
  \delta \int_0^t \sign(Y(s))\, \mathrm{d}  Q^\flat(s) +\gamma \,W^\flat (t)\,=\,
  (\delta \sigma+\gamma \rho )\,W_1(t) +(\delta \rho-\gamma \sigma)\,W_2 (t)\,=\\=\,
  \rho W_1 (t)+\sigma W_2(t)\,=\,W(t)\, =\,Y (t) -y-\lambda \int_0^t \sign(Y(s))\,\mathrm{d}s\,;
\end{multline*}
equivalently, 
\begin{displaymath}
  Y(t)\,=\,y + \int_0^t \sign(Y(s))\,\mathrm{d}\big(\delta Q^\flat(s)+\lambda s\big)+\gamma \, W^\flat (t)\,, \qquad 0 \le t < \infty\,.
\end{displaymath}
Then, since we are assuming $\, \delta=2\rho\sigma >0\,$ and   $\,
\gamma=\rho^2-\sigma^2 \neq 0\,$,  it follows from Proposition
\ref{Prop}~(\ref{prop:8.1:it:3}) that pathwise uniqueness holds for the above
equation; and since the equation admits a weak solution, this solution is
actually strong. In particular, the process $\, Y(\cdot)\, $ is adapted to 
$\,\mathbf{F}^{(Q^\flat,W^\flat)}\,$,  and we arrive at the filtration identities 
$\, \mathbf{F}^{(Q^\flat,W^\flat)}=\mathbf{F}^{(Q^\flat,W^\flat,Y)}
=\mathbf{F}^{(Q^\flat,W^\flat,W)}=\mathbf{F}^{(W_1,W_2)}\, $. Thus $\,W^\flat
(\cdot)\, $ complements the Brownian motion $\,Q^\flat (\cdot)\,$ in
$\,\mathbf{ F}^{\,    (W_1, W_2)}\,$. 

\smallskip
Recall now the one-dimensional Brownian motion $\, Q (\cdot)\,$ of
(\ref{3.5.bb}), which is independent of $\, V^\flat (\cdot)\,$; it is also
independent of  the one-dimensional Brownian motion $\, W (\cdot)\,$, and
indeed we have $\, \mathfrak{ F}^{\, (W_1, W_2)} (t) =  \mathfrak{ F}^{\,
  (W, Q)} (t) \,$, $ \, 0 \le t < \infty\,$  from Proposition
\ref{Proposition 2}. In other words, the one-dimensional Brownian motion $\,
Q (\cdot)\,$ is complementable (by $W(\cdot))$, therefore also maximal, in
the two-dimensional filtration $\,\mathbf{ F}^{\, (W_1, W_2)}\,$. 

\medskip
Let us consider now the one-dimensional Brownian motions $\, V_1 (\cdot)\,$
and $\, V_2 (\cdot)\,$ of (\ref{3.5}); we claim that they are complementable
(by $W_2(\cdot)$ and $W_1(\cdot)$, respectively;  though {\it not by one
  another!}$\,$), therefore also maximal, in the two-dimensional filtration
$\,\mathbf{ F}^{\, (W_1, W_2)}\,$. To see the first claim (the second one is
argued in a completely analogous manner), let us observe that we have the
equalities $\,\d Y (t) + \lambda \mathrm{sgn} \big( Y(t) \big) \, \d t = \d W
(t)= \rho\, \mathrm{sgn} \big( Y(t) \big) \, \d V_1 (t) + \sigma \, \d W_2
(t)\,$  
by virtue of (\ref{3.5.a}), (\ref{3.5}), (\ref{1.3}), therefore
$$
\d Y (t)  
\,=\,  \mathrm{sgn} \big( Y(t) \big) \, \big( \,\rho\, \d V_1 (t) - \lambda  \, \d t \big) + \sigma \, \d W_2 (t)\,.
$$
But this is an equation of the form (\ref{eq: prokaj 3}) with $\, A(\cdot)
\equiv 0\,$, for which Proposition \ref{Prop}~(\ref{prop:8.1:it:3}) is
satisfied, and the Brownian motions $\, V_1(\cdot)\,$ and $\, W_2(\cdot)\,$
are independent. Thus, pathwise uniqueness holds for this equation and, since
the equation admits a weak  solution, this solution is actually strong: 
\[
\mathfrak{ F}^{\, Y} (t)  \, \subseteq \,  \mathfrak{ F}^{\, (V_1, W_2)} (t) \,, \qquad  \, 0 \le t < \infty\,.
\]
Furthermore, we deduce   $\, \mathfrak{ F}^{\, V_2} (t)    \subseteq
\mathfrak{ F}^{\, (V_1, W_2)} (t)\,$ from the second equation in (\ref{3.5}),
therefore also  
\[
\mathfrak{ F}^{\, V} (t)  \, \subseteq \,  \mathfrak{ F}^{\, (V_1, W_2)} (t) \,, \qquad  \, 0 \le t < \infty\,;
\]
it follows from the last two displayed inclusions and (\ref{4.11}) that we
have $$\, \mathfrak{ F}^{\, (W_1, W_2)} (t)  = \mathfrak{ F}^{\, (Y, V)} (t)
\subseteq    \mathfrak{ F}^{\, (V_1, W_2)} (t) \subseteq  \mathfrak{ F}^{\,
  (W_1, W_2)} (t) \,, \qquad  0 \le t < \infty\,,$$ and this proves the
complementability of $\, V_1 (\cdot)\,$ by $\, W_2(\cdot)\,$. As for the third
claim of this paragraph,  %
it is fairly straightforward from the strict inclusion in (\ref{4.11}) that
$\, V_1 (\cdot)\,$ cannot be complemented by $\, V_2(\cdot)\,$ in the
two-dimensional filtration $\,\mathbf{ F}^{\, (W_1, W_2)}\,$.

\medskip
Finally, let us turn to  the Brownian motion $\, V  (\cdot) = \rho\, V_1 (\cdot) + \sigma \, V_2 (\cdot)\,$ in (\ref{3.5.b}). We follow the same approach as for the Brownian motion  $\,Q^\flat (\cdot)\,$, only with the roles of
$\, \delta=2\,\rho\sigma$ and $\gamma=\rho^2-\sigma^2$   interchanged, and with 
$W^\flat (\cdot)$   replaced by $U (\cdot) $ from \eqref{3.5.ab}. A bit more
precisely, we observe   
\begin{multline*}
  \gamma \int_0^t \sign(Y(s))\, \mathrm{d}  V(s) +\delta \,U^\flat (t)\,=\,
 \gamma \,( \rho  \,W_1(t) - \sigma W_2 (t))    + \delta ( \sigma  \,W_1(t) - \rho \,W_2 (t))  \,=\\=\,
  \rho W_1 (t)+\sigma W_2(t)\,=\,W(t)\, =\,Y (t) -y-\lambda \int_0^t \sign(Y(s))\,\mathrm{d}s\,,
\end{multline*}
or equivalently 
\begin{displaymath}
  Y(t)\,=\,y + \int_0^t \sign(Y(s))\,\mathrm{d}\big(\gamma \,V(s)+\lambda s\big)+\delta \, U^\flat (t)\,, \qquad 0 \le t < \infty\,.
\end{displaymath}
Since we are assuming $\, \delta=2\rho\sigma >0\,$,  it follows from
Proposition \ref{Prop}~(\ref{prop:8.1:it:3}) that pathwise uniqueness holds for
the above equation (even when $\gamma=\rho^2-\sigma^2=0)$; and since the
equation admits a weak solution, this solution is actually strong. To wit, the
process $\, Y(\cdot)\, $ is adapted to $\,\mathbf{F}^{(V,U^\flat)}\,$,  and we
arrive at the filtration identities $\,
\mathbf{F}^{(V,U^\flat)}=\mathbf{F}^{(Y, V,U^\flat )} =\mathbf{F}^{(W,
  V,U^\flat )}=\mathbf{F}^{(W_1,W_2)}\, $ (recall Proposition \ref{Proposition
  2}). Thus $\,U^\flat (\cdot)\, $ complements the Brownian motion $\,V
(\cdot)\,$ in $\,\mathbf{ F}^{\,    (W_1, W_2)}\,$. 

\begin{acknowledgements}

Numerous   discussions with Drs.$\,$Vasileios   \textsc{Papathanakos} and
Mykhaylo \textsc{Shkolnikov}, as well as several  exchanges  with Professor
Michel \textsc{Emery}, are very gratefully   acknowledged.  Drs.$\,$Daniel
\textsc{Fernholz} and Johannes \textsc{Ruf} 
provided crucial motivation and encouragement, when it came to  settling
issues of strength and/or uniqueness as in Theorems \ref{Theorem 1},
\ref{Theorem 2} and \ref{Theorem 3}. We are indebted to Drs.$\,$Vasileios
\textsc{Papathanakos} and Adrian \textsc{Banner} for suggesting the
representation of SO(2)  inherent in (\ref{SO2})  and for urging us to
undertake the full analysis that appears now in subsection \ref{sec9},  as
opposed to our initial effort that amounted to Example \ref{51}.  
We are grateful to
Dr.$\,$Johannes \textsc{Ruf} for his very careful reading of the near final
version of the paper, and his detailed and incisive comments. 
Last, but  not least, we are grateful to Professor Marc \textsc{Yor} for
posing the question that led to our work in section \ref{App}.  

The research of the third author was  supported in part by the National
Science Foundation under grant NSF-DMS-09-05754. The research of the fourth
author is supported by the European Union and co-financed by the European
Social Fund (grant agreement no. TAMOP 4.2.1./B-09/1/KMR-2010-0003).
\end{acknowledgements}

\def\doi#1{doi:~\href{http://dx.doi.org/#1}{\nolinkurl{#1}}}
\def\eprint#1{arXiv:~\href{http://arxiv.org/abs/#1}{\nolinkurl{#1}}}

\end{document}